\documentclass[reqno]{amsart}

\usepackage{latexsym}
\usepackage{graphicx}

\newcommand{\poly}{{\operatorname{Poly}}}
\newcommand{\I}{{\mathbf I}}
\newcommand{\todd}{{\operatorname{Todd}}}
\newcommand{\T}{{\mathbf T}^m}

\newcommand{\sm}{\setminus}
\newcommand{\szego}{Szeg\"o }

\newcommand{\Si}{\Sigma}
\newcommand{\inv}{^{-1}}
\newcommand{\kahler}{K\"ahler }

\newcommand{\wt}{\widetilde}
\newcommand{\wh}{\widehat}
\newcommand{\PP}{{\mathbb P}}
\newcommand{\N}{{\mathbb N}}
\newcommand{\R}{{\mathbb R}}
\newcommand{\C}{{\mathbb C}}
\newcommand{\Q}{{\mathbb Q}}
\newcommand{\Z}{{\mathbb Z}}

\newcommand{\CP}{\C\PP}
\renewcommand{\d}{\partial}
\newcommand{\dbar}{\bar\partial}
\newcommand{\ddbar}{\partial\dbar}

\newcommand{\E}{{\mathbf E}}

\newcommand{\Log}{{\operatorname{Log\,}}}
\newcommand{\half}{{\frac{1}{2}}}
\newcommand{\vol}{{\operatorname{Vol}}}

\newcommand{\codim}{{\operatorname{codim\,}}}
\newcommand{\SU}{{\operatorname{SU}}}
\newcommand{\FS}{{{\operatorname{FS}}}}
\newcommand{\supp}{{\operatorname{Supp\,}}}
\renewcommand{\phi}{\varphi}
\newcommand{\eqd}{\buildrel {\operatorname{def}}\over =}
\newcommand{\go}{\mathfrak}
\newcommand{\acal}{\mathcal{A}}

\newcommand{\ccal}{\mathcal{C}}
\newcommand{\dcal}{\mathcal{D}}
\newcommand{\ecal}{\mathcal{E}}
\newcommand{\fcal}{\mathcal{F}}

\newcommand{\hcal}{\mathcal{H}}
\newcommand{\ical}{\mathcal{I}}
\newcommand{\jcal}{\mathcal{J}}

\newcommand{\lcal}{\mathcal{L}}
\newcommand{\mcal}{\mathcal{M}}
\newcommand{\ncal}{\mathcal{N}}
\newcommand{\ocal}{\mathcal{O}}

\newcommand{\rcal}{\mathcal{R}}
\newcommand{\scal}{\mathcal{S}}

\newcommand{\al}{\alpha}
\newcommand{\be}{\beta}
\newcommand{\ga}{\gamma}
\newcommand{\Ga}{\Gamma}
\newcommand{\La}{\Lambda}
\newcommand{\la}{\lambda}
\newcommand{\ep}{\varepsilon}

\newcommand{\De}{\Delta}
\newcommand{\om}{\omega}

\newcommand{\di}{\displaystyle}
\newtheorem{theo}{{\sc Theorem}}[section]
\newtheorem{cor}[theo]{{\sc Corollary}}
\newtheorem{lem}[theo]{{\sc Lemma}}
\newtheorem{prop}[theo]{{\sc Proposition}}

\newenvironment{rem}{\medskip\noindent{\it Remark:\/} }{\medskip}
\newenvironment{defin}{\medskip\noindent{\it Definition:\/} }{\medskip}

\title{Random
polynomials with prescribed Newton polytope}

\author{Bernard Shiffman}
\author{Steve Zelditch}
\address{Department of Mathematics, Johns Hopkins University, Baltimore, MD
21218, USA} \email{shiffman@math.jhu.edu, szelditch@jhu.edu}

\thanks{Research partially supported by NSF grant
DMS-0100474 (first author) and DMS-0071358 (second author).}

\date{May 17, 2003}

\begin{document}

\begin{abstract} The Newton polytope $P_f$ of a polynomial $f$ is
well known to have a strong impact on its behavior. The
Bernstein-Kouchnirenko theorem asserts that even the number of
simultaneous zeros in $(\C^*)^m$ of a system of $m$ polynomials
depends on their Newton polytopes. In this article, we show that
 Newton polytopes further  have a strong impact on the
distribution of  zeros and pointwise norms of polynomials, the basic theme
being that Newton polytopes determine allowed and forbidden
regions in $(\C^*)^m$ for these distributions.

Our results are statistical and asymptotic in the degree of the
polynomials. We equip the space of polynomials of
degree $\leq p$ in $m$ complex variables with its usual SU$(m +
1)$-invariant Gaussian probability measure and then consider the
conditional measure induced on the subspace of polynomials with
fixed Newton polytope $P$. We then determine the asymptotics of
the conditional expectation $\E_{|N P}(Z_{f_1, \dots, f_k})$ of
simultaneous zeros of $k$ polynomials with Newton polytope $NP$ as
$N \to \infty$. When $P = \Sigma$, the unit simplex, it is
clear that the expected zero  distributions
$\E_{|N\Sigma}(Z_{f_1, \dots, f_k})$ are uniform relative to the
Fubini-Study form.  For a convex polytope $P\subset p\Sigma$, we
show that there is an {\it allowed region\/} on which
$N^{-k}\E_{|N P}(Z_{f_1, \dots, f_k})$ is asymptotically uniform
as the scaling factor $N\to\infty$. However, the zeros have an
exotic distribution in the complementary {\it forbidden region\/}
and when $k = m$ (the case of the Bernstein-Kouchnirenko theorem),
the expected percentage of simultaneous zeros in the forbidden
region approaches 0 as $N\to\infty$.

\end{abstract}

\maketitle

\tableofcontents

\section{Introduction}

It is well known that the Newton polytope $P$ of a
polynomial $f(z_1, \dots, z_m)$ of degree $p$ has a crucial
influence on its value distribution and in particular on its zero
set.  Even the number of simultaneous zeros   in
$(\C^*)^m:=(\C\sm\{0\})^m$ of $m$ generic polynomials $f_1, \dots,
f_m$ depends on their Newton polytopes  \cite{Be, Ku1, Ku2}. Our
purpose in this paper  is to demonstrate that the Newton polytope
of a polynomial $f$ also has a crucial influence on its {\it mass
density} $|f(z)|^2 dV$ and on the {\it spatial distribution} of
{\it zeros} $\{f = 0\}$. We will show that there is a {\it
classically allowed region\/}  where the mass almost surely
concentrates  and a {\it classically forbidden region\/} where it
almost surely is exponentially decaying. The classically allowed
region is the inverse image $\mu^{-1}(\frac{1}{p}P^\circ)$ of the
(scaled, open) polytope $\frac{1}{p}P^\circ$ under the standard
moment map $\mu$ of complex projective space $\CP^m$. The simultaneous zeros in
$(\C^*)^m$ of $m$ generic polynomials $f_1, \dots, f_m$ with Newton polytope
$P$ tend to concentrate (in the limit of high degrees) in the
classically allowed region, giving a kind of quantitative localized version
of the Bernstein-Kouchnirenko theorem \cite{Be, Ku1, Ku2}. Polynomials with a given
Newton polytope $P$ are often called {\it sparse} in the
literature, and methods of algebraic (including toric) geometry
have recently been applied to the computational problem of
locating zeros of systems of  such sparse polynomials (e.g., see
\cite{HS, MaR, S, R,V}). Our results give information on the
expected location of zeros when the polynomials are given the
conditional measure described below. To our knowledge, this kind of precise
asymptotic concentration of zeros of sparse systems of random polynomials has not been
observed before.

The Newton polytope has an equally strong (though different)
impact on the common zeros of $k < m$ polynomials $f_1, \dots, f_k$.  The image of
the zero set of $f_1, \dots, f_k$ under the moment map is (up to a
logarithmic re-parameterization) known as an {\it amoeba\/} in the
sense of \cite{GKZ, M}. Results on the expected distribution of
amoebas can be obtained from our results on the expected zero
current; for example, for a polytope $\frac{1}{p}P$ with vertices
in the interior of the standard unit simplex $\Sigma \subset \R^m$,
there is also a {\it very forbidden region} which the amoeba
almost surely avoids (see Corollary~\ref{subtler}).

\subsection{Statement of results}\label{s-results}
 The patterns
in zeros discussed here are statistical---they hold for random polynomials with
prescribed Newton polytope---and are asymptotic as the degree of
the polynomials tends to infinity.
      To state our problems and results precisely, let us recall
      some definitions. Let
\begin{equation}\label{chi}f(z_1, \dots, z_m) = \sum_{\alpha \in \N^m:
|\alpha| \leq p} a_{\alpha} \chi_\al(z_1, \dots, z_m),\;\;\;\;\
\chi_\al(z) = z_1^{\alpha_1} \cdots z_m^{\alpha_m}\end{equation}
($\al=(\al_1,\dots,\al_m),\ |\al|=\al_1+\cdots+\al_m$) be a polynomial of degree
$p$ in
      $m$ complex variables. By the {\it support} of $f$ we mean the set
    \begin{equation} \label{SUPPORT} S_f = \{\alpha\in\Z^m: a_{\alpha} \not
= 0\},\end{equation}
    and by its   {\it Newton polytope} $P_f$ we mean the set
\begin{equation}\label{NEWTONP}  P_f :  =
\mbox{the convex hull in $\R^m$ of } \;\; S_f. \end{equation} Our
aim is to study the statistical patterns in the space
\begin{equation}\label{polyP}\poly(P)=\{f\in\C[z_1,\dots,z_m]: P_f\subset P
\}\end{equation}
of polynomials with support contained in a
fixed Newton polytope $P$. We therefore need to find a natural
measure  on this space of polynomials. Since our purpose is to
compare zero sets and masses as the polytope $P$ varies, we  view
the polytope $P$ as placing a condition on the Gaussian ensemble
$\poly(p\Si)$ of all polynomials of  degree $\le p$ (where $p\Si$ is the dilation by $p$
of the standard unit simplex $\Si$), and we give
$\poly(P)$ the resulting conditional probability measure. Here, $p$ is chosen such that
$p\Si\supset P$.

To describe this Gaussian ensemble, we first identify $\poly(p\Si)$ with  the space of
homogeneous polynomials  of degree $p$ in $m+1$ variables  in the usual way, i.e., by
identifying
$f\in\poly(p\Si)$ with the homogeneous  polynomial $F$ such that
$F(1, z_1, \dots, z_m)=f(z_1, \dots, z_m)$.  Using this identification, we  give the
space $\poly(p\Si)$  the
$\SU(m+1)$-invariant  inner product
\begin{equation}\label{inner} \langle f_1, \bar f_2\rangle:=
\langle F_1, \bar F_2\rangle_{S^{2m+1}} = \frac 1{m!}\int_{S^{2m+1}} F_1 \bar
F_2\, d\nu\;,\end{equation} where $d\nu$ is Haar probability measure on the
$(2m+1)$-sphere
$S^{2m+1}\subset\C^{m+1}$. The monomials $\chi_\al$ are orthogonal, so an orthonormal
basis of
$\poly(p\Si)$ is given by
$\left\{\|\chi_{\alpha}\|\inv
\chi_{\alpha}\right\}_{|\al|\le p}$, where $\|\cdot\|$ denotes the
norm in $\poly(p\Si)$ given by (\ref{inner}).  We emphasize here that the norms
$\{\|\chi_{\alpha}\|\}$ of the monomials $\{\chi_\al\}$ do not depend on the
constraining polytope
$P$---they are given by the $\SU(m+1)$-invariant inner product (\ref{inner}) on
$\poly(p\Si)$.  (Of course they do depend on the choice of the integer $p$.)

The
$\SU(m+1)$-invariant  Gaussian measure $\gamma_p$ corresponding to the inner product
(\ref{inner}) is defined by
\begin{equation}
\label{G} d \gamma_p (f) = \frac{1}{\pi^{k_p}}e^{-|c|^2}
dc,\;\;\;\; f = \sum_{|\alpha|\le p} c_{\alpha}
\frac{\chi_{\alpha}}{\|\chi_{\alpha}\|}\;,\end{equation} where
$k_p=\#\{\al\in\N^m: |\al|\le p\}= {m+p\choose p}$. Thus, the
coefficients $c_{\alpha}$ are independent complex Gaussian random
variables with mean zero and variance one. We then endow the space
$\poly(P)$ with the associated  {\it conditional probability
measure\/} $\gamma_{p|P}$:
 \begin{equation} \label{CG} d \gamma_{p|P} (f) =
\frac{1}{\pi^{\# P}}e^{-|c|^2} dc,\quad f = \sum_{\alpha \in P}
c_{\alpha} \frac{\chi_{\alpha}}{\|\chi_{\alpha}\|}\;,
\end{equation} where the coefficients $c_{\alpha}$ are again
independent complex Gaussian random variables with mean zero and
variance one. (By an abuse of notation, we let $\sum_{\al\in
P}$ denote the sum over the lattice points $\al\in P\cap\Z^m$; $\#
P$ denotes the cardinality of $P\cap\Z^m$.) We observe that $\poly(P)$ inherits
the inner product $\langle f_1,\bar f_2 \rangle$ from  $\poly(p\Si)$,  and that
$\gamma_{|P}$ is the induced  Gaussian measure. Probabilities
relative to $\gamma|_P$  can be  considered as conditional
probabilities; i.e. for any event $E$,
$$\mbox{Prob}_{\gamma}\{ f \in E | P_f = P\} = \mbox{Prob}_{\gamma|_P} (E).$$

The expected distribution of mass and zeros of polynomials with
fixed Newton polytope  turns out to involve the moment map
$\mu:(\C^*)^m\to \R^m$ given by
\begin{equation}\label{momap}\mu(z)=\left(\frac{|z_1|^2}{1+\|z\|^2},\dots,
\frac{|z_m|^2}{1+\|z\|^2} \right)\;.\end{equation} The image of $\mu$ is the
interior of the standard unit simplex $\Si$ in $\R^m$ with
vertices at the points
$$(0,\dots,0),\ (1,0,\dots, 0),\ (0,1,\dots,0),\ \dots,(0,\dots,0,1)\ .$$ The map
$\mu$ is the moment map of $\CP^m\supset (\C^*)^m$ (with the
Fubini-Study symplectic form $\om_\FS$) and plays a role in the
geometric approach in \cite{SZ2, STZ1, STZ2}, where we regard polynomials
with a fixed polytope as sections of a holomorphic line bundle on
a toric variety (see \S\ref{Appendix}).

By a {\it convex integral polytope\/} $P$, we mean the convex hull
in $\R^m$ of a finite set of points in $\Z^m$. We use the moment
map $\mu$ to describe our classically allowed regions:

\begin{defin} Let $P$ be a convex integral polytope
in $ \R^m$ such that $P\subset p\Si$. The
 {\it classically allowed region\/} for polynomials in $\poly(P)$
is the set
\begin{equation}\label{AP}\acal_P:=\mu\inv\left(\frac{1}{p}P^\circ\right) \subset
(\C^*)^m\end{equation} (where $P^\circ$ denotes the interior of $P$), and the {\it
classically forbidden region\/} is its complement
$(\C^*)^m\sm\acal_P$. (If $P$ has empty interior, i.e. if $\dim P
<m$, then $\acal_P=\emptyset$.)
\end{defin}

All of our  results on the asymptotic expected distribution of
zeros  of systems of  random
polynomials with prescribed Newton polytopes follow from our asymptotic result
(Theorem~\ref{MASS}) on the pointwise expected values of random polynomials with
given Newton polytopes.  Our results on zeros have a somewhat different flavor
as the codimension (or number of polynomials in the system)
 varies. We first describe the results for the maximum codimension $m$ where
the zero set is almost surely discrete (Theorem~\ref{probK}). Then
after describing the expected value distribution, we discuss  the
codimension 1 case, where the zero set is a hypersurface, and we
finish with a result for general codimension that encompasses the
previous results.

\subsubsection{Distribution of zeros: the point case.}
Let us first consider  the simultaneous  zero set of $m$
independent random polynomials in $m$ variables. B\'ezout's theorem tells
us that $m$ generic homogeneous polynomials $F_1,\dots,F_m$ of
degree $p$ have exactly $p^m$ simultaneous zeros in $(\C^*)^m$. In
fact, one immediately sees (by uniqueness of Haar probability
measure on $\CP^m$) that the expected distribution of zeros is uniform with
respect to the $\SU(m+1)$-invariant projective volume form, when the
ensemble is given the $\SU(m+1)$-invariant measure $d\ga_p$.

According to the Bernstein-Kouchnirenko theorem \cite{Be, Ku1,Ku2}, the
number of common zeros in $(\C^*)^m$ of $m$ generic polynomials
$\{f_1,\dots, f_m\}$ with given Newton polytope $P$ equals $ m!
\vol(P)$, where $\vol(P)$ is the Euclidean volume of $P$. For
example, if $P=p\Si$, where $\Si$ is the standard unit simplex in
$\R^m$, then $\vol (p\Si)=p^m\vol(\Si)=\frac{p^m}{m!}$, and we get
B\'ezout's theorem.  (More generally, the $f_j$ may have different
Newton polytopes, in which case the number of zeros is given by
the Bernstein-Kouchnirenko formula as a `mixed volume' \cite{Be, Ku1,Ku2}.)

Now consider  the ensemble of $m$ independent random polynomials
with Newton polytope $P$, equipped with  the conditional
probability (\ref{CG}). We let $\E_{|P} (Z_{f_1, \dots, f_m})$
denote the expected density of their simultaneous zeros. More generally, for $N\in\Z^+$
we define the expected zero density measure
$\E_{|NP} (Z_{f_1, \dots, f_m})$ by
\begin{equation}\label{EZ0}\E_{|NP} (Z_{f_1, \dots, f_m})(U) =
\int d\ga_{Np|NP}(f_1)\cdots\int d\ga_{Np|NP}(f_m)\;\big[\#\{z\in U:
f_1(z)=\cdots=f_m(z)=0\}\big]\;,\end{equation} for  $U\subset
(\C^*)^m$, where the integrals are over $\poly(NP)$.  In fact,
$\E_{|NP} (Z_{f_1, \dots, f_m})$ is an absolutely continuous
measure given by a $\ccal^\infty$ density (see
Proposition~\ref{E-cond}).

Our first result shows that, as the polytope
$P$ expands, these zeros
 are expected to concentrate in the classically
allowed region and have (asymptotically) uniform density there. We measure the density
of zeros with respect to the projective volume $d\vol_{\CP^m}= \frac 1 {m!} \om_\FS^m$,
where
$\om_\FS=\frac i{2\pi}\ddbar
\log (1+\|z\|^2)$ is the  Fubini-Study \kahler form on $\C^m\subset\CP^m$.
(Note that with this normalization, the volume of $\CP^m$ is $\frac 1 {m!}$.)

\begin{theo}\label{probK} Suppose that $P\subset p\Si\subset\R^m$
is a convex integral polytope with nonempty interior. Then
$$\lim_{N\to\infty} \frac{1}{(Np)^m}\,\E_{|NP} (Z_{f_1, \dots, f_m})=
\left\{\begin{array}{ll}\om_\FS^m
\ \ & \mbox{\rm on \ }\acal_P\\[10pt]
0 & \mbox{\rm on \ }(\C^*)^m\sm \acal_P\end{array} \right.\ ,$$ in
the measure sense; i.e., for any Borel set $B\subset(\C^*)^m$, we have
$$\frac{1}{(Np)^m}\E_{|NP}\big(\#\{z\in B:
f_1(z)=\cdots=f_m(z)=0\}\big)\to m!\vol_{\CP^m}(B\cap
\acal_P)\;.$$
\end{theo}

Theorem~\ref{probK} is a special case of our general
result on zeros (Theorem~\ref{simultaneous}).
In fact, our results imply that the convergence of the zero
current on the classically allowed region is exponentially fast in
the sense that
$$\E_{|NP} (Z_{f_1, \dots, f_m})= (Np)^m \om_\FS^m  +O\left(e^{-\la N}\right)
\ \  \mbox{\rm on \ }\acal_P\;,$$ for some positive continuous
function $\la$ on $\acal_P$.

Note that we can choose any integer $p$ such that $p\Si \supset P$. Theorem~1 (as well
as our other results) depends quite strongly on the choice of $p$, since the
classically allowed region $\acal_P$ shrinks as $p$ is increased.  In order to minimize
notational clutter, $p$ does not appear $\acal_P,\ \E_{|NP},\
|f(z)|_\FS$, and some other expressions to be define below.

Figure \ref{allowed-square} shows the classically allowed region
(shaded) and the classically forbidden region (unshaded) when $P$
is the unit square in $\R^2$ (and $p=2$):

\begin{figure}[htb]
\centerline{\includegraphics*[bb= 1.9in 7.0in 4in 9.6in]{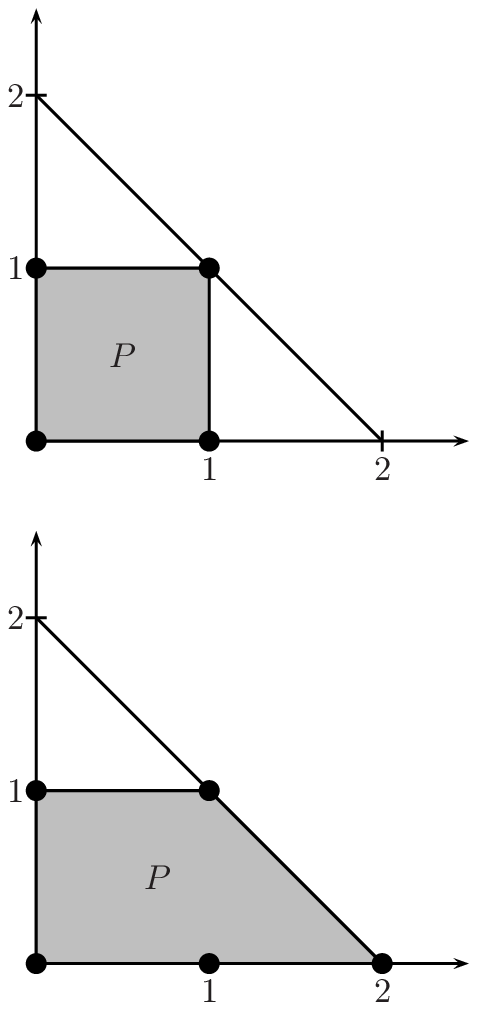}
\includegraphics*[bb= 1.6in 6.7in 5in 9.6in]{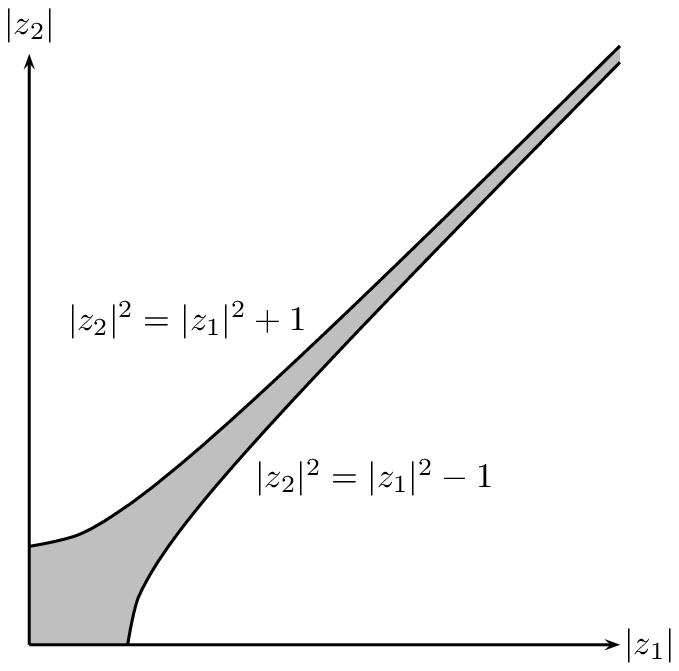}}
\caption{The classically allowed region for $P=[0,1]\times [0,1]$}
\label{allowed-square}\end{figure}

\subsubsection{Mass distribution.} By the mass density of a polynomial $f\in\poly(p\Si)$
at a point
$z\in(\C^*)^m$, we mean the square $|f(z)|_\FS^2$ of   the Fubini-Study norm
\begin{equation}\label{FSnorm}|f(z)|_\FS := \frac {|f(z)|}{(1+\|z\|^2)^{p/2}}
= |F(\zeta)|\;,
\qquad \zeta=\frac 1{\sqrt{1+\|z\|^2}}(1,z_1,\dots,z_m)\in S^{2m+1}\;,\end{equation}
where
$F$ is the homogenization of $f$ described above. Our next result describes the
asymptotics of the expected mass density of $\lcal^2$-normalized polynomials with
Newton polytope
$NP$, i.e.\ the expected density $\E_{\nu_{NP}}(|f(z)|^2_\FS)$
with respect to the Haar probability measure, denoted by
$\nu_{NP}$, on the unit sphere in $\poly(NP)$. We shall show  that
$\E_{\nu_{NP}}(|f(z)|^2_\FS)$ is
asymptotically uniform with respect to Fubini-Study measure in the
classically allowed region (as if there were no constraint at
all); while in the forbidden region the mass decays exponentially.
Thus in the semiclassical limit $N\to\infty$, all the mass
concentrates in the classically allowed region.

Behavior of the expected mass density (as well as the zero
distribution in positive dimension) in the forbidden region is
subtle, so we pause to introduce the relevant concepts.  We show
in \S \ref{s-normalbundle} that the classically  forbidden region
decomposes into the disjoint union of the normal `flow-outs' of
boundary points of $P$,
\begin{equation}\label{flow} \mbox{Flow}(x): = \{e^\tau \cdot z: \tau\in C_x ,\
p\mu(z) =x\}\;,\qquad x\in \d P\cap p\Si^\circ\;,
\end{equation}
where $C_x\subset \R^m$ is the normal cone to $P$ at $x$ (see \S
\ref{s-fans}). Here $r \cdot z = (r_1z_1,\dots,r_mz_m)$ denotes
the $\R_+^m$ action on $(\C^*)^m$.

We use these flow-outs to divide the forbidden region  into
subregions $\rcal_F$ as follows:  Recall that the boundary of a
polytope $P$ decomposes into a disjoint union of its faces (of all
dimensions). For each face $F$ of $P$ not contained in $\d(p\Si)$,
we let
\begin{equation} \rcal_F = \bigcup_{x\in F}\mbox{Flow}(x)
= \{e^\tau \cdot z: \tau\in C_F ,\ p\mu(z) \in F\}\;.
\label{RF}\end{equation}  where $C_F$ is the normal cone to $P$
along the face $F$; i.e., the normal cone to $P$ at any point of
$F$ (see \S \ref{s-fans}).  Note that for the case where $F$ is
the open face $P^\circ$,  $\rcal_F$ is the classically allowed
region $\acal_P$.  We call $\rcal_F$ the {\it flow-out\/} of the
face $F$. The regions $\rcal_F$ all have nonempty interior.  A
sample illustration of these regions if given in Figure~\ref{Fn}
in \S \ref{s-ex3}.

We can now state our result on mass asymptotics:

\begin{theo}\label{MASS}  Suppose that $P\subset p\Si\subset \R^m$
is a convex integral polytope such that $P\not\subset\d (p\Si)$.
Then the expected mass density of random $\lcal^2$-normalized
polynomials with Newton polytope $NP$ is given by the asymptotic
formulas:
$$
\begin{array}{rcll}\E_{\nu_{NP}}\left(|f(z)|^2_\FS\right)&
 \sim &
c_0 + c_1N^{-1} + c_2N^{-2} + \cdots,  \quad & \mbox{for }\ z \in
\acal_P,
\\[10pt] \E_{\nu_{NP}}\left(|f(z)|^2_\FS\right)&=& N^{-s/2} e^{-N b_P(z)}\left[c_0^F(z) +
O(N\inv)\right], &\mbox{for } z\in \rcal_F^\circ, \
F\subset\Si^\circ,
\end{array}  $$ where $c_0=\frac{p^m}{\vol(P)}$,  $s=\codim F$,
$c_0^F\in \ccal^\infty(\rcal_F^\circ)$, and $b_P$ is a positive
$\ccal^1$ function on $(\C^*)^m\sm\overline{\acal_P}$.
Furthermore, the remainder estimates are uniform on compact
subsets of the open regions $\rcal_F^\circ$ and of $\acal_P$.
\end{theo}

We shall give a formula for $b_P(z)$ below.  A more precise
asymptotic formula for the expected mass density is given by
Theorem~\ref{SZEGO} (see also (\ref{Eszego2})).

\begin{figure}[htb]
\centerline{\includegraphics*[bb= 164 295 447 530]
{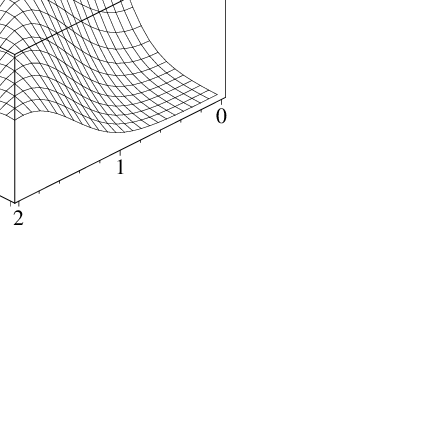}\hspace{-2cm}
\includegraphics*[bb= 164 295 447 530]{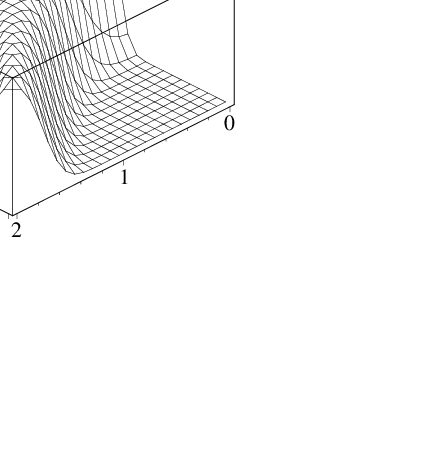}}
\caption{$\frac{1}{4}\E_{\nu_{NP}}\left(|f(z)|^2_\FS\right)$ for
$P=[0,1]\times [0,1]$} \label{mass-graph}\end{figure}

It follows that
$$ \E_{\nu_{NP}}\left(|f(z)|^2_\FS\right) \longrightarrow \left\{
\begin{array}{ll}
\frac{p^m}{\vol(P)}\,,  \quad  &\mbox{for }\ z
\in \acal_P\\[8pt] 0\,, & \mbox{for } z\in
(\C^*)^m\sm\overline{\acal_P}\end{array} \right.\ , $$ as
illustrated  in Figure \ref{mass-graph} below (plotted using
Maple) for the case where $P$ is the unit square. (Recall Figure
\ref{allowed-square} for the depiction of $\acal_P$ for this
case.)

Heuristically, the mass concentration in $\acal_P$ can be
understood as follows: the mass of a monomial $\chi_{\alpha}\in\poly(Np)$
concentrates (exponentially) on the torus $\mu^{-1}(\frac 1{Np}\alpha)$. The
constraint $S_f\subset NP$ thus concentrates all the mass in $\acal_P$, with
exponentially small errors, and the mass there is uniformly
distributed with small errors. Such mass concentration of
monomials is given in formula (\ref{Pidecay})  and in more detail
in \cite{STZ2}.

We also show that $\frac{1}{N}\log
\E_{|NP}\left(|f(z)|^2_\FS\right) \to - b_P(z)$ uniformly on
compact subsets of $(\C^*)^m$ (Proposition~\ref{convergence}),
even at the interfaces of the regions $\rcal_F$. Combined with the
Poincar\'e-Lelong formula, this uniform convergence result leads
to the asymptotics of the expected distribution of zeros of a
single random polynomial with polytope $NP$. Using independence of
sections and the Bedford-Taylor Theorem \cite{BT} on continuity of
the Monge-Amp\`ere operator, we obtain asymptotics for any number
$1 \leq k \leq m$ of sections.

We shall show that $b_P$  is  $\ccal^\infty$ on the regions
$\rcal_F^\circ$, but is not even $\ccal^2$ on the interfaces
between these regions. To state our formula for $b_P$, we let
$q(z)$ denote the unique point in $\d P$ such that
$z\in\mbox{Flow}(q(z))$ and we let $\tau_z\in C_{q(z)}$ be such
that $z= e^{\tau_z/2}\cdot \xi$ where $\xi\in\mu\inv \left(\frac
1p q(z)\right)$. Thus $q(z),\ \tau_z$ are given by the conditions:
\begin{eqnarray} &&
 p\mu(e^{-\tau_z/2}\cdot z) =q(z)\in\d P,\label{cond1}\\
&&\tau_z\in  C_{q(z)}.\label{cond2}\end{eqnarray} The existence
and uniqueness of $q(z)$ and $\tau_z$ are stated in Lemma
\ref{claim2-3}.

We then have:
\begin{equation}
\label{b} b_P(z)= -\langle q(z), \tau_z\rangle + p
\log\left(\frac{1+\|z\|^2}{1+\|e^{-\tau_z/2}\cdot z\|^2}\right)
\qquad \mbox{for }\ z\in(\C^*)^m\sm\acal_P\;.\end{equation}

To understand formula (\ref{b}) better, for any $x\in P$ (not
necessarily a lattice point),  we consider the `monomials'
\begin{equation}\label{mcal} |\chi_x(z)|:=|z|^x =
|z_1|^{x_1}\cdots|z_m|^{x_m}\;,\quad |\wh\chi_x(z)|:=
|\chi_x(z)|_\FS= \frac{|\chi_x(z)|}{(1+\|z\|^2)^{p/2}}\;,\quad
\mcal_x (z):=\frac {|\wh\chi_x (z)|} {\sup |\wh\chi_x|}\;,
\end{equation}
for ${z\in
(\C^*)^m}$.
Thus the normalized monomial $\mcal_x$ has sup-norm $1$; in
fact it takes its maximum on the torus $\mu\inv(\frac 1 p x)$.  One
easily checks that
\begin{equation} \label{interpretb} b_P(z)=-2\log
{\mcal_{q(z)}(z)}\;,\end{equation} and
thus Theorem~\ref{MASS} says that
\begin{equation}\label{interpretM}
\E_{\nu_{NP}}\left(|f(z)|^2_\FS\right)=\left[c_0^F(z) +
O(N\inv)\right]\, N^{-s/2}\, \mcal_{q(z)}(z)^{2N}\;, \qquad z\in
\rcal_F^\circ.\end{equation} We will also obtain an  integral
formula for $b_P$ (see Proposition \ref{b-prop}), which allows us
to interpret $b_P$ as an `Agmon distance' (see \S
\ref{s-remarks}), and $q(z)$ is the closest point of $P$ to
$p\mu(z)$ in this sense.

We do not assume in Theorem \ref{MASS} that $P$ has interior.
However, if $P^\circ=\emptyset$ (i.e., if $\dim P<m$), then we
must assume that  $P\not\subset \d (p\Si)$ since if $P$ is
contained in a face of $p\Si$, the decay function $b_P$ is not
defined. Note that if $P^\circ=\emptyset$, then $\acal_P$ is
vacuous and the first expansion of the theorem does not occur.

\subsubsection{Distribution of zeros: arbitrary codimension.}\label{s-distz}
Our next results concern the zero set of one  or more random
polynomials. Since   the zero set
\begin{equation}\label{zeroset}|Z_{f_1, \dots, f_k}|:=\{z\in(\C^*)^m:f_1(z)=\cdots
=f_k(z)=0\}\end{equation}
 is  a
submanifold of complex codimension $k$ (without multiplicity), for
almost all polynomials $f_1, \dots, f_k$, it defines a {\it
current of integration\/} $Z_{f_1, \dots, f_k}\in
\dcal'{}^{k,k}((\C^*)^m)$.  We recall that this current is given by
\begin{equation}\label{zerocurr}(Z_{f_1, \dots, f_k}, \phi) := \int_{|Z_{f_1, \dots,
f_k}|}\phi\;,\qquad
\mbox{for test forms }\ \phi\in \dcal^{m-k,m-k}((\C^*)^m)\;.\end{equation} The
zero set $|Z_{f_1, \dots, f_k}|$ also carries a natural
$(2m-2k)$-dimensional Riemannian volume measure (induced from the
Fubini-Study metric on $\CP^m$), denoted $\|Z_{f_1, \dots,
f_k}\|$, against which one can integrate scalar functions. Since
the volume form on any complex $n$-dimensional holomorphic
submanifold of $\CP^m$ is given by the restriction of $\frac 1
{n!}\om^n_\FS$, the volume measure on the zero set is given by
$\|Z_{f_1, \dots, f_k}\|= Z_{f_1, \dots, f_k} \wedge \frac
1{(m-k)!}\omega^{m-k}_\FS$; i.e.,
\begin{equation}\label{zerovol}(\|Z_{f_1, \dots, f_k}\|, \phi)
=\int_{|Z_{f_1, \dots, f_k}|}\phi\,d\vol_{2m-2k} =\frac 1{(m-k)!}
\int_{|Z_{f_1, \dots, f_k}|}\phi\; \omega^{m-k}_\FS  \;,\qquad
\mbox{for  } \phi  \in \dcal ((\C^*)^m)\;.\end{equation} As a measure,
$\|Z_{f_1, \dots, f_k}\|(U)=\vol(|Z_{f_1, \dots, f_k}|\cap U)$ for
open sets $U\subset(\C^*)^m$. (For the theory of currents defined by complex algebraic
or analytic varieties, see for example,
 \cite[I.3]{Shabat}.) We shall discuss the expected values of both the
current of integration
$Z_{f_1, \dots, f_k}$ and the measure $\|Z_{f_1, \dots, f_k}\|$.

We first consider the expected distribution of zeros of $1$
polynomial.  We denote by
 $\E_{|P }(Z_{f})=\E_{\ga_p|P }(Z_{f})$ the
conditional expectation of the zero current of a random polynomial
$f \in \poly(P)$  with Newton polytope $P$.  In fact, $\E_{|P
}(Z_{f})$ is actually a smooth $(1,1)$-form on $(\C^*)^m$
(Proposition \ref{E-cond}).

Let us recall what happens when $P=p\Si$. By the uniqueness of the
$\SU(m+1)$-invariant \kahler form $\om_\FS$ on $\CP^m$, the
expected zero current $\E (Z_{f})$ taken over all polynomials of
degree $p$ is given by $p\om_\FS$, where $\om_\FS=\frac{i}{2\pi}
\d\dbar \log\|z\|^2$ is the Fubini-Study \kahler form on $\CP^m$.
Thus the expected distribution of zeros, as well as the tangents
to the zero varieties, is uniform over $\CP^m$. We now describe
how the expectation changes if we add the condition that $P_f =
P$.

\begin{theo} \label{main}  Let $P\subset p\Si\subset \R^m$
be a convex integral polytope  such that $P\not\subset\d (p\Si)$.
Then there exists a closed semipositive $(1,1)$-form $\psi_P$ on
$(\C^*)^m$ with piecewise $\ccal^\infty$ coefficients such that:
\begin{enumerate}
\item[i)\ ]  $N^{-1}\E_{|NP} (Z_f)\to \psi_P$ \
in  $\lcal^1_{\rm{loc}}((\C^*)^m)$.
\item[ii)\ ] $\psi_P =p\om_\FS$ on the classically allowed region
$\mu\inv(\frac{1}{p}P^\circ)$.

\item[iii)\ ] On each region $\rcal_F^\circ$, the $(1,1)$-form
$\psi_P$ is $\ccal^\infty$ and has constant rank equal to $\dim
F$; in particular, if $v\in p\Si^\circ$ is a vertex of
$P$, then $\psi_P|_{\rcal_v^\circ}=0$.
\end{enumerate}
\end{theo}

We see from part (iii) that the zero set $|Z_f|$ of a polynomial with
polytope $NP$ typically  intersects the classically forbidden
region $\mu\inv(\Si\sm\frac{1}{p}P)$ in the semiclassical  limit
$N\to\infty$. However,  there are subtler `very forbidden
regions' that $|Z_f|$ avoids in the case where the polytope has
vertices in the interior of $p\Sigma$, namely the regions $\rcal_v^\circ$
comprising the flow-out of these vertices.

As a corollary, we obtain some statistical results on the
so-called `tentacles' of amoebas in dimension $2$ (see \S
\ref{AM}). Roughly speaking, the (compact) amoeba of a polynomial
$f(z_1,z_2)$ is the image of the Riemann surface $Z_f$ under the
moment map $\mu$ on $(\C^*)^2$, and the tentacles are the
ends of the amoeba. Certain tentacles must end at vertices of the
triangle $\Si$ while others are `free' to end anywhere along the
boundary of $\Si$. In Corollary \ref{amoeba}, we will prove that
(in the limit $N \to \infty$)  almost all of the free tentacles of
typical amoebas tend to end in the classically allowed portion of
$\d\Si$.

We call the form $\psi_P$ in Theorem \ref{main} the {\it limit
expected zero current.\/} Our explicit formula for $\psi_P$ is
\begin{equation}\label{psiP} \psi_P= p\om_\FS -
\frac{\sqrt{-1}}{2\pi} \ddbar b_P\;,\end{equation} where $b_P$ is
given by (\ref{b}). In fact, (\ref{psiP}) holds as an equation of
currents, and  the current $\psi_P\in \dcal'{}^{1,1}((\C^*)^m)$ is
closed and positive. By $\lcal^1_{\rm{loc}}$ convergence in (i),
we mean $\lcal^1$ convergence of the coefficients on every compact
subset of $(\C^*)^m$. (Recall that $\E_{|NP} (Z_{f_1, \dots,
f_k})$ is a $(k,k)$-form with smooth coefficients.) If we write
$\psi_P=\sqrt{-1}\sum \psi_{jk}(z)dz_j \wedge d\bar z_k$, then
$\big(\psi_{jk}(z)\big)$ is a semi-positive Hermitian matrix. By
the rank of $\psi_P$ at $z$, we mean the rank of the matrix
$\big(\psi_{jk}(z)\big)$. Note that if $\rcal_F$ and $\rcal_{F'}$
are adjoining regions (i.e., have a common codimension 1
interface), then $F$ and $F'$ are of different dimensions, so
$\psi_P$ must be discontinuous along the interface.

Boundary points of the $\rcal_F$ are called {\it transition
points\/}. The set of transition points comprises the
discontinuities of $\psi_P$. Points of $\d \acal_P$ are always
transition points, and there will be others whenever $P$ has a
face $F\subset p\Si^\circ$ of codimension at least 2. For example,
consider the case where $P$ is the unit square.  In this case,
there are two interior faces whose flow-outs are the connected
components of the classically forbidden region (see
Figure~\ref{allowed-square} above) and the set of transition
points equals $\d \acal_P$.  We shall also give an example with an
interior vertex (see \S \ref{s-ex3}), where the forbidden region
is connected but decomposes into flow-outs of 3 faces (see
Figure~\ref{Fn}) and there are transition points not in $\d
\acal_P$.

The form $\psi_P$ not only encodes the expected (normalized)
density of the zero set, but also the expected density of tangent
directions to the zero set. In the course of the proof of Theorem
\ref{main}, we will show that in the forbidden region, the limit
tangent directions are restricted.  In particular, as the polytope
expands, the tangent spaces to typical zero sets approach tangency
to the `normal flow' $\{e^{\tau+i\theta}\cdot z^0: \tau,\theta\in
T_F^\perp\subset\R^m\}$. A precise formulations of this fact is
given in Theorem \ref{more}. Thus, while the expected distribution
of zero densities is absolutely continuous, the expected
distribution of zero tangents is singular.

Finally, we consider the general case of  $k\le m$ independently
chosen random polynomials
$$f_j\in \poly(P_j)\;,\qquad P_j\subset p_j\Si\;,\qquad
P_j\not\subset \d (p_j\Si) \qquad (1\le j\le k)\;,$$ and we let
$\E_{|P_1,\dots,P_k}(Z_{f_1,\dots,f_k})$ denote the expected zero
current with respect to the probability measure
$\ga_{p_1|P_1}\times \cdots\times \ga_{p_k|P_k}$ on the product
space. If $P_1=\cdots = P_k=P$, then we also write
$\E_{|P,\dots,P}(Z_{f_1,\dots,f_k}) = \E_{|P}(Z_{f_1,\dots,f_k}) $.

\begin{theo}\label{simultaneous} Let $P_1,\dots,P_k$ be
convex integral polytopes in $\R^m_{\ge 0}$.  Then
$$N^{-k}\E_{|NP_1,\dots,NP_k}(Z_{f_1,\dots,f_k}) \to
\psi_{P_1}\wedge\cdots\wedge\psi_{P_k} \quad \mbox{in\ \ }
\lcal^1_{\rm{loc}}((\C^*)^m)\;, \quad \mbox{as\ \ } N\to
\infty\;.$$ \end{theo}

Theorem \ref{probK} is a consequence of Theorems~\ref{main} and \ref{simultaneous}.
Furthermore, the expected volume of the zero set of a system of
$k$ polynomials has the following exotic distribution law, as given by
the following corollaries:

\begin{cor} \label{zerovolumes} Let $P_1,\dots,P_k$ be convex
integral polytopes in $\R^m_{\ge 0}$. Then
for every relatively compact,
open set $U\subset(\C^*)^m$, we have
$$\frac{1}{N^k}\E_{|NP_1,\dots,NP_k} \vol (|Z_{f_1,\dots,f_k}|\cap U) \to
\frac{1}{(m-k)!}\int_U \psi_{P_1}\wedge\cdots\wedge\psi_{P_k}
\wedge\om^{m-k}_\FS\;.$$
\end{cor}

\begin{cor} \label{subtler} The expected zero current $N^{-k}\E_{|NP}
(Z_{f_1, \dots, f_k})$ tends to $0$ at all points of each
forbidden subregion $\rcal_F$ with $\dim F <k$.
\end{cor}

\medskip

\subsection{Examples}\label{examples}

We illustrate the results of Theorems \ref{MASS}--\ref{main} in
some simple cases in dimension $m=2$. For our examples, we shall
describe the forbidden subregions $\rcal_F$ given by
(\ref{flow})--(\ref{RF}) and compute the decay function $b_P(z)$
and the limit expected zero current $\psi_P$.

To simplify our computation of $\psi_P$ from (\ref{psiP}), we
write
$$u(z) = -b_P(z) +p
\log (1+\|z\|^2)=p\log (1 +\|e^{-\tau_z/2} \cdot z\|^2) + \langle
q(z), \tau_z\rangle\;,$$ so
that
\begin{equation}\label{psi-formula} \psi_P =
\frac{\sqrt{-1}}{2\pi }\partial\bar\partial u =
\frac{\sqrt{-1}}{2\pi}\partial\bar\partial\left[p \log (1
+\|e^{-\tau_z/2} \cdot z\|^2) + \langle q(z),
\tau_z\rangle\right]\;.\end{equation}

\subsubsection{Example 1: the square}\label{example1} For our first example, we let
$P$ be the unit square with vertices
$\{(0,0),\,(1,0),\,(0,1),\,(1,1)\,\}$ and we let $p = 2$. Recalling that
$$\mu(z_1,z_2)=\left(\frac{|z_1|^2}{1+|z_1|^2 +|z_2|^2}\,,\
\frac{|z_2|^2}{1+|z_1|^2 +|z_2|^2}\right)\;,$$ we see that the
classically allowed region is given by
$$\acal_P=\{(z_1,z_2):|z_1|^2 -1 < |z_2|^2< |z_1|^2 +1\}\;,$$ as illustrated in
Figure~\ref{allowed-square}. The forbidden region consists of two
subregions:
$$\begin{array}{ll}\rcal_F=\{(z_1,z_2):|z_2|^2 \ge |z_1|^2 +1\}\;,\qquad &
 F=\{(x_1,1):0\le x_1
\le
1\}\;,\\[8pt]\rcal_{F^*}=\{(z_1,z_2):|z_2|^2 \le |z_1|^2 -1\}\;,\qquad &
 F^*=\{(1,x_2):0\le x_2
\le 1\} \;.\end{array}$$

Suppose that $z$ is a point in the upper forbidden region
$\rcal_F$. Recalling (\ref{cond1})--(\ref{cond2}), we write
$\tau_z=(\tau_1,\tau_2)$; then $\tau_1=0$ since $\tau\perp T_F$.
The boundary point $q(z)$ of $P$ whose normal flow contains $z$ is
given by
$$ q(z)= 2 \mu (e^{-\tau_z/2}\cdot z) = (a,1)\in F\subset \d
P\;.$$
  Writing
$$ |z_1|^2 =s_1,\quad |z_2|^2 = s_2\;,\quad |e^{-\tau_2/2}
z_2|^2=e^{-\tau_2} s_2 =\tilde s_2\;,$$ we have $$\frac{
s_1}{1+s_1+\tilde s_2} =\frac{a}{2}\;,\qquad \frac{\tilde
s_2}{1+s_1 +\tilde s_2}=\frac{1}{2}\;.$$ Therefore
$$\begin{array}{c}\di s_1 =\frac{a}{1-a},\quad \tilde s_2
=\frac{1}{1-a}= \frac{s_1}{a},\quad a= \frac{s_1}{1+s_1} =
\frac{|z_1|^2}{1+|z_1|^2},\\[12pt]\di e^{-\tau_2} =\tilde s_2/s_2
=\frac{s_1}{as_2} = \frac{1+|z_1|^2}{|z_2|^2}\;.\end{array}$$
 We have $$
\log (1+\|e^{-\tau_z/2}\cdot z\|^2) =
 \log\left (1+|z_1|^2 + \frac{1+|z_1|^2}{|z_2|^2}
|z_2|^2\right)= \log (1 +|z_1|^2) +\log 2\;,$$
$$\langle q(z), \tau_z\rangle =\left\langle\left(\frac{|z_1|^2}{1+|z_1|^2}
,1\right), \left(0,\log
\frac{|z_2|^2}{1+|z_1|^2}\right)\right\rangle =\log |z_2|^2- \log
(1+|z_1|^2) \;.$$ Therefore $$u= \log|z_2|^2 + \log
(1+|z_1|^2) +\log 4\;.$$

We conclude that $$  \psi_P= \left\{\begin{array}{lll}
\frac{\sqrt{-1}}{2\pi} \partial\bar\partial \log (1+|z_1|^2)\quad
&
\mbox{for }\ z \in \rcal_F^\circ \quad &(|z_1|^2 +1 < |z_2|^2)\\[14pt]
  \frac{\sqrt{ -1}}{\pi} \partial\bar\partial \log
(1+|z_1|^2+|z_2|^2) = 2\om_\FS\quad &
\mbox{for }\ z\in\acal_P & ( |z_1|^2-1< |z_2|^2< |z_1|^2+1)\\[14pt]
 \frac{\sqrt{ -1}}{2\pi} \partial\bar\partial \log
(1+|z_2|^2)\quad & \mbox{for }\ z \in \rcal_{F^*}^\circ \quad &( |z_2|^2<
|z_1|^2-1)\end{array}\right.$$ (where the third case is by
symmetry).  Note that $\psi_P$ has constant rank 1 in both of the
forbidden regions $\rcal_F,\ \rcal_{F^*}$, as indicated in Theorem
\ref{main}(iii).

On $\rcal_F$, we have:
$$ e^{-b_P(z)} = \frac{4|z_2|^2(1+|z_1|^2)}{(1+|z_1|^2+|z_2|^2)^2} \qquad
\mbox{for }\ |z_1|^2 +1< |z_2|^2\;.$$ On the boundary $\{|z_1|^2 =
|z_2|^2-1$\}, we have $e^{-b_P(z)}=1$ as expected. On $\{|z_1| =
c\}$, we have the growth rate $e^{-b_P(z)} \sim 1/|z_2|^2$ as
$z_2\rightarrow\infty$. To obtain $e^{-b_P(z)}$ on $\rcal_{F^*}$,
we interchange $z_1$ and $z_2$ in the above.

\subsubsection{Example 2: a trapezoid}
We next consider the trapezoidal polytope of Figure~\ref{F2} below.
  Comparing with the case of the square, we see that the
classically allowed region is given by $|z_2|^2 < |z_1|^2 +1$, and
the forbidden region coincides with the upper forbidden region
$\rcal_F=\{|z_2|^2 \ge |z_1|^2 +1\}$ from Example~1. Thus, the map
$z\mapsto (\tau_z,q(z))$ is the same as before when $z$ is in the
forbidden region $\rcal_F$. Hence, $e^{-b_P(z)}$ and $\psi_P$ are
also the same as in Example~1 on $\rcal_F$.  On the classically
allowed region, $e^{-b_P(z)}=1$ and $\psi_P=2\om_\FS$.

\begin{figure}[htb]
\centerline{\includegraphics*[bb= 1.9in 4.9in 4in 6.9in]{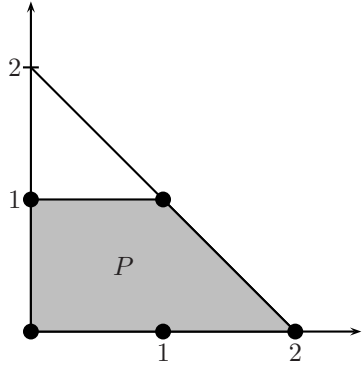}}
\caption{$P\cap\Z^2=\{(0,0),\,(1,0),\,(2,0),\,(0,1),\,(1,1)\}$}
\label{F2}\end{figure}

\subsubsection{Example 3: trapezoids of
higher degree} \label{s-ex3} Now let $n\ge 2$ and let $P$ be the
trapezoid with vertices $(0,0),(n+1,0),(0,1),(1,1)$ given in
Figure~\ref{Fn} below.  In this case, $p=n+1$ and $P$ has an
interior vertex $v=(1,1)\in\Si^\circ$.

We see that the classically allowed region
$\mu\inv(\frac{1}{p}P^\circ)$ is given by
\begin{equation}\label{r1}\acal_P=\left\{(z_1,z_2):|z_2|^2< \min\left\{\frac{|z_1|^2
+1}{n},\,
\frac{1}{n-1}\right\}\,\right\}\;.\end{equation}  This time, the forbidden region
consists of three subregions:
$\rcal_F,\ \rcal_v,\ \rcal_{F'}$, where $$\textstyle
F=\{(x_1,1) :0\le x_1<1\} \;,\qquad
F'=\{(x_1,x_2): x_2=\frac{1}{n}(n+1-x_1),\ 1 <x_1\le
n+1\}\;.$$  (See Figure~\ref{Fn} below.)

\begin{figure}[htb]
\centerline{\includegraphics*[bb= 1.5in 4.9in 4.4in 6.8in]{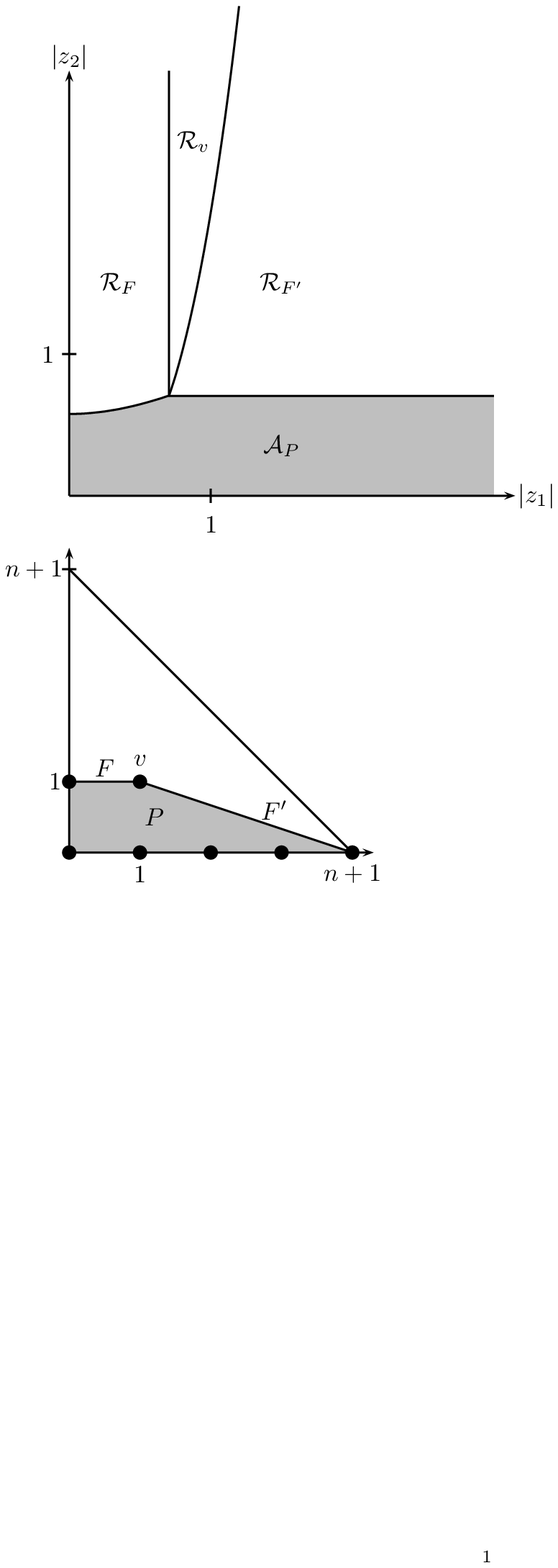}
\includegraphics*[bb= 1.6in 6.8in 5in 9.6in]{f1a.ps}}
\caption{$P\cap\Z^2=\{(0,0),\,(1,0),\,\dots,(n+1,0),\,(0,1),\,(1,1)\}$}
\label{Fn}\end{figure}

Suppose that $z$ is a point in the region $\rcal_F$. Then
$\tau_z=(0,\tau_2)$,  and $$ q(z)= (n+1) \mu (e^{-\tau_z/2}\cdot
z) = (a,1)\in (n+1)F \subset \d P \qquad (0<a<1)\;.$$ Again
writing
$$ |z_1|^2 =s_1,\quad |z_2|^2 = s_2\;,\quad |e^{-\tau_2/2}
z_2|^2=e^{-\tau_2} s_2 =\tilde s_2\;,$$ we have $$\frac{
s_1}{1+s_1+\tilde s_2} =\frac{a}{n+1}\;,\qquad \frac{\tilde
s_2}{1+s_1 +\tilde s_2}=\frac{1}{n+1}\;.$$ Therefore
$$\begin{array}{c}\di s_1 =\frac{a}{n-a},\quad \tilde s_2
=\frac{1}{n-a}= \frac{s_1}{a},\quad a= \frac{ns_1}{1+s_1} =
\frac{n|z_1|^2}{1+|z_1|^2},\\[12pt]\di e^{-\tau_2} =\tilde s_2/s_2
=\frac{s_1}{as_2} = \frac{1+|z_1|^2}{n|z_2|^2}\;.\end{array}$$ In
particular,
$$ a<1 \ \Leftrightarrow\ |z_1|^2 < \frac{1}{n-1}$$
and therefore
\begin{equation}\label{r2}\rcal_F=\left\{(z_1,z_2): |z_2|^2\ge \frac{|z_1|^2 +1}{n},\
|z_1|^2 <\frac{1}{n-1}\right\}\;.\end{equation}

We have $$ \log (1+\|e^{-\tau_z/2}\cdot z\|^2) =
 \log\left (1+|z_1|^2 + \frac{1+|z_1|^2}{n|z_2|^2}
|z_2|^2\right)= \log (1 +|z_1|^2) +\log \frac{n+1}{n}\;,$$
$$\langle q(z), \tau_z\rangle =\left\langle\left(\frac{n|z_1|^2}{1+|z_1|^2}
,1\right), \left(0,\log
\frac{n|z_2|^2}{1+|z_1|^2}\right)\right\rangle = \log |z_2|^2-\log
(1+|z_1|^2) +\log n\;.$$ Therefore $$u= \log|z_2|^2 +n \log
(1+|z_1|^2) +(n+1)\log (n+1) -n\log n\;.$$ Hence,
$$\left.\begin{array}{rcl} \psi_P&=&\di n \frac{\sqrt{-1}}{2\pi}
\partial\bar\partial
\log (1+|z_1|^2)\\[14pt] e^{-b_P(z)} &=&\di \frac{(n+1)^{n+1}}{n^n}
\frac{|z_2|^2(1+|z_1|^2)^n}{(1+|z_1|^2+|z_2|^2)^{n+1}}
\end{array}\right\} \qquad \mbox{for }\ z\in \rcal_F\;.$$

Now suppose that $z$ is a point in $\rcal_{F'}$.  Since $\tau_z
\perp T_{F'}$, we can write $\tau_z=(\tau_1, n\tau_1)$. Let
\begin{equation} q(z)=  (n+1)\mu (e^{-\tau_z/2}\cdot z) =
\left(c, \frac{1}{n}(n+1-c)\right)\in
F'\;.\label{slant}\end{equation} As before, we write
$$s_1=|z_1|^2,\quad s_2=|z_2|^2,\qquad \tilde s_1=|e^{-\tau_1/2}z_1|^2=e^{-\tau_1}
s_1,\quad \tilde s_2=|e^{-n\tau_1/2}z_2|^2=e^{-n\tau_1}s_2\;.
$$  By (\ref{slant}), we have
$$\frac{\tilde s_1}{1+\tilde s_1 +\tilde s_2} = \frac c{n+1}\;,\qquad
\frac{\tilde s_2}{1+\tilde s_1 +\tilde s_2} =
\frac{1}{n}\left(1-\frac c{n+1}\right)\;.$$ Solving for $\tilde s_1,\tilde s_2$, we
obtain
\begin{equation}\label{stilde}\tilde s_1 = \frac{n}{n-1}\; \frac{c}{n+1-c}\;,
\qquad \tilde s_2 = \frac{1}{n-1}\;.\end{equation}  Therefore,
$e^{-n\tau_1}= {\tilde s_2}/{s_2}= {|z_2|^{-2}}/(n-1)$, so we have
$$\tau_1= \frac{1}{n}\log (n-1)|z_2|^2\;,\qquad \tilde
s_1=\frac{|z_1|^2}{(n-1)^{1/n} |z_2|^{2/n}}\;.$$ Thus, $$ \log
(1+\|e^{-\tau_z/2}\cdot z\|^2) =
 \log (1+\tilde s_1 +\tilde s_2)=\log \left(\frac{n}{n-1} +
\frac{|z_1|^2}{(n-1)^{1/n} |z_2|^{2/n}}\right)\;.$$

By (\ref{slant}),
$$\big\langle q(z), \tau_z\big\rangle=
\left\langle\left(c, \frac{1}{n}(n+1-c)\right),
(\tau_1,n\tau_1)\right\rangle = (n+1)\tau_1 = \frac{n+1}{n}\log
(n-1)|z_2|^2\;.$$ Hence by (\ref{psi-formula}),
$$\psi_P=(n+1)\frac{\sqrt{-1}}{2\pi}\ddbar\log \left(\frac{n}{n-1} +
\frac{|z_1|^2}{(n-1)^{1/n} |z_2|^{2/n}}\right) \qquad \mbox{for }\
z\in \rcal_{F'}\;.$$ (Note that $\psi_P$ has constant rank 1 on
$\rcal_{F'}$ as indicated by Theorem \ref{main}(iii), since
$\psi_P =(n+1)g^*\om_{\CP^1}$  on $\rcal_{F'}$, where $g$ is the
multi-valued holomorphic map to $\CP^1$ given by $g(z_1,z_2)= (c_n
z_1, z_2^{1/n})$.)

By Theorem \ref{main}(iii), we know that $\psi_P=0$ on $\rcal_v$.
To complete the description of $\psi_P$, it remains to describe
the regions $\rcal_{F'}$ and $\rcal_v$. We note that a forbidden
point $z$ lies in $\rcal_{F'}$ if and only if $c>1$. By (\ref{stilde}), this is
equivalent to $\tilde s_1 > \frac{1}{n-1}$, or $$|z_2|^2 < (n-1)^{n-1} |z_1|^{2n}\;.$$
Therefore
\begin{equation}\label{r3}\rcal_{F'}=\left\{(z_1,z_2): \frac{1}{n-1} \le |z_2|^2 <
(n-1)^{n-1} |z_1|^{2n},\ \ |z_1|^2 >\frac{1}{n-1}\right\}\;.\end{equation} This leaves
us with
\begin{equation}\label{r4}\rcal_v=\left\{(z_1,z_2):  |z_2|^2 \ge  (n-1)^{n-1}
|z_1|^{2n},\ \ |z_1|^2 \ge\frac{1}{n-1}\right\}\;.\end{equation}
These subregions are illustrated in Figure \ref{Fn}.
To summarize:
$$  \psi_P=
\left\{\begin{array}{lll} (n+1) \frac{\sqrt{ -1}}{2\pi}
\partial\bar\partial \log (1+|z_1|^2+|z_2|^2) = (n+1)\om_\FS \quad & \mbox{for }\
z\in \acal_P \quad & \mbox{(see (\ref{r1}))}\\[14pt]
n \frac{\sqrt{-1}}{2\pi}
\partial\bar\partial
\log (1+|z_1|^2) & \mbox{for }\  z\in\rcal^\circ_F
& \mbox{(see (\ref{r2}))}\\[14pt](n+1)
\frac{\sqrt{-1}}{2\pi}\ddbar\log \left(\frac{n}{n-1} +
\frac{|z_1|^2}{(n-1)^{1/n} |z_2|^{2/n}}\right) & \mbox{for }\
z\in \rcal^\circ_{F'} & \mbox{(see (\ref{r3}))}\\[14pt]
 0 &
\mbox{for }z\in \rcal^\circ_v& \mbox{(see (\ref{r4}))}\end{array}\right.$$
Thus, as stated in Theorem~\ref{main}, $\psi_P$ has constant rank 1 on the
flow-outs  of the one-dimensional faces $F,\ F'$, and vanishes on the flow-out of the
vertex $v$.

\subsection{Methods of proof}\label{s-methods}
Having stated our principal results, we now briefly outline some
key ideas in the  proofs. The key result is Theorem \ref{MASS} on
the mass of polynomials with Newton polytope $P$. To prove it, we
begin with an easy formula
\begin{equation}\label{2measures}\E_{\nu_{NP}}(|f(z)|^2_\FS) =
\frac{1}{\# (NP)}\E_{|NP} \left(|f(z)|^2_\FS\right) =\frac{1}{\#
(NP)}\Pi_{|NP}(z,z)\;,\end{equation} (see \S \ref{s-mass})
relating expected mass to the  {\it conditional \szego kernel\/}
$\Pi_{|NP}$, i.e. the orthogonal projection onto $\poly(NP)$. In
general,  the term `\szego kernel' of a space $\scal$ of
$\lcal^2$ functions refers to the kernel for the
orthogonal projection to $\scal$ from the space of all $\lcal^2$
functions;  i.e., it is of the form $\Pi(x,y)=\sum _j s_j(x)
\overline{s_j(y)}$, where $\{s_j\}$ is an orthonormal basis of
$\scal$.  Precise asymptotics for $\frac{1}{\#
(NP)}\Pi_{|NP}(z,z)$ are given  in Theorem~\ref{SZEGO}. Note
that the normalizing factor $\frac{1}{\# (NP)}$ is straightforward
to evaluate since
$$\# (NP) =\dim \poly(NP) =
 \vol(P) N^m +\dots\,,$$
 is the Ehrhart polynomial \cite{Eh}.

To analyze the conditional \szego kernels $\Pi_{|NP}(z, w)$, we
introduce  the {\it polytope characters}
\begin{equation}\label{PC} \chi_{NP}(e^{i \phi}): = \sum_{\alpha
\in NP} e^{i \langle \alpha, \phi \rangle}\;,\qquad e^{i\phi}=
(e^{i\phi_1},\dots, e^{i\phi_m})\;.
\end{equation} We observe that
\begin{equation} \label{CSCH} \Pi_{|N P} (z,w) = \frac 1{(2\pi)^m} \int_{\T }
\Pi_{Np} ( z,e^{i\phi}\cdot w) {\chi_{NP}(e^{i \phi})}\, d\phi\;,
\end{equation}
where \begin{equation}\label{torus}\T=\{(\zeta_1,\dots\,\zeta_m)\in (\C^*)^m:
|\zeta_j|=1, 1\le j\le m\}\end{equation} is the $m$-torus, and where $\Pi_{Np}(z,w)$ is
the
\szego kernel of $\poly(Np\Si)$,  the space of all
homogeneous polynomials of degree $N p$. To obtain asymptotics, we
need to analyze the behavior of $\chi_{N P}(e^w)$ as $N \to
\infty$.

To do so, we  use  the Euler-MacLaurin formula of
Khovanskii-Pukhlikov \cite{KP}, Brion-Vergne \cite{BV,BV2}, and
Guillemin \cite{Gu2}:
 \begin{equation} \chi_P(e^w) =
 \mbox{Todd}(\fcal_P,\partial/\partial h)
\left. \left(\int_{P(h)} e^{\langle w,x \rangle}
dx\right)\right|_{h = 0} \qquad (w\in \C^m,\ \|w\|<\ep)\;,
\label{KPBV}\end{equation} where $P(h)$ is of the form $\{x:
\langle u_j, x\rangle + \la_j + h_j \geq 0,\ 1 \leq j \leq n\}$
($P(0)=P$) and where Todd$(\fcal_P,\partial/\partial
 h)$ is a certain infinite order differential {\it Todd operator}.
Upon dilating the polytope, one obtains
\begin{equation}\label{CHARACTER} \chi_{N P}(e^w) =  N^m  \;
\mbox{Todd}(\fcal_P,N\inv\partial/\partial h) \left.
\left(\int_{P(h)} e^{ N\langle w,x \rangle} dx\right)\right|_{h =
0}\qquad (\|w\|<\ep).
 \end{equation}
In  Proposition
 \ref{newchar}, we show that the family $\chi_{N P}$ is a complex
 oscillatory integral of the form:
\begin{equation} \label{CHARFOR} \chi_{NP}(e^w) = \int_P e^{N\langle w, x \rangle}[ A_0(x,w)
N^n+ A_1(x,w)N^{n-1} +\cdots +A_n(x,w)]\,d\vol_n(x)\;,
\end{equation} for $\Im w$ sufficiently small (but $\Re w$ arbitrary),
where the $A_l$ are analytic functions that are holomorphic in $w$
and algebraic in $x$

 We then substitute this expression in (\ref{CSCH}) and use the method of stationary phase for
complex oscillatory integrals \cite[Ch.~7]{H} to obtain the
asymptotics of $\Pi_{|N P}(z,z)$.  For the case where $z$ is in
the classically allowed region, we easily find that the critical
point of the phase is given by $\phi=0$ and $x=p\mu(z)$. Since the
phase vanishes and has nondegenerate Hessian at the critical
point, we immediately obtain an asymptotic expansion. The case
where $z$ is in the classically forbidden region is more subtle.
Since $p\mu(z)$ lies outside of $P$ for this case, the phase has
no critical points. To complete the analysis, we must deform the
contour to pick up critical points. In particular, we consider the
complexification $(\C^*)^m$ of $\T$ and deform $\T$ to a contour
of the form $(\log |\zeta_1|,\dots, \log|\zeta_m|) =\tau\in\R^m$,
on which the phase has a `critical point' with $\phi=0$ and
$x=x_c\in\d P$. To be precise, the derivative tangential to the
face of $P$ containing $x_c$ vanishes at $x_c$ (while the normal
derivative is nonvanishing), and furthermore the phase takes its
maximal real part at $x_c$. Indeed, $\tau$ is the unique vector
$\tau_z$ and $x_c$ is the unique point $q(z)\in\d P$ used in
formula (\ref{b}) for the decay rate $b_P(z)$, and the maximal
real part of the phase is $-b_P(z)$. We then obtain an asymptotic
expansion from the method of stationary phase on domains with
boundary.

\subsection{Brief outline and remarks}\label{s-remarks} To
assist the reader in navigating this article, we give a brief
outline:  We begin in Section~\ref{s-background} by reviewing some
notation and terminology involving convex polytopes and \szego
kernels for spaces of polynomials. Section~\ref{s-char} states and
proves our oscillatory integral formula for the characters
$\chi_{NP}(e^w)$. The heart of the paper is Section~\ref{s-mass},
where we derive the diagonal asymptotics of the \szego kernel
(Theorem~\ref{SZEGO}) for the spaces $\poly(NP)$, using our
character formula from Section~\ref{s-char} as one of the
ingredients. (Theorem~\ref{MASS} is an immediate consequence of Theorem~\ref{SZEGO}.)
Then applying our
\szego kernel asymptotics, we prove Theorems \ref{probK}, \ref{main} and
\ref{simultaneous} on the distribution of zeros in Section~\ref{DZ}. We also include an
appendix (Section~\ref{Appendix}), where we describe an
alternative geometric approach to our results using toric
varieties.  An index of notation is included at the end of the
paper.

Before embarking on the proofs, we comment on our  borrowed
terminology from physics  and on the role of toric geometry:

\subsubsection{Tunnelling of zeros.}

Some of the terminology we use---allowed and forbidden regions, mass density---is taken
from the semiclassical analysis of ground states of Schr\"odinger operators $H_\hbar =
-\hbar^2
\Delta + V$ on $\R^n$. The well-known `Agmon estimates' of ground states (cf.\
\cite{Ag}) show that $\lcal^2$-normalized ground states or
low-lying eigenfunctions $H_\hbar \phi = E \phi$ of $H$ are
concentrated as $\hbar \to 0$  in the classically allowed region
$C_E = \{x \in \R^n: V(x) \leq E\}$ and $|\phi(x)|^2$ decays
exponentially
 as $\hbar \to 0$ for points $x$ in the
complement.  In our setting, the Hamiltonian is $\bar{\partial}^*
\bar{\partial}$ on $\lcal^2$-sections of powers  of the hyperplane
section bundle on $\CP^m$, $\hbar = 1/N$, and the ground states
are the holomorphic sections. Replacing $|\phi(x)|^2$ is the
expected  mass density $|f(z)|^2_\FS$ of the polynomials in this
subspace.  Theorems \ref{probK}--\ref{MASS} show that the
Newton-polytope  constraint on the polynomials creates a
`tunnelling theory' for zeros and mass.

The decay function $b_P$ in Theorem \ref{MASS} is analogous to the
Agmon action to the allowed region \cite{Ag}.   We will also
obtain (Proposition \ref{b-prop}) the integral formula
\begin{equation}\label{b-action} b_P(z) =
\int_0^{\tau_z}\left[-q(e^{-\sigma/2}\cdot z) +
p\mu(e^{-\sigma/2}\cdot z)\right] \cdot d\sigma \end{equation}
(where $\int_0^{\tau_z}$ denotes the integral over any path in
$\R^m$ from $0$ to $\tau_z$),  which could be interpreted as an
action, thereby bringing the results closer to classical Agmon
estimates.

\subsubsection{Role of toric varieties}
Since $\poly(NP)$ is naturally isomorphic to the space
$H^0(M_P, L_P)$ of holomorphic sections of
a natural line bundle $L_P$  over a toric variety $M_P$ associated
to $P$ (see \S \ref{Appendix}), the reader may wonder whether toric varieties play any
role in this paper.  The answer is that neither the statements nor proofs of our results
involve toric varieties in any essential way. However, the theory of these varieties
does give an alternative approach to the asymptotics of the conditional \szego kernel,
as will be explained   \S\ref{Appendix}.  It was also the approach in the original
version of this paper \cite{SZ2} and motivated some of the ideas.

The reader may also  wonder how the results of this paper would change
if,  instead of defining  the  Gaussian measures $\gamma_{|P}$ to be the
conditional measures of one fixed ensemble, we defined  the  measures
on  $\poly(NP) \simeq H^0(M_P, L_P^N)$ to be the   Gaussian measures $\ga_N^{M_P}$
induced by the $\lcal^2$ inner product induced by a Hermitian
metric on $L_P$ and its curvature form $\omega_{M_P}$.
As a special case of our results in
\cite[Prop.~4.4]{SZ},  one obtains
(under the added assumption that $M_P$
is smooth or maybe has orbifold singularities)
\begin{equation}\label{action} \textstyle\frac{1}{N^k }
\E_{\ga_N^{M_P}}  (Z_{f_1, \dots, f_k}) =\om_{M_P}^k +
O\left(\frac{1}{N}\right)= \left(\sum_{j=1}^m dI_j \wedge
d\theta_j\right)^k + O\left(\frac{1}{N}\right)\;,\end{equation}
where $I_j,\theta_j$ are the action-angle variables of the  moment
map $\mu_P: M_P \to \R^m$ of the $\T $-action on $M_P$.
Related
results in a somewhat different set-up  have also been obtained by
Malajovich and Rojas \cite{MaR}. This  quite  different law shows
that the  measures $\ga_N^{M_P}$ are singular relative to $\ga_{Np|
N P}$ in the limit as $N \to \infty$.

As mentioned previously,  the polytope $P$  in this article
is only used as  a constraint on the polynomials in creating conditional measures. The
norms of the monomials $z^{\alpha}$  are fixed (as their $\CP^m$ norms). Thus,
the change in the distribution of zeros as $P$ varies is due solely to the choice of
which monomials occur in the polynomials.
In the case of $\gamma^{M_P}_N$,  the norms of the monomials vary as $P$ varies
since they are $\lcal^2$-normalized on the toric variety $M_P$. Hence the
variances of the coefficients of a polynomial in the $\gamma^{M_P}_N$
ensembles  depend on the choice of metric on $L_P$.
This dependence creates complicated biases towards some
monomials and away from others as $P$ varies, making it  difficult to understand what
 a comparison between the  $\gamma^{M_P}_N$ ensembles  would be  measuring.

In fact, our results using conditional measures also apply to
polytopes that are not convex and hence do not correspond to any
toric variety. Indeed, suppose that $P$ is a nonconvex `lattice'
polytope.  Then Theorems \ref{probK}--\ref{simultaneous} hold with
the following modifications:  The point $q(z)$ satisfying
(\ref{cond1})--(\ref{cond2}) is not always unique; instead we
choose $q(z)$ to minimize $b_P(z)$.  Then the function $b_P(z)$ is
$\ccal^0$, not $\ccal^1$, and $\psi_P$ is a positive current,
which has singular support.  Part (iii) of Theorem~\ref{main}
applies only to the absolutely continuous part of $\psi_P$. More
significantly,  the limit measure in Theorem~\ref{probK} does not
vanish on the forbidden region; instead, it is orthogonal to
 volume measure there.  This observation is relevant to the theory of
`fewnomial' systems, which we discuss in a forthcoming paper
\cite{SZ4}.

\subsubsection{Acknowledgements} Our interest  in
polynomials with a fixed Newton polytope was in part stimulated by
a discussion with A. Varchenko at the outset of this work.  We
would like to thank M. Brion for many helpful comments regarding
polytopes and the Euler-MacLaurin formula. We also grateful to T.
Theobald for giving us permission to use Figure \ref{fig-amoeba}
from his paper and A. Carass for providing the illustrations in
Figures \ref{allowed-square}, \ref{F2}, \ref{Fn}, \ref{f-fan}, and \ref{L}.

\section{Background}\label{s-background}

\subsection{Fans}\label{s-fans}
By a {\it convex integral polytope\/}, we mean the convex hull in
$\R^m$ of a finite set  in the lattice $\Z^m$. A convex integral
polytope $P$ with nonempty interior $P^\circ$   can be defined by
linear equations
\begin{equation}\label{Pdef}\ell_j(x): = \langle x, u_j \rangle +\la_j \geq 0, \;\;\; (j
=  1, \dots, d)\;,\end{equation} where $u_j \in \Z^m$ is the
primitive inward-pointing normal to the $j$-th {\it facet\/}
(codimension-one face) $$F_j^{m-1}:=\{x\in P: \ell_j(x)=0,\
\ell_k(x)>0 \ \mbox{for }\ k\ne j\}\;.$$ For each point $x\in P$,
we consider the {\it normal cone to $P$ at $x$\/},
\begin{equation}\label{Cx}C_x =C_x^P: =
\{u\in \R^m: \langle u, x\rangle =\sup _{y\in P}  \langle
u,y\rangle\}\;,\end{equation} which is a closed convex polyhedral cone. We
decompose $P$ into a finite union of {\it faces\/}, each face
being an equivalence class under the equivalence relation $x\sim y
\iff C_x = C_y$.  For each face $F$, we let $C_F$ denote the
normal cone of the points of $F$.  {\bf Note that by our
convention, the faces are disjoint sets.}  We shall use the term
{\it closed face\/} to refer to the closure of a face of $P$.

Each face of dimension $r$ ($0\le r\le m$) is an open polytope in
an $r$-plane in $\R^m$; i.e., the $0$-dimensional faces are the
vertices of $P$, the 1-dimensional faces are the edges with their
end points removed, and so forth. The facets $F_j=F_j^{m-1}$ and their
normal cones are given by:
$$\bar F_j=\{x\in P:
\ell_j(x)=0\}\;,\qquad C_{F_j}=\{-tu_j:t\ge 0\}\;.$$ The
$m$-dimensional face is the interior  $P^\circ$ of the polytope
with normal cone $C_{P^\circ}=\{0\}$.

For each $x\in P$, we let \begin{equation}\label{J}\jcal(x) =
\{j\in\Z: 1\le j\le d,\ \ell_j(x)=0\}\;.\end{equation} One easily
sees that $C_x=C_y \iff \jcal(x)=\jcal(y)$, and hence we can write
$\jcal(F)=\jcal(x)$, where $x$ lies in the face $F$. In
particular, $\jcal(F^{m-1}_j)=\{j\}$ and
$\jcal(P^\circ)=\emptyset$. The polytope $P$ is called {\it
simple\/}
 if  $\#\jcal(v)=m$ for each vertex $v$ of $P$.  In this case,
$\{u_j:j\in \jcal(v)\}$ is a basis for $\R^m$ for each vertex $v$
(since $P^\circ\ne\emptyset$), and furthermore $\#\jcal(F)=\codim
F$ for all faces $F$. The polytope $P$ is said to be {\it
Delzant\/}
 if  $\{u_j:j\in \jcal(v)\}$ generates the lattice $\Z^m$ for each vertex $v$ of $P$.

\begin{figure}[htb]
\centerline{\includegraphics*[bb= 1.0in 7.5in 5in 9.6in]{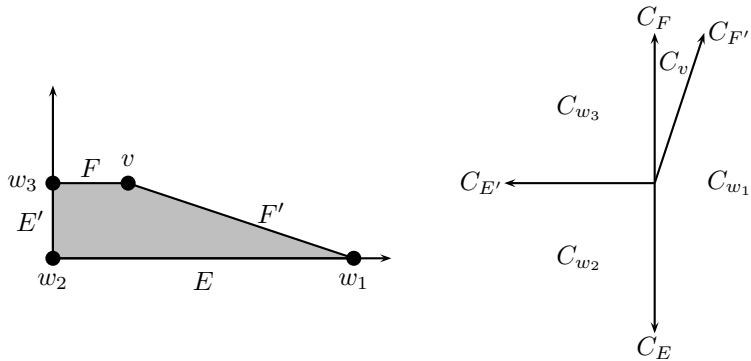}}
\caption{A convex polytope and its fan} \label{f-fan}\end{figure}

The convex integral polytope $P$ determines the fan
$\fcal_P:=\{C_F:F\ \mbox{is a face of} \ P\}$.  A {\it fan\/}
$\fcal$ in $\R^m$ is a collection of closed convex rational
polyhedral cones such that a closed face of a cone in $\fcal$ is
an element of $\fcal$ and the intersection of two cones in $\fcal$
is a closed face of each of them.  (Fans are used in the algebraic
construction of toric varieties; see \cite[\S 1.4]{F}.) An example
of a convex integral polytope and its fan is given in
Figure~\ref{f-fan}.

We shall also consider convex integral polytopes $P$ with empty
interior.  In this case, the faces of $P$ and the normal cones
$C_x^P$ are defined exactly as above. If $P^\circ =\emptyset$,
then the normal cones of $P$ all contain the linear subspace of
$\R^m$ orthogonal to $P$ (contrary to the case
$P^\circ\ne\emptyset$, where the normal cones are {\it pointed\/},
i.e., do not contain any lines in $\R^m$).  A convex polytope in
$\R^m$ of dimension $n<m$ is also said to be {\it simple\/} if it
is linearly isomorphic to a simple polytope in $\R^n$ or
equivalently, if every vertex is the intersection of precisely $n$
closed faces of dimension $n-1$.

\subsection{\szego kernels }\label{s-szego}
By homogenizing the polynomials in $\poly(P)$, we obtain a finite
dimensional subspace $\scal$ of $\lcal^2(S^{2m+1})$.  By the
`\szego projector,' we mean the orthogonal projection
$\Pi_{\scal}:\lcal^2(X)\to {\mathcal S}$. We recall that the
projector $\Pi_{\scal}$ is given by a kernel of the form
\begin{equation}\label{szego}\Pi_{\scal}(x,y) = \sum_{j = 1}^{k_N}
s_j^N(z) \overline{s_j^N(y)}\;,\end{equation} where $\{s_j\}$ is
an orthonormal basis of ${\mathcal S}$.

We now describe  two sequences of \szego kernels, which play a
crucial role in  our main results:

\subsubsection{The projective \szego kernels}

As a first example, we recall that $\poly(p\Si)$ can be identified
with the space of degree-$p$ homogeneous polynomials
 in $m+1$ variables by
identifying  the (non-homogeneous) polynomial
$$f(z_1,\dots,z_m)= \sum _{|\al|\le p} c_\al
z^\al \qquad (z^\al=z_1^{\al_1} \cdots z_m^{\al_m})$$ with the
homogeneous polynomial
$$F(\zeta_0, \dots,
\zeta_m) = \sum_{|\al| \le p} c_\al
\zeta_0^{p-|\al|}\zeta_1^{\al_1} \cdots \zeta_m^{\al_m}\;.$$ We
equip the space $\poly(p\Si)$ with the $\lcal^2$ inner product on
$S^{2m+1}$:
\begin{equation}\label{IP}\langle f, \bar g \rangle = \frac 1{m!}\int _{S^{2m+1}}
F\overline{G} \,d\nu=  \frac{1}{m!} \int
_{\C^m}\frac{f(z)\overline{g(z)}}{(1+\|z\|^2)^p}
\,\om_\FS^m(z),\quad f,g\in  \poly(p\Si),\end{equation}
where \begin{equation}\label{FSkahler}\om_\FS=\frac {\sqrt{-1}}{2\pi}\;\ddbar \log
(1+\|z\|^2)\end{equation} is the Fubini-Study \kahler form on $\C^m\subset\CP^m$.

A basis for $\poly(p\Si)$ consists of the monomials $\chi_\al(z) =
z_1^{\alpha_1} \cdots z_m^{\alpha_m}$, $|\al|\le p$.  The
monomials $\{ \chi_{\alpha}\}$ are orthogonal but not normalized.
Their $\lcal^2$ norms given by the inner product (\ref{IP}) are:
\begin{equation}\label{IP2}\|\chi_\al\|
=\left[\frac{(p-|\al|)!\al_1! \cdots\al_m!}{(p+m)!}\right]^\half\;.
\end{equation}  Note that the norm
$\|\cdot\|$ depends on $p$. Thus we have an orthonormal
basis for $\poly(p\Si)$ given by the monomials
\begin{equation}\label{normchi} \frac{1}{\|\chi_\al\|}\,\chi_\al=
\left[\frac{(p+m)!}{(p-|\al|)!\al_1!\cdots\al_m!}\right]^\half
\chi_\al= \sqrt{\frac{(p+m)!}{p!} {p\choose \al}}\ \chi_\al \
,\qquad |\al|\le p\;.\end{equation}
where
\begin{equation}\label{multinom}
{p\choose\al}= \frac{p!}{(p-|\al|)!\alpha_1!\cdots \alpha_m!}\;.
\end{equation}

  We let
$\wh\chi_\al^p:S^{2m+1}\to \C$ denote the homogenization of
$\chi_\al$:
\begin{equation}\label{chihat}\wh\chi_\al^p(x) =  x_0^{p-|\al|}x_1^{\al_1}\cdots
x_m^{\al_m} \;.\end{equation} Hence the \szego kernel $\Pi_p$ for
the orthogonal projection to $\poly(p\Si)$ is given by:
\begin{equation}\label{szego-proj}\Pi_{p}(x,y)
=\sum_{|\al|\le p}\frac{1}{\|\chi_\al\|^2}\wh \chi_\al(x)
\overline {\wh \chi_\al(y)}  =\frac{(p+m)!}{p!}\langle x,\bar
y\rangle^p \;,\end{equation} for $x,y\in S^{2m+1}$. (The sum
$\sum_{p=0}^\infty \Pi_{p}$ is the usual \szego kernel for the
sphere; see \cite[\S 1.3.1]{BSZ}.)

We also identify the point $z\in(\C^*)^m$ with the lift $x=\frac
1{(1+\|z\|^2)^{1/2}}(1,z_1,\dots,z_m)\in S^{2m+1}$, and we write
\begin{equation}\label{mchi}\wh
\chi^p_\al(z)= \frac{z^\al}{(1+\|z\|^2)^{p/2}}\;.\end{equation}
The \szego kernel can then be written explicitly as
\begin{eqnarray}\Pi_{p}(z,w) =
\frac{(p+m)!}{p!} \frac{\sum_{|\alpha|\le p }{p\choose\al}z^\al
\bar w^\al}{(1+\|z\|^2)^{p/2} (1+\|w\|^2)^{p/2}}  =
\frac{(p+m)!}{p!} \left[\frac{1+\langle z,\bar w\rangle}
{(1+\|z\|^2)^{1/2} (1+\|w\|^2)^{1/2}}\right]^p\;.\label{Sz}
\end{eqnarray}

\subsubsection{The conditional \szego kernels associated to a
polytope $P$}

In this case, the relevant space of polynomials is the subspace
$\poly(P)\subset \poly(p\Si)$ of polynomials with Newton
polytope $P$. Here, we may choose any $p\ge \deg
P:=\max\{|\al|:\al\in P\}$, but we normally choose $p=\deg P$.  In
the conditional \szego kernel, we retain the Fubini-Study inner
product on this subspace. Hence
this example is very similar to the previous one.
 The main difference is that an orthonormal basis of $\poly(P)$ is given by
$$\left\{\frac{1}{\|\chi_\al\|}\,\chi_\al: \al\in P\right\}\;.$$

\begin{defin}The  conditional \szego kernel $\Pi_{|P}$ is the kernel for the orthogonal
projection to $\poly(P)$ with respect to the induced Fubini-Study
inner product:
\begin{equation}\label{Sz2} \Pi_{|P}(x,y) = \sum_{\alpha \in P}
\frac{1}{\|\chi_\al\|^2} \wh\chi^p_\al(x) \overline{\wh
\chi^p_{\alpha}(y)}\;.\end{equation}\end{defin}

When defining the term `random polynomial with fixed Newton
polytope $P$', we wish to use an $\lcal^2$-norm on monomials which
is defined independently of $P$. This explains why the conditional
\szego kernel is the essential one in our problem.  The
conditional \szego kernel can be written explicitly on $\C^m$
as
\begin{equation}\label{Sz2e} \Pi_{|P}(z,w) =\frac{(p+m)!}{p!} \frac{\sum_{\alpha
\in P}{p\choose\al}z^\al \bar w^\al}{(1+\|z\|^2)^{p/2}
(1+\|w\|^2)^{p/2}} \;.\end{equation}

\bigskip
\section{Polytope character}\label{s-char}

This section is devoted to the asymptotic  analysis of the
polytope characters (\ref{PC}). Our interest is in the asymptotics
of the `ray' of characters $\chi_{N P}(e^{w})$ defined on
$(\C^*)^m$ by
\begin{equation} \chi_{N P} (e^w) = \sum_{\alpha \in N P} e^{\langle w,\al\rangle}
\,\;\;\;w \in \C^m. \label{char}\end{equation}
We shall derive a formula
expressing  the  character $\chi_{NP}(e^w)$
 as an oscillatory integral over the original polytope $P$
with the same phase $\langle w, x\rangle$ as in formula
(\ref{KPBV}) of \cite{KP,BV,Gu2}, which is the starting point of
our derivation.  Although (\ref{KPBV}) holds only for small $w$, our
formula holds by analytic extension for all $w$ in a `strip' of
the form
\begin{equation}\label{strip}S(\ep):=\{w\in\C^m: |\Im w_j| <\ep \ \mbox{for }\ 1\le j\le
m\}\;.\end{equation}

Our integral formula for $\chi_{N P} (e^w)$ is given by
the following proposition, which is also of independent interest.

\begin{prop}\label{newchar} Let $P$ be a simple integral polytope in $\R^m$ of
dimension $n$ ($1\le n\le m$).  Then there exists $\ep>0$ such
that the characters $\chi_{NP}$ can be given as integrals over the
polytope $P$ of the form
$$\chi_{NP}(e^w) = \int_P e^{N\langle w, x \rangle}[ A_0(x,w) N^n+
A_1(x,w)N^{n-1} +\cdots +A_n(x,w)]\,d\vol_n(x)\;,\quad \mbox{for all }\
w\in S(\ep)\;,$$ where the $A_l$ are analytic functions on
$P\times S(\ep)$ that are holomorphic in $w$ and algebraic in $x$,
and
\begin{enumerate}
\item[i)] $A_l(x,w)\in\R$ whenever $w\in \R^m$,
\item[ii)] $A_0(x,0)=1$,
\item[iii)] $A_0(x,w)\in\R^+$ whenever $w$ is in the normal cone  to
$P$ at $x$, i.e., whenever $w\in C_x\subset\R^m$.
\end{enumerate}\end{prop}

As discussed in the
introduction, Proposition \ref{newchar}  will later  be
used to obtain asymptotics of the conditional \szego kernel
through  formula (\ref{CSCH}). When $w = 0$,  the character (\ref{char})
equals the number
of lattice points of $NP$ and Proposition~\ref{newchar} reduces to the Ehrhart formula
\cite{Eh}:
\begin{equation}\label{RR}  \# (NP) =\sum_{j=0}^n a_j N^{n-j}\;,\qquad n=\dim P\;,\quad
a_0=\vol_n(P)\;.\end{equation}
We note that one can
easily extend (\ref{RR}) to an arbitrary convex integral polytope $P$ (not
necessarily simple) by triangulating $P$ into $n$-simplices $\{\De_1,\dots,\De_k\}$
whose vertices are vertices of $P$ (see  \cite[pp.~214--216]{GKZ} for an elementary
proof of this fact),  and then applying the inclusion-exclusion principle.

\begin{rem} The
 characters  $\chi_{N P} (e^w)$ have a natural
interpretation in  toric geometry, and $\chi_{N P}(1)=\#(NP)$ can also be given by the
Riemann-Roch formula on the toric variety associated to the polytope $P$.  See \S
\ref{Appendix} for further details and references.\end{rem}

\subsection{Proof of Proposition \ref{newchar}}

We first reduce to the case where $P$ has nonempty interior: If
$\dim P = n<m$, then by a translation we may assume that $0\in P$.
Next by an orthogonal projection to the $n$-plane $\langle
P\rangle$ spanned by $P$, we may assume that $w\in \langle
P\rangle$.  Since $\langle P\rangle \cap\Z^m\approx \Z^n$, we may
assume that $\langle P\rangle = \R^n$.  Hence we may assume that
$n=m$.

Recalling (\ref{Pdef}), we  express $P$ in the form
\begin{equation} \label{POLY} P = \{x: \langle x, u_j
\rangle+\la_j \geq 0,\;\; j = 1, \dots, d\}\;,
\end{equation} where the vectors $u_j$ are the inward-pointing normal vectors
to the facets of $P$. The normals $u_j$ are normalized so that
 $u_j$ is a primitive element of $\Z^m$ (i.e., $ru_j\in\Z^m \iff r\in\Z$).

We shall use the  formula given independently by Brion and Vergne
\cite[Theorem 3.12]{BV} and by Guillemin \cite[Theorem 4.1]{Gu2} :
 \begin{equation} \chi_P(e^w) = \todd(\fcal, \partial/\partial h)
 \left.\left(\int_{P(h)} e^{\langle w,x \rangle} dx\right)\right|_{h =
0}\;, \qquad w\in\C^m\;,\ \|w\|< \ep_\fcal.
\label{character}\end{equation}
Here,
\begin{equation}\label{P(h)} P(h) = \{x: \langle u_j, x\rangle + \lambda_j + h_j \geq
0,\ \ 1\leq j
\leq d\}\;,\end{equation} $\fcal=\fcal_P$ is the fan associated to $P$, and  $\ep_\fcal$
is a positive constant. (E.g., if $P$ is Delzant, then $\ep_\fcal =2\pi$.) Also,
$\todd(\fcal,
\partial/\partial
 h)$ is the (generalized) Todd operator constructed by Brion and Vergne
\cite[pp.~374--376]{BV} as follows:

For each cone $C_F\in\fcal$, we consider the open $k$-parallelogram
$$U_F:=\textstyle\left\{\sum_{j\in \jcal(F)} t_j u_j : 0<t_j<1\right\}\;,$$
where $k=\codim F=\dim C_F$.  (Recall (\ref{J}).)  We let $U$
denote the (disjoint) union of the $U_F$, where $F$ runs over the
faces of $P$ (including the open face $P^\circ$, where
$U_{P^\circ}=\{0\}$). Equivalently, $U$ is a `bouquet of
$m$-parallelograms':
$$U=\bigcup_{{\rm vertices }\ v}
\textstyle\left\{\sum_{j\in \jcal(v)} t_j u_j : 0\le
t_j<1\right\}\;,$$ and is an open neighborhood of $0$. We consider
the finite set
$$\Gamma_{\fcal}:= U\cap \Z^m\;.$$
(It is easy to see that $\Gamma_{\fcal}=\{0\}$ if and only if the
polytope $P$ is Delzant.)  For each lattice point $\ga\in
\Gamma_{\fcal}$, we define a unique vector
$g(\ga)=(g_1(\ga),\dots,g_m(\ga))\in \big([0,1)\cap\Q\big)^m$ as
follows:  Let $F$ be the  face such that $\ga\in U_F$; then
$$\ga=\sum_{j=1}^m g_j(\ga) u_j,\qquad g_j(\ga)\ne 0\iff j\in\jcal(F)\;.$$
We note that the map $\ga\mapsto g(\ga)$ is injective. The Todd
operator is then given by:
 \begin{equation} \todd(\fcal, \partial/\partial h) = \sum_{\gamma
 \in \Gamma_{\fcal}}   \left[\prod_{j = 1}^d
 \todd (e^{i2\pi
g_j(\ga)}, \partial/\partial h_j)\right] ,\quad \todd(a, z) =
\frac{z}{1 - a e^{-z}}\label{toddf}\;.\end{equation}

As noted in \cite{BV} (Remark 2.16), if $P$ is a Delzant polytope,
then $\todd(\fcal,
\partial/\partial h)$ is the usual Todd operator
\begin{equation} \todd(\partial/\partial h) = \Pi_{j = 1}^d
\todd(\partial/\partial h_j),\;\; \todd (z) = \todd(1, z) =z (1 -
e^{-z})^{-1}. \end{equation} and formula (\ref{character}) is due
to Khovanskii and Pukhlikov \cite{KP} for that case. (Formula
(\ref{character}) has a generalization to arbitrary convex rational polytopes;
see \cite{BV2}.)  Since the term in (\ref{toddf}) corresponding to
$\ga=0$ is the usual Todd operator and the other terms  vanish at
the origin, the constant term  in the series expansion of
$\todd(\fcal,\d/\d h)$ equals 1 for all simple polytopes.
Furthermore, the coefficients of  the expansion of
$\todd(\fcal,\d/\d h)$ are real since if we regard the $g(\ga)$ as
elements of $\Q^m/\Z^m$, then the set $\{g(\ga)\}$ is a union of
subgroups and hence is invariant under multiplication by $-1$.

Let us check how formula (\ref{character}) dilates. The fan does
not change when the polytope is dilated, so the Todd operator is
the same for $P$ and for $N P$. Noting that
\begin{equation} (N P) (Nh) = \{x: \langle u_j, x\rangle + N
\lambda_j + N h_j \geq 0\} = N (P(h))\;,
\end{equation} we change variables $h \mapsto N h,\ x\mapsto Nx$ in
(\ref{character}) to obtain
\begin{equation} \chi_{N P}(e^{w}) =  N^m  \;  \todd(\fcal, N^{-1}
\partial/\partial h)  \left.\left(\int_{P(h)} e^{ N \langle w, x
\rangle} dx\right)\right|_{h = 0}, \qquad \|w\|< \ep_\fcal\;.
\label{charN} \end{equation}

\medskip To prove Proposition \ref{newchar}, we use the following
lemma to change variables so that the integral is over the fixed
polytope $P$ while the integrand depends on $h$.

We say that a map $f:P\to\R^m$ is a {\it regular rational map\/}
if it can be written in the form $f=(p_1/q_1,\dots,p_m/q_m)$,
where the $p_j$ and $q_j$ are polynomials on $\R^m$ and the $q_j$
are nonvanishing on $P$.

\begin{lem}\label{Ph}Let $P$ be a simple polytope in $\R^m$ with nonempty interior.
Then there exist regular rational maps $f_j:P\to\R^m$, $1\le j\le
d$, such that for $h\in\R^d$ sufficiently small, the map
$\ga_h:P\to\R^m$ given by
$$\ga_h(x)=x+\sum_{j=1}^dh_jf_j(x)$$
is a diffeomorphism from $P$ onto $P(h)$ mapping each face of $P$
onto the corresponding face of $P(h)$. \end{lem}

\begin{proof} Recall that the assumption that $P$ is simple means that there are precisely $m$ hyperplanes $\{x: \langle x, u_j
\rangle+\la_j =0\}$ through each vertex, and hence for small $h$,
the faces of $P(h)$ correspond to the faces of $P$.

We first construct an `algebraic partition of unity' $\{\psi_k\}$
over $P$ as follows. Let $\{v^1,\dots,v^s\}$ denote the set of
vertices of $P$.  Recalling (\ref{Pdef}), we let $\wt\psi_k:P\to
[0,+\infty)$ be given by
$$\wt\psi_k= \prod\{\ell_j: {j\not\in\jcal(v^k)}\}\;,\qquad k=1,\dots,s\;.$$
We then let
$$\psi_k= \frac{\wt\psi_k}{\sum_{j=1}^s \wt\psi_j}\;.$$
We easily see that $0\le \psi_k \le 1$, $\sum_{k=1}^s  \psi_k = 1$
and that $\psi_k\inv(0)$ is the union of the closed facets of $P$
that do not contain $v^k$.

The vertices $v^k(h)$ of $P(h)$ are of the form
$$v^k(h)=v^k+ L(k)h\;,\qquad1\le k\le s,\ \|h\|<\ep\;,$$ where $L(k):\R^d\to\R^m$ is a
linear transformation.  (In fact, it is easy to see that $L(k)$ is
given by first projecting onto the $m$ coordinates corresponding
to the facets incident to $v_k$, and then applying an element of
GL$(m,\R)$.) Let
$$w^{kj} = L(k)E^j\in\R^m, \qquad 1\le k\le s,\ 1\le j\le
d\;,$$ where $\{E^j\}$ is the standard basis of $\R^d$,  so that
\begin{equation}\label{wkj}v^k(h)=v^k+\sum_{j=1}^d h_jw^{kj}\;,\qquad \|h\|<\ep\;.
\end{equation}

For each $k$, let $v^{\nu(k,1)},\dots v^{\nu(k,m)}$ be the
vertices connected to $v^k$ by an edge (1-dimensional face) of
$P$. For $1\le k\le s,\ 1\le j\le d$, let $f_{kj}:\R^m\to\R^m$ be
the unique affine map such that
$$f_{kj}(v^l)=w^{lj} \quad \mbox{for }\ l=k,\ \nu(k,1),\ \dots,\ \nu(k,m)\;. $$
So by (\ref{wkj}), if $l=k$ or  if $v^l$ is a vertex connected to
$v^k$ by an edge of $P$, then
\begin{equation}v^l(h)=v^l+ \sum _{j=1}^dh_j f_{kj}(v^l)\;,\qquad  \mbox{for }\ l=k,\
\nu(k,1),\ \dots,\ \nu(k,m)\;,\quad
\|h\|<\ep\;.\label{fkj}\end{equation}

We claim that the maps
$$f_j:= \sum_{k=1}^s \psi_k f_{kj}$$ satisfy the conclusions of the lemma. To verify
the claim, we let
$$\ga^k_h(x) = x+\sum_{j=1}^dh_jf_{kj}(x)\;,\qquad 1\le k\le s\;,$$
so that
$$\ga_h=\sum\psi_k\ga^k_h\;.$$ It follows from (\ref{fkj}) that when $x$ is a vertex we
have
$$\ga_h(v^k)=\ga_h^k(v^k)=v^k(h),$$ since $\psi_l(v^k)=0$ when $k\ne l$.

Next we consider points on the 1-skeleton of $P$. Suppose that
$x=tv^1+(1-t)v^2$, where $0< t < 1$ and $v^1,v^2$ are joined by an
edge of $P$. Then $\psi_k(x)=0$ for $k>2$ and we again conclude
from (\ref{fkj}) that $\ga_h^k(v^1)= v^1(h)$ and $\ga_h^k(v^2)=
v^2(h)$ for $k=1,2$. Since the $\ga_h^k$ are affine maps, it follows that
$$\ga_h(x)= tv^1(h)+(1-t)v^2(h)$$ so that $\ga_h$ takes each edge of $P$ linearly onto
the corresponding edge of $P(h)$.

In general we do not have linearity on faces of dimension $>1$, so
we need to modify the argument for the higher skeletons. Suppose
now that $x$ lies in a 2-dimensional face $F$ with vertices
$v^1,\dots,v^r$ ($r\ge 3$) and  $x=t_1v^1+t_2v^2+t_3v^3$, where
$\min t_j\ge 0$ and $\sum t_j=1$.  Let $1\le l\le r$. Since the
points $v^1(h),\dots,v^r(h)$ lie on a plane $S_h$ and the
$\ga_h^k$ are linear, it follows from the above argument that
$\ga_h^k(v^l)\in S_h$ for $1\le k\le r$. (For example, if
$v^2,v^3$ are each connected to $v^1$ by an edge, then
$\ga^1_h(v^l)=v^l(h)$ for $l=1,2,3$, but for $3<l\le r$, we can
conclude only that $\ga^1_h(v^l)$ lies on the plane through
$v^1(h),v^2(h),v^3(h)$.)  Since for $k>r$, there exists a closed
facet containing $v_1,v_2,v_3$ but not containing $v_k$, and hence
$\psi_k(x)=0$.  Therefore,
$$\ga_h(x)= \sum_{1\le j,k\le r} t_j\psi_k(x)\ga_h^k(v^j)\in S_h\;.$$ Furthermore,
for $h$ small, $D\ga_h \approx D\ga_0 \equiv I$, so $\ga_h|_F:F\to
S_h$ is a local homeomorphism.  Therefore, $\d_{S_h}
\big(\ga_h(F)\big)=\ga_h(\d F)$, which is the boundary of the
convex polygon $F(h)$ with vertices $v^1(h),...,v^r(h)$. It
follows that $\ga_h(F)=F(h)$ and hence $\ga_h|_{\overline
F}:\overline F \to \overline {F(h)}$ is a (global) homeomorphism for sufficiently small
$h$.

Continuing by induction on $\dim F$, we conclude that
$\ga_h:\overline F\approx \overline {F(h)}$ for all faces $F$ of $P$, and in particular,
$\ga_h$ maps $P$ diffeomorphically onto $P(h)$, for $h$ small.
\end{proof}

\noindent{\it Completion of the proof of Proposition
\ref{newchar}:\/} Let $\ga_h:P\approx P(h)$ be as in Lemma
\ref{Ph}.  Making the change of variables $x \mapsto \ga_h(x)$, we
have
\begin{equation} \int_{P(h)}
e^{N{\langle w, x \rangle}}\, dx =\int_P e^{N\langle w, \ga_h(x)
\rangle}\det D\ga_h(x)\, dx =  \int_P e^{N[\langle w, x \rangle
+\Si_j\langle w, f_j(x)\rangle h_j]}G(x,h)\, dx \;, \label{change}
\end{equation} where $G$ is a regular rational function on $P\times\R$ of the form
$$G(x,h)=\sum_{\al\in\N^d,|\al|\le m} G_\al(x)h^\al\;,\qquad
 G_0\equiv 1\;.$$

Hence by (\ref{charN}) and (\ref{change}), we have
\begin{equation} \chi_{N P}(e^{w}) =  N^m  \int_P  \todd(\fcal, N^{-1}
\partial/\partial h) \left.\left(e^{N[\langle w, x \rangle +\Si_j\langle w,
f_j(x)\rangle h_j]}G(x,h)\right)\right|_{h=0}\, dx \;,\quad
\|w\|<\ep\;. \label{change2}\end{equation}

To justify interchanging integration and  Todd differentiation in
the above and to simplify (\ref{change2}), we use the following
identity from  \cite[Lemma 1]{KP}): Suppose that $G(h)\in
\C[h_1,\dots,h_k]$  is a polynomial and $F$ is a convergent power
series about $0\in\C^k$ with domain of convergence $\Omega$. Then
\begin{equation} F\left(\partial/\partial h\right)\left[G(h) e^{\Si_{j = 1}^k
q_j h_j}\right] = G \left(\partial/\partial q\right) \left[F(q)
e^{\Si_{j = 1}^k q_j h_j} \right],\quad q\in\Omega,\ h\in\C^k.
\label{switch}\end{equation} (Equation (\ref{switch}) is easily
verified by noting that both sides equal
$F\left(\partial/\partial h\right)G
\left(\partial/\partial q\right) \exp({\Si_{j = 1}^k q_j h_j})$.)

We apply (\ref{switch}) with $F(q)=\todd(\fcal, N\inv q),\
G(q)=G(x,q)$ to the integrand of (\ref{change2}):
\begin{eqnarray} &&\left.\todd(\fcal, N^{-1}
\partial/\partial h)\left[
e^{N[\langle w, x \rangle +\Si_j\langle w, f_j(x)\rangle
h_j]}G(x,h)\right]\right|_{h = 0} \nonumber \\&&\hspace{2cm} =\
e^{N\langle w,x\rangle} \left.\todd(\fcal, N^{-1}
\partial/\partial h)\left[G(x,h)
e^{\Si_{j = 1}^d q_j h_j}\right]\right|_{h = 0,\ q_j=N
\langle w, f_j(x)\rangle}\nonumber \\
&&\hspace{2cm} =\ e^{N\langle w,x\rangle}G(x,\d/\d q)
\left.\left[\todd(\fcal, N^{-1}q) e^{\Si_{j = 1}^d q_j
h_j}\right]\right|_{h = 0,\ q_j=N
\langle w, f_j(x)\rangle}\nonumber \\
&&\hspace{2cm} =\ e^{N\langle w,x\rangle}G(x,\d/\d q)
\left.\todd(\fcal, N^{-1} q)\right|_{q_j=N
\langle w, f_j(x)\rangle}\nonumber \\
&&\hspace{2cm} =\ e^{N\langle w,x\rangle}G(x,N\inv\d/\d \tilde q)
\left.\todd(\fcal, \tilde q)\right|_{\tilde q_j= \langle w,
f_j(x)\rangle}\label{amp0}\;.\end{eqnarray} In particular, if we
let $\todd_k$ denote the sum of the terms of the Todd series of
degree $\le k$, then
$$\begin{array}{l}\left.\todd_k(\fcal, N^{-1}
\partial/\partial h)\left[
e^{N[\langle w, x \rangle +\Si_j\langle w, f_j(x)\rangle
h_j]}G(x,h)\right]\right|_{h = 0} \\[10pt] \hspace{3cm} =\ e^{N\langle
w,x\rangle}G(x,N\inv\d/\d \tilde q) \left.\todd_k(\fcal, \tilde
q)\right|_{\tilde q_j= \langle w,
f_j(x)\rangle}\\[10pt] \hspace{3cm} \to
\left.\todd(\fcal, N^{-1}
\partial/\partial h)\left[
e^{N[\langle w, x \rangle +\Si_j\langle w, f_j(x)\rangle
h_j]}G(x,h)\right]\right|_{h = 0}\end{array}$$ uniformly for $x\in
P$ and $w$ sufficiently small, which justifies the interchange of
integration and Todd differentiation in (\ref{change2}).

Let $\La(x)$ denote the $d\times m$ matrix whose rows are the
$f_j$; i.e., $\La(x):\R^m\to\R^d$ is the linear map given by
$(\La(x)y)_j=\langle y, f_j(x)\rangle$.   By (\ref{amp0})
\begin{equation}\label{amp} \left.\todd(\fcal, N^{-1}
\partial/\partial h)\left[
e^{N[\langle w, x \rangle +\Si_j\langle w, f_j(x)\rangle
h_j]}G(x,h)\right]\right|_{h = 0} = e^{N\langle w, x
\rangle}\sum_{l=0}^m N^{-l} T_l(x,\La(x)w) \;,\end{equation} where
$T_l(x,q)$ is a convergent power series in $q$ with coefficients
that are regular rational functions on $P$.  More specifically,
$$T_0(x,q)=\todd(\fcal,q)\;,\qquad T_l(x,q)=G_l(x,\d/\d q)\todd(\fcal, q) \quad (1\le
l\le m)\;,$$ where $G_l(x,\d/\d q)$ is a homogeneous differential
operator of order $l$ whose coefficients are regular rational
functions of $x\in P$.

It follows from (\ref{change2}) and (\ref{amp}) that the integral
formula of Proposition~\ref{newchar} holds for small $w$, with
\begin{equation}\label{Al}A_l(x,w)=T_l(x,\La(x)w)\;.\end{equation}  Moreover, from
the construction of the $T_l$, we see that the integrand has no
poles if $\|\La(x)\Im w\| <\ep_\fcal$ for all $x\in P$. Thus by
analytic continuation, the identity holds on $S(\ep)$, where $\ep=
\ep_\fcal/\sup_{x\in P}\|\La(x)\|$.

Statements (i) and (ii) follow from (\ref{Al}) and the fact that
the coefficients of the Todd function are real and the constant
term is 1.

To verify (iii), we fix $x\in P$ and $w\in C_x$, and we let  $y_j=
\langle w, f_j(x)\rangle$, for $j=1,\dots,d$.

{\medskip\noindent{\it Claim 1:\/} }
\begin{enumerate}
\item[i)] $y_j \ge 0$, for $j=1,\dots,d$;
\item[ii)]  $y_j = 0$, for $j\not\in\jcal(x)$.
\end{enumerate}

\medskip\noindent{\it Proof of Claim 1:\/}
(i)  Let $F$ be the face containing $x$. Then  $x- t f_j(x)$ is on
the face $F(h)$ corresponding to $F$  of the polytope $P(h)$,
where $h_j= -t$ (for $t>0$ sufficiently small) and $h_k=0$ for
$k\ne j$.  In particular, $x- t f_j(x) \in P$. Since $w$ is in the
normal cone to $F$, we have $\langle w, x- t f_j(x) \rangle\le
\langle w,  x\rangle$; i.e., $y_j=\langle w, f_j(x) \rangle\ge 0$.

(ii) If $j\not\in\jcal(x)$, then $F(h)\subset F$ for $h$ as above, and hence
$f_j(x)$ is tangent to $F$, so that $y_j=\langle w, f_j(x)
\rangle= 0$. \qed

\medskip\noindent{\it Claim 2:\/}  Let $G_F= \Gamma_{\fcal}\cap \overline {U_F}$,
where  $F$ is the face containing $x$. Then
$$ A_0(x,w) = \sum_{\gamma
 \in G_F} \; \prod_{j\in\jcal(F)}
 \todd (e^{i2\pi
g_j(\ga)}, y_j)\;.$$

\medskip\noindent{\it Proof of Claim 2:\/} We have
$$A_0(x,w)=\todd(\fcal,y)\;,\quad y=(y_1,\dots,y_d)\;,\quad y_j= \langle w,
f_j(x)\rangle\;. $$   We first note that
$$G_F=\{\ga\in\Z^m: \ga=\sum_{j\in\jcal(F)}t_ju_j,\ 0\le t_j <1\}\;.$$
Thus, if $\ga\in\Ga_\fcal\sm G_F$, then there exists $j\not\in
\jcal(F)$ such that $g_j(\ga)\ne 0$.  But by Claim 1, the
corresponding $y_j =0$.  Hence the term corresponding to $\ga$ in
the sum (\ref{toddf}) for $\todd(\fcal,y)$ vanishes.

Now suppose that $\ga\in\Ga_\fcal\cap G_F$.  Then by definition,
$g_j(\ga)=0$, for all $j\not\in\jcal(F)$. Therefore $$\todd
(e^{i2\pi g_j(\ga)}, y_j)=\todd(1,0) =1\qquad \mbox{for }\ j\not\in \jcal(F)\;.$$
Hence the sum  (\ref{toddf}) for $\todd(\fcal,y)$ reduces to the
sum in Claim 2.\qed

\bigskip We easily see that the set
$$Q_F:=\{(g_j(\ga))_{j\in\jcal(F)}\in \Q^n/\Z^n:\ga\in G_F\}$$
is a subgroup of  $\Q^n/\Z^n$, where $n=\# \jcal(F)=\codim F$.
The final conclusion (iii) of the Proposition is then a
consequence of Claims 1--2 and the following elementary lemma.

\begin{lem} \label{lastlemma} Let $\Ga$ be a finite subgroup of ${\mathbf T}^n$ and let
$y_j\ge 0$, for $1\le j \le n$.  Then
$$ \sum_{a\in \Ga} \; \prod_{j=1}^n \todd (a_j, y_j) >0\;,\qquad a=(a_1,\dots,a_n)\;.$$
\end{lem}
\begin{proof}
We can assume without loss of generality that the $y_j$ are
positive.  For suppose that $y_n=0$, for example.  Then
$\todd(a_n,y_n)=1$ if $a_n=1$, and $\todd(a_n,y_n)=0$ if $a_n\ne
1$. Then the sum reduces to the corresponding sum over the group
$\Ga':=\{a'\in \mathbf T^{n-1}:(a',1)\in \Ga\}$.

Let $r_j=e^{-y_j}<1$.  Then  $$S:=\frac 1 {y_1\cdots y_n}
\sum_{a\in \Ga} \; \prod_{j=1}^n\; \todd (a_j, y_j) = \sum _{a\in
\Ga}\;\prod_{j=1}^n (1-a_jr_j)\inv = \sum _{a\in
\Ga}\;\sum_{\be\in\N^n} a^\be r^\be = \sum_{\be\in\N^n} c_\be
r^\be\;,$$ where $$c_\be= \sum _{a\in \Ga} a^\be\;.$$

We claim that $c_\be$ is either 0 or $o(\Ga)$.  Indeed, if
$a^\be=1$ for all $a\in\Ga$, then $c_\be = o(\Ga)$.  On the other
hand, if there exists $b\in\Ga$ with $b^\be \ne 1$, then $c_\be =
\sum _{a\in \Ga} a^\be =  \sum _{a\in \Ga} (b\cdot a)^\be = b^\be
c_\be$, and hence $c_\be = 0$.  Since $c_0=o(\Ga)$, it then
follows that $S>0$.
\end{proof}

\section{Mass asymptotics: proof of Theorem~\ref{MASS}} \label{s-mass}

In this section, we prove a precise asymptotic formula for the
conditional \szego kernel on the diagonal
(Theorem~\ref{SZEGO}), which yields the mass asymptotics of
Theorem~\ref{MASS}.

First, we verify formula (\ref{2measures}) for the expected mass
density.
 Let $P$ denote a convex integral polytope in $\R^m$. Recalling the definition
(\ref{CG}) of the conditional probability measure $\gamma_{p|P}$
on the space $\poly(P)$ of polynomials with Newton polytope $P$,
we see that the expected value of the mass density with respect to
$\gamma_p|_P$ is given by:
\begin{equation*}
\E_{|P}\left(|f(z)|^2_\FS\right) = \sum_{\alpha,\be \in P}
\frac{\E(\la_{\alpha}\bar \la_\be
)\chi_{\alpha}(z)\overline{\chi_{\be}(z)}}
{\|\chi_{\alpha}\|\|\chi_{\be}\| (1+\|z\|^2)^{p/2}}
\;.\end{equation*} Since the $\la_\al$ are independent complex
random variables with variance 1 (i.e., $\E_{|P}(\la_{\alpha}\bar
\la_\be)=\delta_\al^\be$), we have by (\ref{Sz2}):
\begin{equation}\label{Eszego}
\E_{|P}\left(|f(z)|^2_\FS\right) = \sum_{\alpha \in P} \frac{
|\chi_\al(z)|_\FS^2}{||\chi_{\alpha}||^2}=\Pi_{|P}(z,z)\;.\end{equation}
It then follows by expressing the Gaussian in spherical
coordinates that $$\E_{\nu_{P}}(|f(z)|^2_\FS)  =\frac{1}{\#
P}\E_{|P}\left(|f(z)|^2_\FS\right)= \frac{1}{\#
P}\Pi_{|P}(z,z)\;.$$ Replacing $P$ with $NP$, we obtain
\begin{equation}\label{Eszego2}\E_{\nu_{NP}}(|f(z)|^2_\FS)
=\frac{1}{\#(NP)}\Pi_{|NP}(z,z)\;.\end{equation} Recall that by
(\ref{RR}), the number $\#(NP)$ of lattice points in the polytope
$NP$ is a polynomial in $N$ of  degree equal to the dimension of
$P$. The mass asymptotics of Theorem~\ref{MASS} is an immediate
consequence of (\ref{Eszego2}) and the asymptotic expansion of the
conditional \szego kernel on the diagonal given in
Theorem~\ref{SZEGO} below.

\begin{theo}\label{SZEGO}   Suppose that $P$ is a convex integral polytope in
$\R^m$ such that $P\subset p\Sigma$ and  $P\not\subset \d (p\Si)$.
Then:
\begin{enumerate}
\item[i)] If $P^\circ\ne \emptyset$,
then for $z$ in the classically allowed region $\acal_P$, we have
$$\Pi_{|NP}(z,z) = \prod_{j=1}^m (Np+j) + R_N(z)
\;,\qquad \|R_N\|_{\ccal^k(K)}=O(e^{-\la_K N}) \quad \mbox{for }\ k=1,2,3,\dots,
$$ for all compact $K\subset \acal_P$, where $\la_K>0$.
\item[ii)] On each open forbidden region $\rcal_F^\circ$,
$$\Pi_{|NP}(z,z)= N^{\frac{m+r}{2}} e^{-N b_P(z)}\big[c_0^F(z) +
c_1^F(z)N\inv+\cdots + c_k^F(z)N^{-k} +R_k^F(z)\big],$$ where
$r=\dim F$ and
\begin{enumerate}
\item \ $c_j^F\in\ccal^\infty(\rcal_F^\circ)$ and $c_0^F>0$ on $\rcal_F^\circ$;
\item \ $\|R_k^F\|_{\ccal^j(K)}=O(N^{-k-1})$, for all compact $K\subset
\rcal_F^\circ$ and for all $j,k$;
\item \ $b_P>0$ on
$(\C^*)^m \sm \overline{\acal_P}$;
\item \ $b_P$ is given by formula (\ref{b});
\item \  $b_P\in\ccal^1_\R((\C^*)^m)$  (with  $b_P=0$ on
${\acal_P}$), and $b_P$ is $\ccal^\infty$ on each closed region
$\overline{\rcal_F}$.
\end{enumerate}
\end{enumerate}
\end{theo}

Note that $b_P$ fails to be $\ccal^2$ at transition points, since
the limit expected zero current $\psi_P$ of Theorem~\ref{main} is
discontinuous at the transition points and is given in terms of
the second derivatives of $b_P$  by formula (\ref{psiP}).

\begin{rem}
For the case  where  $P=\{\be\}$ is a single lattice point in
$p\Si^\circ$, the allowed region $\acal_{\{\be\}}$ is all of $(\C^*)^m$. Since the point
$q(z)$ given by  (\ref{cond1}) lies in $\d P$, in this case $q(z)=\be$ for all $z\in
(\C^*)^m$.  Recalling (\ref{mcal})--(\ref{interpretb}), we then have
\begin{equation}\label{bpoint} b_{\{\be\}}(z) =-2\log
{\mcal_{\be}(z)}=\textstyle\log
\left|\wh\chi^p_\be\big(\mu\inv(\frac 1p \be)\big)\right|^2 -\log
\left|\wh\chi^p_\be(z)\right|^2\;.\end{equation} Theorem~\ref{SZEGO} then says that
\begin{equation}\label{Pidecay} \Pi_{|N\{\be\}}(z,z) = \frac 1
{\|\chi_\be\|^2}|\chi_\be(z)|^2_\FS  = N^{\frac m2} (c_0+c_1 N\inv +c_2 N^{-2}+ \cdots)
\,\mcal_{\be}(z)^{2N} \;.\end{equation}
(The asymptotic formula (\ref{Pidecay}) is also an immediate consequence of
formula (\ref{IP2}) for the $\lcal^2$ norms of monomials and Stirling's formula.
A more precise formula and
general results on the decay of monomials normalized over toric
varieties are given in \cite{STZ2}.)

Replacing $\be$ with an arbitrary (non-lattice) point $x\in P^\circ$ in (\ref{bpoint}),
we define
\begin{equation}\label{bx1}b_{\{x\}}:=-2\log \mcal_{x}(z)\;.\end{equation}
Recalling (\ref{interpretb}), we see that for a convex lattice polytope $P$, we have
$$b_P(z)=b_{q(z)}(z)\;,\qquad \Pi_{|NP}(z,z)= N^{\frac{m+r}{2}} e^{-N
b_{q(z)}(z)}\big[c_0^F(z) + c_1^F(z)N\inv+\cdots  \big]\;.$$
To compute $b_{\{x\}}$, we use (\ref{b}) with $q(z)=x$ and
$\mu(e^{-\tau/2}\cdot z)= \frac 1p x$, which yields
$$\tau_j= \log |z_j|^2 -\log x_j+\log (p-\textstyle\sum_{j=1}^m x_j)\;,$$
and hence
\begin{equation}\label{bx2} b_{\{x\}}(z)= \sum_{j=0}^m x_j\log \frac{x_j}{p}
-\log\frac{|z|^{2x}}{(1+\|z\|^2)^p}\qquad
(x_0=p-\textstyle\sum_{j=1}^m x_j)\;.\end{equation}
\end{rem}

\medskip The asymptotics of Theorem~\ref{SZEGO} are uniform away from the
transition points only.  To take care of the transition points, we
shall also prove the following local uniform convergence result on
all of $(\C^*)^m$:

\begin{prop}\label{convergence}
Let $P$ be as in Theorem~\ref{SZEGO}. Then
$$\frac{1}{N} \log \Pi_{|NP}(z,z)
\to -b_P(z) $$ uniformly on all compact subsets of $(\C^*)^m$.
\end{prop}

Proposition \ref{convergence} is the main ingredient in our proofs of Theorems
\ref{main}--\ref{simultaneous} on the asymptotics of  zeros in \S \ref{DZ}.

The remainder of this section is devoted to the proofs of Theorem~\ref{SZEGO}
and Proposition~\ref{convergence}. To motivate the argument, we begin by using
our integral formula for the
 polytope character (Proposition~\ref{newchar}) to show that the integral is rapidly
decaying when $z$ is in the classically forbidden region, while
for $z$ in the allowed region, the phase is of positive type with
a nondegenerate critical point and hence $\Pi_{|N P}(z, z)$ has an
asymptotic expansion in powers of $\frac 1N$. To obtain our
expansion in the forbidden region,
 we use Proposition~\ref{newchar} to give  an oscillatory integral formula
(\ref{lattice3})--(\ref{amplitude}) for  the conditional \szego
kernel $\Pi_{|N P}(z, z)$.  We then seek deformations of the
contour of integration so that the (analytic continuation of) the
phase picks up a critical point.   In \S \ref{s-critical}, we
formulate necessary and sufficient geometric conditions for the
existence of a critical point where the phase has maximal real
part, and in \S \ref{s-normalbundle} we show that there is a
contour for which these conditions are satisfied. Assuming that $P$ is simple, we
obtain in \S \ref{s-asymptotic} the asymptotic expansion (ii) by performing a
stationary phase analysis over the polytope when the critical
points of the phase lie on the boundary. In \S \ref{s-decay}, we show that the
decay function satisfies properties (c) and (e).  In \S \ref{s-precise}, we verify the
precise asymptotic formula (i), and in \S \ref{s-general} we
complete the proof of part (ii) of Theorem~\ref{SZEGO} by extending our results
to non-simple polytopes. We prove
Proposition~\ref{convergence} in \S \ref{s-convergence}.

We start our proof of Theorem~\ref{SZEGO} with a simple integral
formula for the conditional \szego kernel:
\begin{equation}\label{lattice}\Pi_{|N P} (z,z) = \frac{1}{(2\pi )^m}
\int_{\T } \Pi_{Np} (z,e^{i\phi}\cdot z)
{\chi_{NP}(e^{i\phi})}\,d\phi\;,\end{equation} where we recall
that
$$\T=\{(\zeta_1,\dots\,\zeta_m)\in (\C^*)^m: |\zeta_j|=1,\ 1\le j\le m\}.$$
The identity (\ref{lattice}) follows immediately from the explicit
formulas (\ref{Sz}) and (\ref{Sz2e}) for the \szego kernels and
the definition (\ref{char}) of the character $\chi_{NP}$.

Substituting formula (\ref{Sz}) for $\Pi_{Np} (z,e^{i\phi}\cdot
z)$ in (\ref{lattice}), we obtain:
\begin{equation}\Pi_{|NP}(z,
z) =\frac{(Np+m)!}{(Np)!(2\pi)^m}\int_{\T} \left[\frac{
1+\sum_{j=1}^m e^{-i\phi_j}|z_j|^2}{1+\|z\|^2}\right]^{Np}
{\chi_{NP}(e^{i\phi})}\,d\phi\;.\label{lattice*}\end {equation}

{\bf We now  assume that $P$ is a simple polytope.}  In \S
\ref{s-general}, we shall derive the general case as a corollary
of our results for simple polytopes.

\subsection{The classically allowed region} In order to
motivate our argument, we give a short proof of part (i) of
Theorem~\ref{SZEGO} with the exponentially small error term
replaced with an $O(N^{-\infty})$ error.  The actual proof of part
(i) will be given in \S \ref{s-precise} and does not depend on the
results of this section.

Assume for simplicity that $P^\circ \neq \emptyset$. We first
choose a cut-off function $\xi\in \ccal^\infty(\R^m)$ vanishing
outside the $\ep_\fcal$-ball and with $\xi\equiv 1$ on the
$(\ep_\fcal/2)$-ball.  Since $\Pi_{Np}(z,e^{i\phi}\cdot z)$ decays
exponentially (in $N$) when $\|\phi\|\ge \ep_\fcal$, it follows
from Proposition \ref{newchar} and (\ref{lattice*}) that
\begin{equation} \Pi_{|NP}(z, z) \sim N^{2m}\frac{\sigma(N)}{(2\pi
)^{m}} \;   \int_{{\bf T}^m}\int_{P}\xi(\phi)
e^{N\Psi_{\!\acal}(\phi, x;
z)}\big[A_0(x,\phi)+A_1(x,\phi)N\inv+\cdots\big]\, dx\, d\phi\;,
\label{mass}
\end{equation} where $\sigma(N)= (p+\frac
1N)\cdots(p+\frac mN)= p^m +\sigma_1 N\inv +\cdots+ \sigma_M
N^{-m}$, $A_0(x,0)=1$ and the phase $\Psi_{\!\acal}$ is given by
\begin{equation} \Psi_{\!\acal}(\phi, x; z) = i  \langle \phi, x
\rangle +p\log\left(\frac{ 1+\sum_{j=1}^m
e^{-i\phi_j}|z_j|^2}{1+\|z\|^2}\right).\label{phase}
\end{equation} Here, ``$\sim$'' means modulo $O(N^{-\infty})$.

The integral in (\ref{mass}) is an `oscillatory integral' with `complex phase'
$\Psi_\acal$; to evaluate it by the method of stationary phase \cite[Ch.~7]{H}, we must
find the critical points of
$\Psi_{\!\acal}$. The (interior) critical point equations
$$d_x\Psi_{\!\acal}=i\phi = 0,\quad d_\phi\Psi_{\!\acal}|_{\phi=0} =ix
-\frac{ip}{1+\|z\|^2}(|z_1|^2,\dots,|z_m|^2) =ix-ip\mu_{\Sigma}(z)
= 0$$ yield \begin{equation}\label{allowedcrit} \phi=0,\qquad
x=p\mu(z)\;.\end{equation} We note that $\Re\Psi_{\!\acal} \le 0$
and $\Psi_{\!\acal}=0$ at the critical point.

The Hessian $\hcal \Psi_{\!\acal} $ of $\Psi_{\!\acal}$ (with
respect to the variables $\phi,x$) is of the form
\begin{equation}\label{Hessian}\hcal\Psi_{\!\acal}|_{(0,x)} =
\left(\begin{array}{ll} {\mathbf C} & i\I
\\ &
\\ i\I & 0
\end{array}
\right)\end{equation} where $\I$ denotes the $m\times m$ identity
matrix.  The matrix ${\mathbf C}$ is given by
$$C_{jk} = \left. \frac{\d^2\Psi_{\!\acal}}{\d\phi_j\d\phi_k}\right|_{(0,x)} =p(
-I_j\delta_j^k +I_jI_k)\;,$$ where
$$I_j=\frac{|z_j|^2}{1+\|z\|^2}\;;\ \mbox{ i.e., }\
\mu(z)=(I_1,\dots,I_m)\;.$$

Hence \begin{equation}\label{detH} \det
\hcal\Psi_{\!\acal}|_{(0,x)} =1\,\end {equation}
 and the inverse Hessian is given by
$$\hcal\Psi_{\!\acal}|_{(0,x)}\inv = \left(\begin{array}{cc} 0 & -i\I \\ & \\ -i\I
& {\mathbf C}
\end{array} \right).$$ Hence, the  Hessian operator is given by
$$ {H} = -i\sum_j \partial^2/\partial
\phi_j \partial x_j +\sum_{jk}C_{jk}
\partial^2/
\partial x_j \d x_k.$$

By (\ref{mass}) and \cite[Theorem 7.7.5]{H},  we have a stationary
phase expansion of the form
$$\Pi_{|NP}(z,z)
\sim p^m\sigma(N)\sum_{2k\ge 3j\ge 0}\left.\frac{1}{j!k!2^k}\,
N^{m+j-k}\,{H}^kg_3^j\right|_{\phi= 0, x=p\mu(z)}\;,$$
 where $g_3$ is the third order Taylor remainder of the
phase $\Psi_{\!\acal}$.  We note that $g_3$ is a function of
$\phi$ alone, but that every $\phi$ derivative is accompanied by
an $x$-derivative in ${ H}$. It follows that only the very first
term of the summation is non-zero. Thus, the stationary phase
expansion is simply $p^m\sigma(N)N^m$.

\subsection{The classically forbidden region} \label{s-critical}
By (\ref{allowedcrit}), $\Psi$ has a critical point if and only if
$z$ is in the classically allowed region. Thus, in the classically
forbidden region, there are no critical points and the integral
(\ref{mass}) is rapidly decaying.

We now let $z$ be a point in the classically forbidden region.  To
evaluate the integral by the method of stationary phase,  we need
to deform the contour of integration to pick up a  critical point
with maximal real part along the contour. We complexify the real
torus $\T $ to $(\C^*)^m$ with variables $\zeta_j =  \rho_j e^{i
\phi_j}$ so that (\ref{lattice*}) may be written as
\begin{equation}\label{lattice1}\Pi_{|N P} (z,z) = \frac{(Np+m)!}{(Np)!(2\pi i)^m}
\int_{\T } \left(\frac{ 1+\sum_{j=1}^m
\zeta_j\inv|z_j|^2}{1+\|z\|^2}\right)^{Np} {\chi_{NP}({\zeta})}\,
\prod_{j = 1}^m \frac{d \zeta_j}{\zeta_j}\;.\end{equation} Since
the integrand is holomorphic in $\zeta\in (\C^*)^m$, we  can deform
the contours in (\ref{lattice1}) and instead integrate over a
torus of the form
$$\T_\tau:=\{\zeta\in(\C^*)^m:|\zeta_1| = e^{\tau_1}, \cdots,  |\zeta_m| = e^{\tau_m}\}
\;.$$ We shall show below how to choose the parameter
$\tau=(\tau_1,\dots,\tau_m)\in\R^m$ to obtain our required
critical point.

Deforming (\ref{lattice1}) to $\T_\tau$, we have
\begin{equation}\label{lattice2}\Pi_{|N P} (z,z) =
\frac{(Np+m)!}{(Np)!(2\pi)^m} \int_{\T
}e^{N\Psi_0(\phi;\tau,z)}{\chi_{NP}(e^{\tau+i\phi})}\,d\phi\;.
\end{equation} where
\begin{eqnarray}\Psi_0(\phi;\tau, z) =
p\log\left(\frac{ 1+\sum_{j=1}^m
e^{-\tau_j-i\phi_j}|z_j|^2}{1+\|z\|^2}\right)\;.
\label{PHASEcx}\end{eqnarray} Using a $\ccal^\infty$ partition of
unity $\{\xi_1,\xi_2\}$  of $\T$, where $\xi_1(\phi)$ has support
in a small neighborhood of 0 and is identically 1 near 0, we
decompose (\ref{lattice2}) into two integrals:
\begin{equation}\label{lattice3}\Pi_{|N P} (z,z) =
\frac{(Np+m)!}{(Np)!}(\ical_N^1+\ical_N^2)\;,\qquad \ical_N^j:=
\frac{1}{(2\pi )^m} \int_{\T } \xi_j(\phi)
e^{N\Psi_0(\phi;\tau,z)}{\chi_{NP}(e^{\tau+i\phi})}\,d\phi\;.\end{equation}
We shall show that $\ical_N^1$ has the desired asymptotic
expansion (when $z$ is not a transition point) and that
$\ical_N^2=O(e^{-\delta}\ical_N^1)$ for some $\delta >0$.

Let $\dim P=n$. Under our assumption that $P$ is simple,
Proposition \ref{newchar} lets us write the primary part
$\ical^1_N$ of the \szego kernel given in (\ref{lattice3}) as an
oscillatory integral over $\T\times P$:
\begin{equation} \ical_N^1= N^n \int_{{\bf T}^m}\int_P  e^{N\Psi(\phi,x;\tau,z)}A(\phi,x,N;\tau)\, dx\, d\phi\;, \label{mass-osc}
\end{equation}
where the phase is given by
\begin{equation}\label{tauphase}\Psi(\phi,x;\tau,z) =\langle \tau+i\phi, x\rangle
+p\log\left(\frac{ 1+\sum_{j=1}^m
e^{-\tau_j-i\phi_j}|z_j|^2}{1+\|z\|^2}\right)\;,\end {equation}
and the amplitude is given by \begin{equation}\label{amplitude}
A(\phi,x,N;\tau)= \xi_1(\phi)[A_0(x,\tau+i\phi)+
A_1(x,\tau+i\phi)N\inv+\cdots+ A_N(x,\tau+i\phi)N^{-n}]\;.\end
{equation}

 We therefore look for  complex critical points of $\Psi$. We
note that as before, by the triangle inequality,
\begin{equation}\label{maxatphi0}\Re
\Psi(\phi,x;\tau,z) <\Re\Psi(0,x;\tau,z)\;,\qquad \mbox{for\ \ }
\phi\ne 0\;.\end{equation} Thus for each fixed $\tau$ and $z$, the
maximal value of $\Re\Psi$ occurs where $\phi=0$. {From}
(\ref{tauphase}), we obtain
\begin{equation} d_\phi\Psi = ix-\frac {ip}
{1+\sum_je^{-\tau_j-i\phi_j}|z_j|^2}\left(e^{-\tau_1-i\phi_1}|z_1|^2,
\dots,
e^{-\tau_m-i\phi_m}|z_m|^2\right)\;.\label{dphicx}\end{equation}
Hence the critical-point equation $d_\phi\Psi(0,x;\tau,z)=0$ is
equivalent to:
\begin{equation} \frac{1}{p}\,
x =\mu(e^{-\tau/2}\cdot z)\;.\label{ii''}\end{equation}

{From} the second critical point equation
\begin{equation} \label{i''}d_x\Psi= \tau + i\phi=0\;,\end{equation} it follows
that there are no critical points on the interior of $P$, since by
(\ref{ii''}), we must choose $\tau\ne 0$ to obtain a critical
point when $z$ is in the classically forbidden region.

Suppose $\tau\in\R^m$ and $x^0\in P$.  Recall that $\tau$ is in
the normal cone to $P$ at $x^0$ if $\langle \tau,x^0 \rangle \ge
\langle\tau, y\rangle$ for all $y\in P$.  Since
$\Psi(\tau,0,x,z)=\langle\tau,x\rangle +c_{\tau,z}$, where
$c_{\tau,z}$ is independent of $x$, it follows that  $\Re\Psi$
takes its maximum along the contour $\T_\tau$  at the point
$\{\phi=0,\ x=x^0\}$ if (and only if) $\tau$ is in the normal cone
$C^{P}_{x^0}$ to $P$ at $x^0$. Recall that if $\tau$ is in the
normal cone at $x^0$, then $\tau$ must be orthogonal to the face
of $P$ containing $x^0$.

Therefore, to obtain a `critical point' that determines the
asymptotics of the integral (\ref{mass-osc}), we must find
$\tau\in\R^m, x\in \d P$ such that (\ref{ii''}) holds and $\tau\in
C_x$.  We show the existence of such $\tau, x$ in the next
section.

\subsubsection{Existence of critical points}\label{s-normalbundle}

Recall that the normal cone $C_F$ of a face $F$ of a convex
polytope $Q\subset\R^m$ is given by
\begin{equation}\label{ncone}C_F =\{u\in\R^m:\langle u,
x\rangle =\sup _{y\in Q}  \langle u,y\rangle,\ \forall x\in
F\}\;,\end{equation}  so that $C_F=C^Q_x$ for all points $x\in F$.
The purpose of this section is to prove the following general
existence result:

\begin{lem}\label{claim2-3} Let $Q$ be a convex polytope
contained in $\Si$ and suppose that $Q\not\subset \d\Si$. Then for
each $z\in(\C^*)^m$, there exists unique $\tau_z\in \R^m$ and $x_z\in
Q$ so that \begin{itemize}

\item $\ \mu(e^{-\tau_z/2}\cdot z)=x_z$,
\item
$\ \tau_z\in C_{x_z}^Q$.
\end{itemize}  \end{lem}

\begin{figure}[htb]
\centerline{\includegraphics*[bb= 0.7in 6.9in 4in 9.6in]{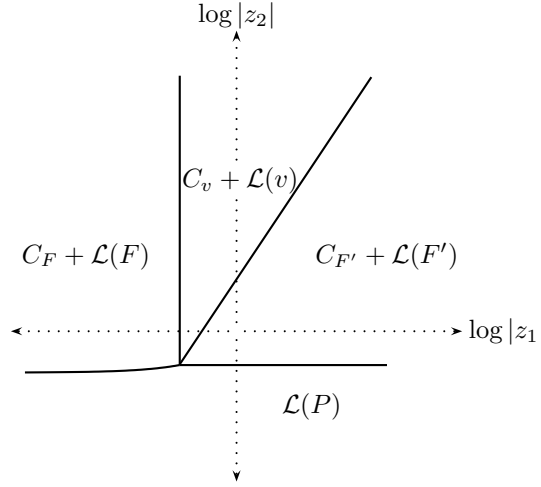}}
\caption{Transition points in log coordinates for the polytope of
Figure \ref{f-fan}} \label{L}\end{figure}

Here we do not assume that $Q$ is simple or even that $Q$ is
integral.  We note that even though $\tau_z$ is orthogonal to $\d
Q$ at $x_z$, the orbit $\{\mu(e^{r\tau_z}\cdot z):r\in\R\}$ is in
general not orthogonal to $\d Q$, even when $x_z$ lies in a facet.
(See also Figure \ref{L} and the associated remark.) Our proof of
Lemma \ref{claim2-3} uses the following elementary, but not so
well known, fact about the invertibility of Lipschitz maps.

\begin{lem}\label{lipeo} {\rm \cite{Fan}}
Let $f:U\to\R^m$ be a Lipschitz map, where $U$ is open in $\R^m$.
Then $f$ has a local orientation-preserving Lipschitz inverse with
Lipschitz constant $L$ at a point $x_0\in U$ if and only if there
exists $\ep>0$  such that
\begin{enumerate}
\item[i)] \  $\liminf_{v\to 0} |f(x+v) - f(x)|/|v|  \ge 1/L$ for all $x \in
B_\ep(x_0)$,
\item[ii)] \ $\det f'(x)>0$ for all $x \in
B_\ep(x_0)$ such that $f$ is differentiable at $x$,
\item[iii)] \ $f(x_0)\not\in f(\d B_\ep(x_0))$ and $\deg\big[f:\d B_\ep(x_0) \to
\R^m\sm\{f(x_0)\}\big] = 1$.
\end{enumerate}\end{lem}

\noindent Here, we let $B_\ep(x_0)=\{x\in\R^m:|x-x_0|<\ep\}$
denote the $\ep$-ball about $x_0\in\R^m$. By the degree in (iii),
we mean the degree of $f_*(\d B_\ep(x_0)) \in H_{m-1}(
\R^m\sm\{f(x_0)\},\Z)$.  The hypotheses (i) and (iii) in
\cite{Fan} differ slightly from those above, but Fan's proof of
sufficiency uses only (i)--(iii) above. (Necessity is obvious.)
Hypothesis (i) is given as a lemma in \cite{Fan}. Hypothesis (iii)
above is replaced in \cite{Fan} by the condition that the index of
$f$ at $x_0$ is 1, which is easily seen to be equivalent to (iii)
under the assumptions (i) and (ii).

Our proof of Lemma \ref{claim2-3} also makes use of the `normal
bundle' of a polytope, which we define below as the collection of
normal cones at the points of the polytope.

\begin{defin} The {\it normal bundle\/}
$\ncal(Q)$ of a convex polytope $Q\subset\R^m$ is the subset of
$T_{\R^m}=\R^m\times \R^m$ consisting of pairs $(x,v)$, where
$x\in Q$ and $v\in C^Q_x$. (Note that $\ncal(Q)$ is not a fiber
bundle over $Q$.)
\end{defin}

The normal bundle $\ncal(Q)$ is a piecewise smooth submanifold of
$T_{\R^m}$; it is homeomorphic to $\R^m$ via the `exponential map'
\begin{equation}\label{EQ}\ecal_Q:\ncal(Q)\to
\R^m\;, \quad \ecal_Q(x,v)= x+v,\qquad x\in F,\ c\in C_F\qquad (F\
\mbox{a face of} \ Q)\;.\end{equation} It is easily seen that $\ecal_Q$ is a
homeomorphism and is a $\ccal^\infty$ (in fact, linear)
diffeomorphism on each of the `pieces' $\bar F\times C_F\subset
\ncal(Q)$.

\medskip
\noindent{\it Proof of Lemma \ref{claim2-3}.\/} We note that the
$\R_+^m$ action on $(\C^*)^m$ descends via the moment map $\mu$ to
an $\R_+^m$ action on $\Si^\circ$ given by $r\odot
\mu(z)=\mu(r\cdot z)$. We let
$$\ncal^\circ=\textstyle\ncal(Q)
\cap(\Si^\circ\times\R^m)\;,$$ and we consider the map
$$\Phi:\ncal^\circ \to \Si^\circ\;,\qquad \Phi(x,\tau)=
e^{\tau/2}\odot x\;.$$ We observe that the conclusion of the lemma
is equivalent to the bijectivity of $\Phi$, since if
$z\in(\C^*)^m$, we have
$$\Phi(x,\tau)\eqd e^{\tau/2}\odot x =\mu(z) \iff
x=e^{-\tau/2}\odot\mu(z)\eqd \mu (e^{-\tau/2}\cdot z)\;,$$ for
$(x,\tau)\in\ncal^\circ$, i.e., for $x\in\Si^\circ$ and $\tau\in
C^Q_x$.

To show that $\Phi$ is a bijection, we let $\lcal:\Si^\circ
\rightarrow\R^m$ be the diffeomorphism given by
\begin{equation}\label{lcal}\lcal(x) =\left(\log \frac{x_1}{1-\sum x_j},\ldots,
\log \frac{x_n}{1-\sum x_j}\right)\;,\end{equation} so that
$$\lcal \circ \mu (w) = (\log |w_1|^2, \ldots, \log |w_m|^2)\;.$$ Thus,
writing $x = \mu(w)$ we have
\begin{eqnarray}\lcal \circ \Phi (x,\tau) &=&e^{\tau/2}\odot \,u_\Si(w) \ =\  \lcal
\circ\mu (e^{\tau/2}\cdot w)\nonumber\\ &=&(\tau_1 +\log
|w_1|^2,\ldots, \tau_m +\log |w_m|^2)\nonumber\\ &=& \tau + \lcal
\circ\mu (w)\ =\  \tau+\lcal (x)\;.\label{lcalphi}\end{eqnarray}
We first show that $\lcal \circ \Phi:\ncal^\circ\to\R^m$ is
proper: suppose on the contrary that the sequence $\{(x^n,
\tau^n)\}$ is unbounded in $\ncal^\circ$, but $\lcal \circ\Phi
(x^n, \tau^n) = \lcal(x^n)+\tau^n\rightarrow a\in \R^m$.  By
passing to a subsequence, we can assume that $x^n\rightarrow x^0
\in Q$.  Then $x^0 \in \partial \Si$, since otherwise
$\tau^n\rightarrow a-\lcal (x^0)$.  Write $\tau^n = r^n u^n$,
where $r^n>0$, $|u^n|=1$. We can assume without loss of generality
that $u^n\rightarrow u^0$. We first consider the case where $\sum
x^0_j<1$. If $x^0=0$, then $u^0_j\geq 0$ $(1\leq j\leq m)$, since
otherwise $\lcal \circ\Phi (x^n, \tau^n)$ would diverge.  But
$u^0\in C_{\{0\}}$; hence $Q\subset \{x:\langle u^0, x\rangle \leq
0\}\subset \R^m\backslash \Si^\circ$, a contradiction. Now suppose
that $x^0_1=\cdots=x^0_k=0$, $x^0_l>0$ for $k<l\leq m$. Then we
conclude as before that $u^0_j\geq 0$ for $1\leq j\leq k$.
Furthermore, in this case, $\tau_1^n \to +\infty$ and
$\tau_l^n\to\tau_l^\infty\in\R$ and thus $u^0_l = 0$ for $k<l\leq
m$, and we again obtain a contradiction. Finally, if $\sum
x^0_j=1$, we can suppose that $x^0_1\ne 0$ and make the coordinate
change $\tilde x_1=1-\sum x_j\ ,\tilde x_2=x_2,\ \dots,\ \tilde
x_m=x_m$ (which corresponds to permuting the homogeneous
coordinates in $\C\PP^m$) to reduce to the previous case.
Therefore, $\lcal\circ\Phi: \ncal^\circ \to \R^m$ is a proper map.

Let $\ecal=\ecal_Q:  \ncal(Q) \buildrel{\approx}\over\to \R^m$,
and let $U=\ecal (\ncal^\circ)\subset \R^m$. We consider the map
$f:U\rightarrow \R^m$ given by $f\circ
\ecal|_{\ncal^{\circ}}=\lcal \circ \Phi$, i.e., by the commutative
diagram:

$$\begin{array}{ccc} \ncal^\circ  &\buildrel{\Phi} \over \longrightarrow
& \Sigma^\circ \\ {\scriptstyle\approx}\!\downarrow
{\scriptstyle\ecal} & & {\scriptstyle\approx}\!\downarrow
{\scriptstyle\lcal}  \\ U &\buildrel{f} \over \longrightarrow &
\R^m \end{array}$$ Since $\lcal\circ\Phi$ is a proper map, $f$ is
also proper.  Hence to show that $\Phi$ is a bijection, it
suffices to show that $f$ is a local homeomorphism and is
therefore a (global) homeomorphism.

To describe the map $f$, for each $x\in \Si^\circ$, we let $q_x$
denote the point in $Q$ that is closest to $x$.  We note that
$q_x\in \Si^\circ$; if $x\not\in Q$, then $q_x\in
\partial Q$; if $x\in  Q$, then $q_x=x$.
Furthermore, $\ecal^{-1}(x) = (q_x, x-q_x)$ and hence
$$f(x) = x-q_x +\lcal (q_x)\;.$$

We shall show that $f$ satisfies the hypotheses of Lemma
\ref{lipeo}. Let $F$ be an arbitrary face of $ Q\cap \Si^\circ$.
To verify (i) and (ii), it suffices to show that $\det Df>0$ on
the (noncompact) polyhedron $$U_F:=\ecal(\bar F\times C_F)\cap
U=\ecal[(\bar F\cap \Si^\circ)\times C_F]\;.$$ To compute the
determinant, we let $T_F\subset \R^m$, $N_F=T_F^\perp\subset \R^m$
denote the tangent  and normal spaces, respectively,  of $F$.  Let
$x^0\in F$ be fixed.  For $y\in U_F$, we have
$$y=x + v \ \buildrel{f}\over \mapsto \ \lcal(x)+v\;, \qquad x-x^0\in T_F,\
v\in N_F\;.$$  Choose orthonormal bases $\{Y_1,\dots, Y_r\}$,
$\{Y_{r+1},\dots, Y_m\}$ of $T_F,\ N_F$, respectively. We let $Df$
denote the matrix of the derivative $(f|_{U_F})'$ with respect to
the basis $\{Y_1,\dots,Y_m\}$ of $\R^m$.  We have:
\begin{equation}\label{Df}Df=\left(\begin{array}{cc}
T^t\lcal'(x)T&0\\ *&I\end{array}\right)\;,
\end{equation} where $T$ is the $m\times r$ matrix $[Y_1 \cdots Y_r]$.
We have by (\ref{lcal}),
\begin{equation}\label{lcal'}\Big(\lcal'(x)\Big)_{jk} = \frac{1}{x_0} +
\delta^k_j \frac{1}{x_j}\;,\quad x_0=1-\sum_{j=1}^m
x_j>0\;,\end{equation} for $x\in\Si^\circ$. Hence $\lcal'(x)$ is a
positive definite symmetric matrix, it being the sum of a
semipositive matrix (all of whose entries are $\frac{1}{x_0}$) and
a positive definite diagonal matrix.  Therefore, $T^t\lcal'(x)T$
is positive definite, and hence $\det Df(x)=\det (T^t\lcal'(x)T)
>0$, completing the proof that hypotheses (i) and (ii) of Lemma
\ref{lipeo} are satisfied.

Note that (\ref{lcal'}) implies that the eigenvalues of
$\lcal'(x)$ are $>1$ for $x\in\Si^\circ$ and hence by (\ref{Df}),
$Df(x)$ is a diagonalizable matrix whose eigenvalues are   real
and $\ge 1$.

We verify (iii) by a homotopy argument:  Choose a point $q^0\in
Q^\circ$. We contract $ Q$ to $q^0$; i.e., for $0\le t\le 1$, we
let
$$Q_t= (1-t)\{q^0\} +  {t}Q\;,$$ so that $Q_0=\{q^0\},\ Q_1=
Q$. For $0\le t<1$, $Q_t\subset \Si^\circ$, and hence we have a
map
$$\Phi_t:\ncal(Q_t)\to \Si^\circ\;.$$ For $0\le t <1$, we define $f_t:\R^m\to
\R^m$ by the commutative diagram:
$$\begin{array}{ccc} \ncal(Q_t)  &\buildrel{\textstyle\Phi_t} \over
\longrightarrow & \Sigma^\circ \\[10pt] \ {\scriptstyle\approx}\!\downarrow
\ecal_{Q_t} & & {\scriptstyle\approx}\!\downarrow {\lcal}
\\ \R^m  &\buildrel{\textstyle f_t} \over \longrightarrow & \R^m \end{array}$$
As before we have
\begin{equation}\label{ft}f_t(x) = x-q^t_x +\lcal (q^t_x)\;,\end{equation}
where $q^t_x$ is the (unique) point of $Q_t$ closest to $x$. The
above argument shows that the maps $f_t$ also satisfy (i) and (ii)
of Lemma \ref{lipeo}.

We write
$$F:\big(\R^m\times [0,1)\big)\cup (U\times\{1\})\to \R^m\;,\qquad (x,t)\mapsto
f_t(x)\;,$$ where $f_1=f:U\to\R^m$.  One easily sees that
$$\left|q^t_x -q^{t'}_{x'}\right| \le |x-x'|+|t-t'|\;,$$  and
hence $F$ is continuous. Furthermore, $F$ is uniformly continuous
on $\R^m\times[0,t]$, for each $t<1$.

Let $H$ denote the set of $t\in [0,1)$ such that $f_t:\R^m\to\R^m$
is an orientation preserving homeomorphism.   We shall show that $H=[0,1)$. We first
observe that
$0\in H$ since
$$f_0(x)=x - q+\lcal(q)= x + \mbox{constant}\;.$$
Note that  (\ref{Df}) also holds for the maps $f_t$ and hence as
above, the eigenvalues of $Df_t$ are   real and $\ge 1$. Hence
$f_t$ satisfies conditions (i) and (ii) of Lemma \ref{lipeo} with
$L=1$, for all $t\in [0,1]$.

We now show that $H$ is open in $[0,1)$. Suppose $t_0\in H$.
Then by Lemma \ref{lipeo}, Lip$(f_{t_0}\inv)\le 1$ and hence
$$|x-x_0|=1 \Longrightarrow |f_{t_0}(x)-f_{t_0}(x_0)|\ge 1\;.$$  By the uniform
continuity of $F$ on $\R^m \times [0, \frac{t_0+1}{2}]$, we can
choose $\ep>0$ such that for all $x_0\in\R^m$, we have
\begin{equation}\label{x}|f_t(x)-f_t(x_0)|\ge \half\;, \qquad \mbox{for} \  x\in
\d B_1(x_0),\ |t-t_0|<\ep\;.\end{equation} To simplify notation,
we shall write
$$\deg(f,x_0,r):=\deg\big[f:\d B_r(x_0) \to
\R^m\sm\{f(x_0)\}\big] \;.$$ We conclude from (\ref{x}) that for
$|t-t_0|<\ep$,
$$\deg(f_{t},x_0,1) = \deg(f_{t_0},x_0,1)) = 1\;,$$ and hence
$f_t$ is a local homeomorphism at $x_0$,  by Lemma \ref{lipeo}. By
the very first part of the argument that $f:U\to \R^m$ is proper,
we easily see that $f_t$ is proper. Since $x_0\in\R^m$ is
arbitrary, it follows that $f_t:\R^m\to \R^m$ is a covering map
and therefore is a homeomorphism.

Next we show that $H$ is closed in $[0,1)$ and hence $H=[0,1)$.
Let $t_n\in H$ such that $t_n\to t_0\in [0,1)$. Since $f_{t_0}$
satisfies condition (i) of Lemma \ref{lipeo}, we can choose
$\ep>0$ so that $f_{t_0}(x_0)\not\in f_{t_0}(\d B_\ep(x_0))$. Then
for $n$ sufficiently large,
$$\deg(f_{t_0},x_0,\ep) = \deg(f_{t_n},x_0,\ep) = 1\;.$$  It follows as
above that $f_{t_0}:\R^m\to \R^m$ is  a homeomorphism.

We have shown that $f_t$ is a homeomorphism for $0\le t <1$. To
complete the proof of the lemma, we must show that $f=f_1$ is a
local homeomorphism.  So we let $x_0\in U$ be arbitrary, and we
choose $\ep>0$ such that $\overline {B_\ep(x_0)}\subset U$ and
$f(x_0)\not\in f(\d B_\ep(x_0))$. Then for $t<1$ sufficiently
close to 1, we have as before $$\deg(f,x_0,\ep) =
\deg(f_t,x_0,\ep) = 1\;.$$ Thus by Lemma \ref{lipeo}, $f$ is a
local homeomorphism.  This completes the proof of Lemma
\ref{claim2-3}.\qed

\bigskip

We have thus shown (as a consequence of  (\ref{maxatphi0}) and Lemma \ref{claim2-3} with
$Q=\frac 1p P$) that for each $z\in(\C^*)^m$,
there exist (unique) $\tau_z\in\R^m$ and $q(z)=px_z\in  P$ such
that the phase function $\Psi(\cdot, \cdot;\tau_z, z):\T\times P
\to \C$ given by (\ref{PHASEcx}) takes its maximum real part at
the point $(0,q(z))$, and furthermore $(0,q(z))$ is a critical
point of  $\Psi|_{\T\times F}$, where $F$ is the face of $P$
containing $q(z)$.  The latter statement follows from
(\ref{ii''}), (\ref{i''}) and the fact that vectors in the normal
cone $C_F$ are orthogonal to $F$; i.e., $C_F\subset F^\perp$.

By the {\it boundary\/} of the normal cone $C_F$, we mean the boundary $\d C_F$ of $C_F$
in
$F^\perp$; the complement $C_F\sm \d C_F$ is called the {\it interior\/} of $C_F$.  We
note the following consequence of the proof of Lemma~\ref{claim2-3}:

\begin{lem} A point $z\in (\C^*)^m$ is a  transition
point if and only if  $\tau_z \in\d C^P_{q(z)}$. \label{tpoint}\end{lem}

\begin{proof} It follows from (\ref{lcalphi}) that $2\,\Log(\rcal_F)=\lcal(F)+C_F$,
where $\Log(z)=(\log|z_1|,\dots,\log |z_m|)$ and that
\begin{equation}\label{Log}\d
[2\,\Log(\rcal_F)]=2\,\Log(\d\rcal_F)=\lcal(F)+\d
C_F\;.\end{equation} Thus
$$z\in\d\rcal_F \iff 2\,\Log(z)\in\lcal(F)+\d C_F \iff 2\,\Log(z)= \lcal(x_z) +\tau_z, \
\tau_z\in \d C_F\;.$$ \end{proof}

\begin{rem} The identity (\ref{Log}) tells us that we can decompose
$\R^m$ as the disjoint union of the sets $\lcal(F)+C_F$.  If $z$
is a transition point, then $2\,\Log(z)$ must lie  in the common
boundary of at least two of the sets $\lcal(F)+C_F$.
Figure~\ref{L} below shows the transition points in log
coordinates (as solid lines) for the case where $P$ is the
polytope of Figure~\ref{f-fan}. Note that although $C_{F'}$ is
orthogonal to $F'$ in Figure~\ref{f-fan}, $C_{F'}$ is not
orthogonal to $\lcal(F')$.  Similarly, in Lemma~\ref{claim2-3}, if
$q(z)\in F'$, then the orbit $\{\mu(e^{r\tau_z}\cdot z):r\in\R\}$
is not orthogonal to $F'$. \end{rem}

\medskip \subsubsection{Continuation of the proof of Theorem~\ref{SZEGO}.}
\label{s-asymptotic}  Since the critical point of $\Psi$ on the
contour $\T_\tau$ always occurs on the boundary of $P$ when $z$ is in the
forbidden region, we shall
 use the method of  stationary phase for complex
oscillatory integrals on quadrants.  Since this general method is
not so well known, we state here the basic result we need:

\begin{lem} \label{halfplane} Let $\Psi(\xi,y),A(\xi,y)\in
\ccal^\infty(\R_{\ge0}^k\times \R^s)$ such that $d_y\Psi(0,0)=0$,
$A$ has compact support and
\begin{enumerate}
\item[1)] $\frac{\d\Psi}{\d \xi_j}\ne 0$ on $\supp(A)$, for $1\le j\le k$,
\item[2)] $d_y\Psi(0,y) \ne 0$ for $(0,y)\in\supp(A)\sm\{0\}$,
\item[3)] $\det
\hcal_y\Psi(0,0)\ne 0$ (where $\hcal_y$ denotes the Hessian with
respect to $y$),
\item[4)] $\Re\Psi\le \Re \Psi(0,0)$ on $\supp(A)$.
\end{enumerate}
Then $$\int_{\R^s}\int_{\R_{\ge 0}^k}e^{N\Psi(\xi,y)} A(\xi,y) \,
d\xi\,dy =
N^{-k-s/2}e^{N\Psi(0,0)}[c_0+c_1N\inv+c_2N^{-2}+\cdots+c_lN^{-l}
+O(N^{-l-1})]$$ for $l=1,2,3,\dots,$ where
$$c_0=\left.\frac {(2\pi)^{s/2}\,
A}{\sqrt{\det(-\hcal_y\Psi)}\;\prod_{j=1}^k
\d\psi/\d\xi_j}\right|_{\xi=y=0}\,.$$ (If $s=0$, then hypotheses
(2) and (3) are vacuous, and the determinant in the formula for
$c_0$ is 1.)

\end{lem}

\begin{proof} (See also \cite{AGV}.) Let us first consider the case $k=1$.
Integrating by parts,
\begin{equation}\label{parts1}\int_{\R^s}\int_0^{+\infty}e^{N\Psi} A \,
d\xi_1\,dy =\frac 1N \int_{\R^s}\left.e^{N\Psi} \frac{A
}{\d\Psi/\d{\xi_1}}\right|_{\xi_1=0}\,dy -\frac 1N
\int_{\R^s}\int_0^{+\infty}e^{N\Psi} \frac
\d{\d\xi_1}\left[\frac{A }{\d\Psi/\d{\xi_1}}\right] \,
d\xi_1\,dy\;.\end{equation} Applying the stationary phase expansion
\cite[Th.~7.7.5]{H} to the first term of (\ref{parts1}) and
iterating, we obtain the desired expansion.

The case $k>1$ follows by induction on $k$, integrating by parts
as above.\end{proof}

\medskip We  consider a point $z$ in the classically forbidden
region. Let $\tau_z,x_z$ be as in Lemma \ref{claim2-3}, and let
$q(z)=px_z$.  Hence $\tau=\tau_z,x=q(z)$ satisfies the critical
point equation (\ref{ii''}), and as we observed above, $q(z)\in
\partial P$.  The decay function $b_P(z)$ from (\ref{b}) is given by
\begin{equation} b_P(z)=-\Psi(0,q(z))\;.\end{equation}

Recalling (\ref{mass-osc}), we further decompose the integral
$\ical^1_N$ into two parts:
\begin{equation} \ical_N^1=\ical_N'+\ical_N'',\qquad \ical_N'= N^n \int_{{\bf T}^m}\int_P
e^{N\Psi(\phi,x;\tau_z,z)}\rho(x)A(\phi,x,N;\tau_z)\, dx\,
d\phi\;, \label{mass-osc2}
\end{equation} where
$A(\phi,x,N;\tau_z)$ is given by (\ref{amplitude}) and $\rho$ is
supported in a small neighborhood of $q(z)$ and is $\equiv 1$ near
$q(z)$. Note that it follows from (\ref{dphicx})--(\ref{ii''})
that $d_\phi \Psi(0,x;\tau_z,z)\ne 0$
 for $x\ne q(z)$. Since $d_\phi\Psi$ does not vanish on the support of
$1-\rho(x)$ and  $\sup _{\T\times P}\Psi=\Psi(0,q(z))=-b_P(z)$, we
conclude by performing the $\phi$ integration first, that
\begin{equation}\label{smallI}\ical_N''=e^{-Nb_P(z)} \cdot O(N^{-\infty})\end{equation}
(see \cite{H}).  Hence, we need consider only $\ical'_N$.

Since we are not assuming that $P$ has interior, we let $n=\dim P$
and we let $\langle P\rangle$ denote the $n$-plane containing $P$.
The polytope $P$ can be given by
$$P=\{x\in\langle P\rangle: \ell_j(x)\ge 0 \ \mbox{for }\ 1\le
j\le d\}\;, \qquad \ell_j(x)=\langle u_j,x\rangle +\la_j\,,$$
where the $u_j$ are primitive vectors tangent to $P$. Suppose that
$q(z)$ is in a face $F\subset P$ of dimension $r$. After permuting
the indices of the normal vectors $\{u_j\}$, we may assume that
$$\bar F = \{x\in P:  \ell_j(x)= 0, \quad \mbox{for }\
1\le j\le n-r\}\;,$$  We consider the translated cone
$$\Ga:=\{x\in\langle P\rangle:\ell_j(x)\ge 0\}\supset P$$ with vertex
$q(z)$, which coincides with $P$ in a sufficiently small
neighborhood of $q(z)$. We then choose (non-homogeneous) linear
functions $t_1,\dots,t_r$ on $\R^m$ such that the map
$$L=(\ell_1,\dots,
\ell_{n-r}, t_1,\dots,t_r):\langle P\rangle\ \to \ \R^n$$ is
bijective and $L(q(z))=0$. We note that $L:\Ga\ {\buildrel
{\approx}\over \to} \ \R_{\ge 0}^{n-r} \times \R^{r}$ and
$L:F\hookrightarrow \{0\}\times \R^{r}$.

We shall evaluate $\ical'_N$ by applying Lemma \ref{halfplane}
with $\xi_j=\ell_j(x)$ for $1\le j\le n-r$ and
$y=(t_1,\dots,t_r,\phi_1,\dots \phi_m)$ and with the phase
$$\wt\Psi(\xi,y):=\Psi\left(\phi,L\inv(\xi;t_1,\dots,t_r);\tau,z\right)$$ and
amplitude $\rho(x)A(\phi,x,N;\tau)$. Since $d_\phi
\Psi(0,q(z);\tau_z,z)= 0$ and $\tau_z$ is perpendicular to $F$, it
follows (using (\ref{tauphase})) that $d_y\wt\Psi(0,0)=0$.   Since
$\tau_z$ is in the normal cone at $q(z)$, the real part of the
phase $\Psi$ takes its maximum at $q(z)$, and hence hypothesis (4)
of the lemma holds. Also, hypothesis (2) is satisfied since
$d_\phi \Psi(0,x;\tau_z,z)\ne 0$
 for $x\ne q(z)$, and hence
$d_y \wt\Psi(0,y)\ne 0$ for $y\ne 0$.

We now verify (3): Recalling (\ref{Hessian}), we note that the
Hessian matrix (with respect to $\phi,x$) of $\Psi$ is given by
\begin{equation}\label{d2cx} \hcal\Psi(\phi,x;\tau,z)= \hcal
\Psi_{\!\acal}(\phi,x;e^{-\tau/2}\cdot z)\;,\end{equation} and
hence
\begin{equation}\label{Hessian2}\hcal\Psi|_{(0,x)} =
\left(\begin{array}{ll} {\mathbf C} & i\I
\\ &
\\ i\I & 0
\end{array}
\right)\end{equation} where
\begin{equation}\label{C}C_{jk} = \left.
\frac{\d^2\Psi}{\d\phi_j\d\phi_k}\right|_{(0,x)} =p(
-I_j\delta_j^k +I_jI_k)\;,\qquad
(I_1,\dots,I_m)=\mu(e^{-\tau/2}\cdot z)\;.\end{equation}

\begin{lem}\label{claim} The matrix ${\mathbf C}$ given by (\ref{C}) is strictly negative
definite.\end{lem}

\begin{proof}  Suppose that $\la=(\la_1,\dots,\la_m)\in\R^m\sm\{0\}$ is
arbitrary. Let $v=(v_1,\dots,v_m),\ w=(w_1,\dots,w_m) $ be given
by $v_j=\sqrt{I_j},\ w_j= \sqrt{I_j}\,\la_j$.  Since
$e^{-\tau/2}\cdot z\in(\C^*)^m$, $I_j\ne 0$ for $1\le j\le m$, and
hence $w\ne 0$.  Furthermore $\|v\|^2=\sum I_j <1$. Therefore
$$\frac 1p \sum_{j,k=1}^m C_{jk}\la_j\la_k= -\|w\|^2+(v\cdot
w)^2 < -\|v\|^2\|w\|^2+(v\cdot w)^2\le 0\;.$$ \end{proof}

\bigskip
Changing variables, we see that $\hcal_y\wt\Psi(0,0)$ is of the
form
$$\left(\begin{array}{cc} 0
& i{\mathbf B}
\\ &
\\ i{\mathbf B}^t & {\mathbf C}
\end{array}
\right)$$ where ${\mathbf B}$ is an $r\times m$ matrix of maximal
rank $r$ and ${\mathbf C}$ is nonsingular; hence
$\hcal_y\wt\Psi(0,0)$ is nonsingular, verifying hypothesis (3).

To verify (1), we must assume that {\bf  $z$ is not a transition
point}. By Lemma \ref{tpoint}, $\tau_z$ is not in the boundary of
the normal cone to $F$, and hence $$\langle \tau_z, x-q(z)\rangle
<0\quad \mbox{for all }\ x\in P\sm\bar F\;.$$
 Let $q^j=L\inv(E^j)\in\Ga$, for $1\le j\le n-r$ (where
$E^j=(\delta^j_1, \delta^j_2,\dots,\delta^j_n)\in\R^n$); i.e.,
$q^j$ lies on the $j$-th edge of $\Ga$. Let $\eta_j=q^j-q(z)$.
Then
$$\d\wt\Psi/\d\xi_j= \nabla_{\eta_j}\Psi \equiv \langle \tau_z, v_j\rangle<0\qquad (1\le
j\le n-r)\;,$$  and hence hypothesis (1) holds.

Therefore, by (\ref{mass-osc}), Lemma \ref{halfplane} (with $k=n-r,\ s=m+r$) and
(\ref{smallI}),
\begin{equation}\label{asympt}\ical^1_N(z)= N^{\frac {r-m}{2}}
e^{N\Psi(0,q(z);\tau_z,z)} \big[c_0+c_1N\inv+\cdots+c_lN^{-l}
+O(N^{-l-1})\big]\;.\end{equation} We let
$$b_P(z)=-\Psi(0,q(z);\tau_z,z)=-\langle q(z), \tau_z\rangle + p
\log\left(\frac{1+\|z\|^2}{1+\|e^{-\tau_z/2}\cdot z\|^2}\right)$$
as in (\ref{b}).

We now show that $\ical_N^2$ decays faster than $\ical_N^1$:
Recall that
$$\ical_N^2=
\frac{1}{(2\pi )^m} \int_{\T } \xi_2(\phi) \left(\frac{
1+\sum_{j=1}^m
e^{-\tau_j-i\phi_j}|z_j|^2}{1+\|z\|^2}\right)^{pN}{\chi_{NP}(e^{\tau_z+i\phi})}\,
d\phi\;,\quad \xi_2(\phi)=0 \ \mbox{for }\
\mbox{dist}_{\T}(\phi,0)\ge \ep\;,$$ where we write
$\tau_z=(\tau_1,\dots,\tau_m)$. We choose $\la=\la(\ep,z)>0$ such
that
$$\left|\frac{ 1+\sum_{j=1}^m
e^{-\tau_j-i\phi_j}|z_j|^2}{1+\|z\|^2}\right| \le e^{-\la}\ \frac{
1+\sum_{j=1}^m e^{-\tau_j}|z_j|^2}{1+\|z\|^2}\quad \mbox{for } \
\mbox{dist}_{\T}(\phi,0)\ge \ep\;.$$ Therefore,
\begin{equation}\label{s1} |\ical_N^2| \le e^{-Np\la}\left(\frac{ 1+\sum_{j=1}^m
e^{-\tau_j}|z_j|^2}{1+\|z\|^2}\right)^{pN}{\chi_{NP}(e^{\tau_z})}\;.\end{equation}
Furthermore,
\begin{equation}\chi_{NP}(e^{\tau_z})= \sum_{\al\in NP}e^{N\langle\tau_z,\al/N\rangle}\le
\#(NP)\,e^{N\langle\tau_z,q(z)\rangle} \le
(Np)^ne^{N\langle\tau_z,q(z)\rangle}\;,\label{s2}\end{equation}
since $\al/N\in P$ and $\tau_z$ is in the normal cone to $P$ at
$q(z)$. Combining (\ref{s1}) and (\ref{s2}), we obtain
\begin{equation} |\ical_N^2| \le (Np)^n e^{-Np\la} e^{-Nb_P(z)}
\;.\label{s3}\end{equation}

The asymptotic expansion of part (ii) of Theorem~\ref{SZEGO} (for
$P$ simple) with $b_P(z)$ given by (\ref{b}) now follows from
(\ref{lattice3}), (\ref{asympt}) and (\ref{s3}).

\subsubsection{The decay function} \label{s-decay} In this section
we verify that the function $b_P$ given by (\ref{b}) satisfies
parts (c) and (e) of Theorem~\ref{SZEGO}(ii) for any convex
polytope $P$.

We first note that when $z\in{\acal_P}$, we have $\tau_z=0$ and
hence if we extend (\ref{b}) to all of $(\C^*)^m$, we then have
$b_P=0$ on $\acal_P$. Furthermore, since $z\mapsto \tau_z$ is
easily seen to be continuous on $(\C^*)^m$ and  $\ccal^\infty$ on
each $\overline{\rcal_F}$, the same holds for $b_P(z)$. We now
show that $b_P$ is $\ccal^1$  by computing its derivative.
Recalling that $b_P(z)=-\Psi(0,q(z);\tau_z,z)$, we have
$$ - db_P = \big[d_x  \Psi(0,x;\tau,z)\cdot D q(z) +d_\tau \Psi(0,x;\tau,z) \cdot
D\tau_z +  d_z
 \Psi(0,x;\tau,z)\big]_{x=q(z),\tau=\tau_z}\;. $$
Here, $D\tau_z$ and $Dq(z)$ are only piecewise continuous, being
discontinuous at transition points. However, from (\ref{tauphase})
and the fact that $(0,q(z))$ is a critical point of $\Psi$, we
have
$$d_\tau \Psi(0,x;\tau,z)|_{\{x=q(z),\tau=\tau_z\}} = -i\,d_\phi
\Psi(\phi,x;\tau,z)|_{\{\phi=0,x=q(z),\tau=\tau_z\}}=0\;.$$
Furthermore, since $d_x\Psi|_{\phi=0,\tau=\tau_z}=\tau_z$ and
$\tau_z$ is perpendicular to the face containing $q(z)$, we also
have $$d_x  \Psi(0,x;\tau,z)|_{\{x=q(z),\tau=\tau_z\}}\cdot D q(z)
=\tau_z\cdot Dq(z)=0\;.$$ Therefore
\begin{equation}\label{Db} db_P= - d_z
\Psi(0,x;\tau,z)|_{\{x=q(z),\tau=\tau_z\}}\in \ccal^0((\C^*)^m)\;,
\end{equation} completing the proof of (e).

Property (c) is a special case of the following  lemma on the
dependence of the decay function $b_P(z)$ on the polytope.

\begin{lem} Let $P'\subset P \subset p\Si$ be convex polytopes in
$\R^m$.  Suppose that $z\in (\C^*)^m$ and let $q(z)\in P$
be given as in (\ref{b}). Then $$b_{P'}(z) \ge b_P(z)\;,$$ with
equality if and only if $q(z) \in P'$. \label{monotone}\end{lem}

\begin{proof} Fix $z\in(\C^*)^m$ and consider the function
\begin{equation}\label{bx}b (x;z):= -\Psi(0,x;\tau(x),z)=-\langle x,
\tau(x)\rangle + p
\log\left(\frac{1+\|z\|^2}{1+\|e^{-\tau(x)/2}\cdot
z\|^2}\right)\;, \quad x\in \Si^\circ\;,\end{equation} where
$\tau(x)\in \R^m$ is given by
\begin{equation}p\mu(e^{-\tau(x)/2}\cdot z)=x\;.\label{taux}\end{equation}
Note that $\tau(q(z))=\tau_z$ and $b (q(z);z)=b_P(z)$. Hence it
suffices to show that $b (x;z)
>b (q(z);z)$ for all $x\in P\sm \{q(z)\}$.

Equation (\ref{taux}) yields $e^{-\tau_j(x)}|z_j|^2 =
\frac{x_j}{1-\sum_k x_k}$, and therefore
\begin{equation}\label{tauj} \tau_j(x)=\log|z_j|^2 +\log(1-\sum_k
x_k)-\log x_j \;.\end{equation} Since
$$d_\tau\Psi|_{\{\phi=0, \tau=\tau(x)\}} = \frac 1i d_\phi\Psi|_{\{\phi=0,
\tau=\tau(x)\}}= x-p\mu(e^{-\tau(x)/2}\cdot z)=0\;,$$ we conclude
from (\ref{bx}) that \begin{equation} \label{dxb} d_xb=
-d_x\Psi|_{\{\phi=0, \tau=\tau(x)\}}=-\tau(x)\;,\end{equation} and
hence
\begin{equation}\label {Hb} \frac{\d^2 b}{\d x_j \d x_k} =
-\frac{\d\tau_k}{\d x_j}= \frac 1 {1-\sum_{l=1}^m x_l} +\delta^k_j
\frac 1 {x_j}\;.\end{equation} Therefore, the Hessian $(\frac{\d^2
b}{\d x_j \d x_k})$ is positive definite. (It coincides with the
derivative $\lcal'(x)$ from \S \ref{s-normalbundle}.)

Now suppose that $x^0\in P$, $x^0\ne q(z)$.  Let
$$f(t)= b (q(z)+tv;z)\;,\qquad v=x^0-q(z).$$  Then by (\ref{dxb}),
$f'(0) = -\langle v, \tau_z\rangle \ge 0$, and by the positivity
of the Hessian, $f''(t) >0$.  Therefore $b (x^0;z)-b (q(z);z) =
f(1)-f(0)>0$.
\end{proof}

Now let $z\in (\C^*)^m\sm \overline{\acal_P}$.  We apply Lemma
\ref{monotone} with $P',P$ replaced with $P,p\Si$, respectively.
Since $q_{p\Si}(z) =p\mu(z)\not\in P$ by
assumption, we conclude by Lemma \ref{monotone} that
$$b_P(z)>b_{p\Si}(z) = 0\;.$$  Thus we have shown that the
function $b_P$ satisfies (c) and (e) for all convex polytopes $P$.

\subsubsection{Precise asymptotics on the classically allowed region.}
\label{s-precise} In this section, we prove part (i) of
Theorem~\ref{SZEGO}. We let $P$ be a convex integral polytope
(not necessarily simple). By  (\ref{szego-proj}),
$$\Pi_{Np\Si}(z,z)= \sum_{\al\in Np\Si}\
\frac{1}{\|\chi_\al\|^2}
|\wh\chi^{Np}_\al(z)|^2=\frac{(Np+m)!}{(Np)!} =\prod_{j=1}^m
(Np+j)\;,$$ and we then have by (\ref{Sz2}),
$$\Pi_{|NP}(z,z) = \prod_{j=1}^m (Np+j) - R_N(z)\;,$$ where
\begin{equation} \label{Piallow}R_N(z)= \sum_{\al\in Np\Si\sm NP}\
\frac{1}{\|\chi_\al\|^2} |\wh\chi^{Np}_\al(z)|^2\;.\end{equation}

To simplify our computations, we introduce the renormalized monomials
\begin{equation}\label{m}\wh
m_\al^{Np}(z):= \left[
\frac{(Np)!}{(Np+m)!}\right]^\half \frac{\wh
\chi_\al(z)}{\|\chi_\al\|}={Np\choose \alpha}^\half \wh\chi^{Np}_\al(z)={Np\choose
\alpha}^\half
\frac{z^\al}{(1+\|z\|^2)^{Np/2}}\;, \quad |\al|\le Np\,.\end{equation}
Let $\al\in  Np\Si\sm NP$ be fixed. We easily check that
\begin{equation}|\wh m^{Np}_\al(z)|^2=\frac{(Np)!}{(Np+m)!}\Pi_{|\{\al\}}
(z,z)= \frac 1 {(2\pi)^m} \int_{\T}
e^{\Psi_N(\phi,\al;\tau,z)}\,d\phi\;,\label{mon}\end{equation}
where the phase $\Psi_N$ is given by\begin{equation}\label{tauphaseN}\Psi_N(\phi,\al;\tau,z)
=\langle \tau+i\phi, \al\rangle +Np\log\left(\frac{ 1+\sum_{j=1}^m
e^{-\tau_j-i\phi_j}|z_j|^2}{1+\|z\|^2}\right)\;.\end {equation}
(Equations (\ref{mon})--(\ref{tauphaseN}) also follow from (\ref{Sz2}) and
(\ref{lattice2}) with $P=\{\al\}$ and $p$ replaced by $Np$; $\Psi_N$ is
the phase (\ref{tauphase}) with $p$ replaced by $Np$.)
Since $\Re \Psi_N(\phi,x;\tau,z) \le \Psi_N(0,x;\tau,z)$, we have
\begin{equation}|\wh m^{Np}_\al(z)|^2 \le
e^{\Psi_N(0,\al;\tau,z)}\;.\label{mbound}\end{equation}   Choosing $\tau$ such that
$Np\mu(e^{-\tau/2} \cdot z)=  \al$ and recalling (\ref{bx}), we have
\begin{equation}-\Psi_N(0,\al;\tau,z)  =N b\left(\frac \al N
;z\right) = Nb_{\{ \frac \al N\}}(z)\;.\label{Nb}\end{equation} By
Lemma \ref{monotone}, $b_{\{x\}}(z) >b_P(z)=0$ for $x\in p\Si\sm
P^\circ$ and $z\in \acal_P$ (i.e.\ $p\mu(z)\in P^\circ$).

Now let $K$ be a compact subset of $\acal_P$, and let
$$\la_K=\half\inf \{b_{\{x\}}(z):x\in
p\Si\sm P^\circ,\ z\in K\}>0\;.$$  Therefore,
\begin{equation} |\wh m^{Np}_\al(z)|^2 \le e^{ Nb_{\{
\frac \al N\}}(z)} \le e^{-2\la_K N}\;,\quad \mbox{for }\ z\in K,\
\al\in  Np\Si\sm NP,\ N\ge 1. \label{mbd}\end{equation} Since
$\#(Np\Si\sm NP)\le (Np)^m$, we conclude from (\ref{Piallow}) and
(\ref{mbd}) that $\|R_N\|_{\ccal^0(K)}=O(e^{-\la_K N})$.

To verify the $\ccal^1$ estimate, we recall from (\ref{mchi}) that
$$|\wh m^{Np}_\al(z)|^2 = {Np\choose \al}h(z)^N|z^\al|^2\;,
\qquad h(z)=(1+\|z\|^2)^{-p}\;.$$  Differentiating with respect to
the variables $\rho_j=\log |z_j|$, we obtain
$$d_\rho (|\wh m^{Np}_\al(z)|^2) = {Np\choose \al}h(z)^N|z^\al|
^2\left[ 2\al +N d_\rho\log h\right] = |\wh m^{Np}_\al(z)|^2
\left[ 2\al +N d_\rho\log h\right] \;,$$ and hence the $\ccal^1$
estimate follows  from the $\ccal^0$ estimate. Differentiating
repeatedly, we obtain all the $\ccal^k$ estimates,  completing the
proof of part (i) of Theorem~\ref{SZEGO}.

\subsubsection{Non-simple polytopes} \label{s-general} In this
section, we prove the general case of part (ii) of Theorem~\ref{SZEGO} by a
reduction to the simple case.

Let $P$ be a non-simple convex integral polytope of dimension $n$.
Consider a face $F$ of $P$
 of dimension $r$ and let $z\in \rcal_F^\circ$.
 We say that $P$ is simple at $F$ if
$\#\jcal(F) = n-r$. (Recall the notation from \S \ref{s-fans}.
Note that we always have $\#\jcal(F) \ge n-r$.  A polytope is
simple $\iff$ it is simple at each of its vertices $\iff$ it is
simple at each of its faces.)

We first consider the case where $P$ is simple at $F$. (This is
always the case when $\codim F\le 2$.)   We  construct a simple
integral polytope $Q\supset P$ coinciding with $P$ near $q(z)$ as
follows: Recall that $q(z)\in F$ (by the definition of $\rcal_F$),
and consider the `barrier cone' $B_F(P)$ generated by elements of
the form $x-q(z)$, where $x\in P$. Choose $L\in \mbox{GL}(m,\Q)$
such that $$L:B_F(P)\approx \R_{\ge 0}^{n-r}\times \R^r\times
\{0\}^{m-n}\;.$$ (E.g., let $L$ map the interior normals
$\{u_j:j\in\jcal_F\}$ at $F$ to the first $n-r$ standard basis
vectors in $\R^m$ and map a basis for the tangent space of $F$ to
the next $r$ standard basis vectors.) We then let
$$Q=q(z)+L\inv\big([0,M]^{n-r}\times [-M,M]^r\times \{0\}^{m-n}\big)\;,$$
where $M\in\Z^+$ is chosen so that the vertices of $Q$ lie in
$\Z^m$. By the construction, $B_F(Q)=B_F(P)$ and hence $Q$
coincides with $P$ near $q(z)$. We then replace $M$ with a
sufficiently high multiple if necessary, so that $Q\supset P$.

Since the integral polytope $Q$ may not be contained in $p\Sigma$,
we consider the convex rational polytope $Q':=Q\cap p\Sigma$. Recalling (\ref{m}), we
have
$$ \Pi_{|NP}(z,z)= \Pi_{|NQ'}(z,z)
-\frac{(Np+m)!}{(Np)!}\sum_{\al\in  NQ'\sm NP}|\wh
m_\al^{Np}|^2\;,$$ where by the \szego kernel $\Pi_{|NQ'}$, we
mean $\Pi_{|K_N}$, where $K_N$ is the convex hull of $\Z^N \cap
NQ'$. By (\ref{mbound})--(\ref{Nb}), for $\al\in (NQ'\sm
NP)\cap\Z^m$, we have
$$|\wh m_\al^{Np}|^2 \le \exp\left(-N b_{\{\al/N\}}(z)\right)\le
e^{-NC}\;,$$ where
$$C= \inf \{b_{\{x\}}(z): x\in Q'\sm P\}\;.$$  Since $Q'$ coincides with $P$ at $q(z)$,
$\tau_z$ is in the normal cone to $Q'$ at $q(z)$ and hence
$q(z)=q_{Q'}(z)$. Therefore by Lemma \ref{monotone},
$b_{\{x\}}(z)>b_{Q'}(z)=b_P(z)$ for all $x\in Q'$, $x\neq q(z)$,
and hence $C>b_P(z)$. Since $\#(NQ')\le \frac 1{m!}(Np)^m$, it
follows that
$$\|\Pi_{|NP}(z,z)- \Pi_{|NQ'}(z,z)\|_{\ccal^0(K)} = O(e^{-NC'}),
\qquad b_P(z) < C' < C, \ K \subset \!\subset \rcal^\circ_F.$$ The
similar $\ccal^l$ estimates follow as in the proof of
Theorem~\ref{SZEGO}(i) in \S \ref{s-precise}.

Thus, to obtain the asymptotic expansion of
Theorem~\ref{SZEGO}(ii)  for $P$, it suffices to obtain the
same expansion for $Q'$. To do this, we make the simple
observation that
\begin{equation}\label{latticeQ}\Pi_{|NQ'} (z,z) =
\frac{1}{(2\pi )^m} \int_{\T } \Pi_{Np} (z,e^{i\phi}\cdot z)
{\chi_{NQ'}(e^{i\phi})}\,d\phi =
 \frac{1}{(2\pi )^m}
\int_{\T } \Pi_{Np} (z,e^{i\phi}\cdot z)
{\chi_{NQ}(e^{i\phi})}\;,\end{equation} since lattice points in
$NQ\sm Np\Si=NQ\sm NQ'$ do not contribute to the integral. By
repeating word for word the argument at the beginning of \S
\ref{s-critical} and in \S \ref{s-asymptotic}, with
(\ref{latticeQ}) in place of (\ref{lattice}), we obtain the
expansion $$\Pi_{|NQ'}(z,z)= N^{\frac{m+r}{2}} e^{-N
b_P(z)}\big[c_0^F(z) + c_1^F(z)N\inv+\cdots + c_k^F(z)N^{-k}
+R_k^F(z)\big],$$ with $c_0^F>0$, and hence part (ii) of
Theorem~\ref{SZEGO} holds for $P$.

Now suppose that $P$ is not simple at $F$.  We subdivide the
barrier cone $B_F(Q)$ into simple rational cones
$\{B^k:k=1,\dots,s\}$ so that $B^k\supset B_F(F)$ for all $k$ and
intersections of any number of the $B^k$ are faces of each of
them. (Note that $B_F(F)=T_F$.) For a non-empty subset $\ecal$ of
$\{1,\dots,s\}$, let $P_\ecal=P\cap\bigcap_{k\in\ecal}(q(z)+B^k)$,
which is simple at $F$. Since $\tau_z$ is in the interior of the
normal cone of $P$ at $q(z)$, $\tau_z$ is also in the interior of
the normal cone at $q(z)$ of each of  the $P_\ecal$. (In
particular, $z$ is not a transition point for any of the
$P_\ecal$.) Thus by the case proven above (which does not use the
fact that the vertices of $P$ are lattice points, only that the
vertices of $Q'$ lie in $\Z^m$), we have the expansion
\begin{equation}\label{ecal}\Pi_{|NP_\ecal}(z,z)=
N^{\frac{m+r}{2}} e^{-N b_P(z)}\big[c_0^\ecal(z)+
c_1^\ecal(z)N\inv+\cdots + c_k^\ecal(z)N^{-k} +R_k^F(z)\big]\;,
\qquad z\in\rcal_F^\circ\;,\end{equation}  where $c_0^\ecal
>0$. Note that $b_P(z)$ is the same for all of the $P_\ecal$. Then by
the inclusion-exclusion principle
\begin{equation}\label{i-e}\Pi_{|NP}(z,z) = \sum_{\ecal} (-1)^{\#\ecal -1}
\Pi_{|NP_\ecal}(z,z) \qquad (\emptyset \ne \ecal \subset
\{1,\dots,s\})\;.\end{equation}  The desired asymptotic expansion
for $\Pi_{|NP}(z,z)$ now follows from (\ref{ecal})--(\ref{i-e}),
by noting that $\Pi_{|NP}(z,z)\ge \Pi_{|NP_{\{1\}}}(z,z)$ and
hence $c_0^F(z) \ge c_0^{\{1\}}(z) >0$.

The proof of Theorem~\ref{SZEGO} is now complete.

\subsubsection{Proof of Proposition \ref{convergence}.}
\label{s-convergence}

By Theorem~\ref{SZEGO}(ii),
\begin{equation}\label{limit}\frac{1}{N} \log \Pi_{|NP}(z,z)
\to -b_P(z)\;,
\end{equation} for all non-transition points $z$.
We use the log coordinates $\rho_j + i\theta_j=\log z_j$, so that
\begin{equation}\label{uN}u_N(z)=\frac{1}{N}\log
\Pi_{|NP}(z,z) +p \log(1+\|z\|^2)=\frac{1}{N} \log \sum_{\al\in
NP} {\textstyle{Np \choose \al}} e^{2\langle
\al,\rho\rangle}\;,\quad z\in(\C^*)^m\;.\end{equation} Thus for all
non-transition points $z$, we have
\begin{equation}\label{u0}
u_N(z)\to u_\infty(z):= p \log(1+\|z\|^2) -b_P(z)\;.\end{equation}

We must show that convergence of (\ref{u0}) also holds at the
transition points and is uniform on compact subsets of $(\C^*)^m$.
From (\ref{uN}), we obtain
\begin{equation}d_\rho u_N = \frac{1}{N}d_\rho \log \sum_{\al\in NP}
{\textstyle{Np \choose \al}} e^{2\langle \al,\rho\rangle}=
\frac{2\sum_{\al\in NP} {\textstyle{Np \choose \al}}e^{2\langle
\al,\rho\rangle}\al}{N\sum_{\al\in NP} {\textstyle{Np \choose
\al}}e^{2\langle \al,\rho\rangle}}\;. \label{DuN}\end{equation}
Since $\|\al\|\le Np$ for all $\al\in NP$, we therefore have the
uniform upper bound
$$\|d_\rho u_N\|_{\lcal^\infty((\C^*)^m)} \le\frac{\sum_{\al\in NP} {\textstyle{Np \choose \al}}e^{2\langle
\al,\rho\rangle}2\|\al\|}{N\sum_{\al\in NP} {\textstyle{Np \choose
\al}}e^{2\langle \al,\rho\rangle}} \le 2p\;.$$  Since
$\{u_N(z^0)\}$ converges for any non-transition point $z^0$, it
follows that $\{u_N\}$ is uniformly bounded and uniformly
equicontinuous on compact sets.  Therefore it converges uniformly
on compact sets in $(\C^*)^m$.\qed

\medskip

\section{Distribution of zeros}\label{DZ}

\subsection{Expected zero current and the \szego kernel}\label{EDZ}

In order to deduce Theorems \ref{main}--\ref{simultaneous} from
Proposition~\ref{convergence}, we need to relate the expected zero current to  the
conditional \szego kernel
$\Pi_{| P}$.  This relationship is given by the following result:

\begin{prop} \label{E-cond}The expected zero current of  $k$ independent random
polynomials $f_j\in \poly(P_j)$, $1\le j\le k$, is given by
$$\E_{|P_1,\dots,P_k}(Z_{f_1,\dots,f_k}) =\bigwedge_{j=1}^k \E_{|P_j}(Z_{f_j})=
\bigwedge_{j=1}^k\left(
\frac{\sqrt{-1}}{2\pi}
\partial
\bar{\partial} \log \Pi_{| P_j}(z, z)+p_j
\omega_{\FS}\right)\qquad \mbox{\rm on}\ \ (\C^*)^m\;.
$$
\end{prop}

Recall that $\om_\FS=\frac{\sqrt{-1}}{2\pi}\ddbar \log(1+\|z\|)^2$.
 The case $k=1$ of the proposition is
essentially the same as in  \cite[Prop.~3.1]{SZ}, but neither the
statement nor the proof there cover the application we need. We remark that Proposition
\ref{E-cond} is a special case of a more general statement for the zeros of general
linear systems on compact
\kahler manifolds, which we discuss in \cite{SZ3}.

We recall that the space $\dcal'{}^{k,l}(Y)$ of $(k,l)$-currents on an $n$-dimensional
complex manifold
$Y$ consists of the continuous linear functionals on the space $\dcal^{n-k,n-l}(Y)$ of
$\ccal^\infty$ compactly-supported $(n-k,n-l)$-forms on $Y$.  We say that a sequence $
\Psi_j\in \dcal'{}^{k,l}(Y)$ {\it converges weakly\/} to $\Psi_0$ if $(\Psi_j,\phi)\to
(\Psi_0,\phi)$ for all test forms $\phi\in \dcal^{n-k,n-l}(Y)$.

We begin the proof with the following lemma, which covers the $k=1$ case.

\begin{lem}\label{EZ} Let $P$ be a convex polytope in $p\Si$ (or more generally,
an arbitrary subset of $p\Si\cap\Z^m$). Let $Y$ be an algebraic submanifold of
$(\C^*)^m$, and let $h\in\ccal^\infty(Y)$ be given by
$$h(z)= \log \Pi_{|P}(z, z) +p  \log(1+\|z\|^2)\;,\qquad z\in Y\;.$$
Then
\begin{eqnarray*}\E_{|P}(Z_{f|Y})  &=&\frac{\sqrt{-1}}{2\pi}
\partial
\bar{\partial} h \in \dcal^{1,1}(Y) \;.\end{eqnarray*}
\end{lem}

\begin{proof} The proof is essentially the same as
the proof of Proposition~3.1 in \cite{SZ}. Write $f_Y=f|Y$ for
$f\in\poly(P)$. If $f_Y\not\equiv 0$, then the current of
integration over the zeros of $f$ is  given by the
Poincar\'e-Lelong formula:
\begin{equation} Z_{f_Y} =
\frac{\sqrt{-1}}{ \pi } \partial \bar{\partial}\log | f_Y|\in \dcal'{}^{1,1}(Y)\;.
\label{Zs}
\end{equation}

Recalling (\ref{G}), we can write
$$f_Y=\sum_{\al\in P}\frac{c_\al}{\|\chi_\al\|} \chi_\al=\langle c,G\rangle\;,\qquad
c=\big(c_\al\big)_{\al\in P}\;,\quad
G=\left(\frac{1}{\|\chi_\al\|} \chi_\al|_Y\right)_{\al\in P}\;.$$
(The $\lcal^2$ norms $\|\chi_\al\|$ are computed over all of
$\C^m$ as before; see (\ref{IP})--(\ref{IP2}).) We then write
$G(z)= |G(z)| u(z)$ so that $|u| \equiv 1$ and
\begin{equation}\label{logG}\log  | \langle c, G \rangle| = \log |G| +
\log  | \langle c, u \rangle|\;.\end{equation}

By (\ref{Zs})--(\ref{logG}), we have for any test form
$\phi\in\dcal^{n-1,n-1}(Y)$ (where $n=\dim Y$),
\begin{eqnarray*}\big(\E_{|P}(Z_{f_Y}),\phi\big)&=&\frac{\sqrt{-1}}{ \pi}
\left(\int_{\C^{k}}\log  | \langle c,G \rangle| \,d\gamma(c)
,\ddbar \phi\right)\\
&=&\frac{\sqrt{-1}}{ \pi} \left(\partial \bar{\partial}
\int_{\C^{k}}\log  |G| \, d\gamma(c), \phi\right)
+\frac{\sqrt{-1}}{ \pi}\left(\ddbar \int_{\C^{k}} \log  | \langle
c, u \rangle|\, d\gamma(c) , \;\phi\right)\;,
\end{eqnarray*}
where $k=\#P$ and $d\ga(c)=\frac 1 {\pi^k}e^{-|c|^2}\,dc$. Upon
integration in $c$, the second term becomes constant in $z$ and
$\ddbar$ kills it. The first term is independent of $c$ so we may
remove the Gaussian integral.  Thus
\begin{equation}\label{Egamma}\E_{|P}(Z_{f_Y})=\frac{\sqrt{-1}}{2 \pi} \partial
\bar{\partial}\log  |G|^2  \;.\end{equation} Recalling
(\ref{normchi}) and (\ref{Sz2e}), we have $$\Pi_{|P}(z,z)=
(1+\|z\|^2)^{-p}|G(z)|^2\;,\qquad z\in Y\;,$$ and the formula of
the lemma follows.
\end{proof}

\medskip\noindent {\it Proof of Proposition \ref{E-cond}:\/}  We first note that the
expected current is well defined, since for almost all choices of
the $f_j$, $$\int_{|Z_{f_1,\dots,f_k}|}\om_\FS^{m-k}\le p_1\cdots
p_k\;,$$ and hence for each test form
$\phi\in\dcal^{m-k,m-k}((\C^*)^m)$, the function
$(f_1,\dots,f_k)\mapsto (Z_{f_1,\dots,f_k},\phi)$ is in
$\lcal^\infty(\poly(P_1)\times\cdots\times \poly(P_k))$.

We shall verify the current identity by induction on $k$.  Proposition \ref{EZ} with
$Y=(\C^*)^m$ gives the case  $k=1$, so assume that $k>1$ and the proposition has been
verified for
$k-1$ polynomials.  By the inductive assumption, it suffices to show that
\begin{equation} \label{ind} \E_{|P_1,\dots,P_k}(Z_{f_1,\cdots, f_k}) =
\E_{|P_1,\dots,P_{k-1}}(Z_{f_1,\cdots, f_{k-1}})
\wedge \E_{|P_k}(Z_{f_k})\;.\end{equation}

 Write $Y=
|Z_{f_1,\dots,f_{k-1}}|$. By Bertini's Theorem, we know that the
$Z_{f_j}$ are smooth and intersect transversally in $(\C^*)^m$ and
hence $Y$ is smooth (and of codimension $k-1$) in $(\C^*)^m$ for
almost all $f_1,\dots,f_k$.  Therefore for
$\phi\in\dcal^{m-k,m-k}((\C^*)^m)$, we have
\begin{equation}\label{Bertini}(Z_{f_1,\dots,{f_k}},\phi )=
\int_{Z_{f_1,\dots,{f_k}}}\phi =\int_{Y\cap Z_{f_k}}\phi =
(Z_{f_k|Y},\phi|_Y)\;,\end{equation} for almost all $f_1,\dots,f_k$.
Averaging (\ref{Bertini}) over $f_k$ and applying Lemma \ref{EZ}, we obtain
\begin{equation}\label{Bertini2}\int_{\poly(P_k)} (Z_{f_1,\dots,{f_k}},\phi )\,
d\ga_{p_k|P_k}(f_k)= \E_{|P_k}(Z_{f_k|Y},\phi|_Y)=\int_Y
\E_{|P_k}(Z_{f_k})\wedge \phi
=(Z_{f_1,\dots,f_{k-1}},\,\E_{|P_k}(Z_{f_k})\wedge \phi)
\;,\end{equation} where $$
\E_{|P_k}(Z_{f_k})=\frac{\sqrt{-1}}{2\pi}\log \Pi_{|P_k}(z, z)
+p\om_\FS\in\dcal^{1,1}((\C^*)^m)\;.$$ Averaging (\ref{Bertini})
over $f_1,\dots,f_{k-1}$ and applying (\ref{Bertini2}), we obtain
$$\E_{|P_1,\dots,P_k}(Z_{f_1,\cdots,
f_k},\phi)=\E_{|P_1,\dots,P_{k-1}} \big(Z_{f_1,\dots,f_{k-1}},
\,\E_{|P_k}(Z_{f_k})\wedge \phi\big) =(\E_{|P_1,\dots,P_{k-1}}
(Z_{f_1,\cdots, f_{k-1}}) \wedge \E_{|P_k}(Z_{f_k}) ,\,\phi)\;.$$
\qed

\subsection{Asymptotics of zeros: proofs of Theorems \ref{main} and
\ref{simultaneous}}\label{s-proofs}
To prove Theorem \ref{main}, we let $u_N,\ u_\infty$ be as in the
proof of Proposition~\ref{convergence}. By Proposition
\ref{E-cond} (recalling that $\om_\FS=
\frac{\sqrt{-1}}{2\pi}\ddbar \log (1+\|z\|^2)$), we have
\begin{equation}\label{ddbaruN}\frac{\sqrt{-1}}{2\pi}\d\dbar
u_N=\frac{1}{N}\E_{|NP}(Z_f) \qquad \mbox{on }\ (\C^*)^m
\;,\end{equation} and hence $u_N$ is plurisubharmonic.  By
(\ref{u0}) and (\ref{ddbaruN}), we conclude that
\begin{eqnarray} \frac{1}{N}\E_{|NP}(Z_f) \to  \frac{\sqrt{-1}}{2\pi}
\ddbar u_\infty=p\om_\FS - \frac{\sqrt{-1}}{2\pi} \ddbar
b_P\;,\label{psi0}\end{eqnarray} where differentiation is in the
sense of currents (see e.g., \cite[Chapter~I]{Shabat}).

We now show that the current $\ddbar b_P\in
\dcal'{}^{1,1}((\C^*)^m)$ is given by a $(1,1)$-form with
piecewise smooth coefficients; in fact $\ddbar b_P$ is
$\ccal^\infty$ on each of the regions ${\rcal_F}$ given by
(\ref{RF}). (Equivalently, the Radon measure $\ddbar b_P\wedge
\om_\FS$ does not charge the set of transition points, and hence
formula (\ref{psi0}) can be interpreted as differentiation in the
ordinary sense on the regions $\rcal_F$.) By Theorem~\ref{SZEGO},
$b_P\in\ccal^1((\C^*)^m)$; i.e., if $z^0\in\d{\rcal_F}\cap
\d\rcal_{F'}$, then the values of $d b_P(z^0)$ computed in the two
regions $\overline{\rcal_F}$ and $\overline{\rcal_{F'}}$ agree.
Then for a test form $\phi\in\dcal^{m-1,m-1}((\C^*)^m)$, we have
\begin{equation}\label{ddbarb1}(\ddbar b_P, \phi) = \sum_F\int_{\rcal_F}b_P\d\dbar
\phi =\sum_F\int_{\rcal_F^\circ} \ddbar b_P\wedge \phi -
\sum_F\int_{\d\rcal_F}(\dbar b_P\wedge \phi
+b_P\wedge\d\phi)\;.\end{equation} Note that $\d\rcal_F$ consists
of $\ccal^\infty$ real hypersurfaces (consisting of those points
of $\d{\rcal_F}$ that are contained in the boundary of only one
other region $\rcal_{F'}$) together with submanifolds of real
codimension $\ge 2$, and hence  Stokes' Theorem applies (see e.g.
\cite[4.2.14]{Fe}). Since $b_P$ and $\dbar b_P$ are continuous on
$(\C^*)^m$, the boundary integral terms  in (\ref{ddbarb1}) cancel
out and we obtain
\begin{equation}\label{ddbarb2}(\ddbar b_P, \phi)  =\sum_F\int_{\rcal_F^\circ}
\ddbar b_P\wedge \phi  =\int_{(\C^*)^m\sm E}\ddbar b_P\wedge
\phi\;,\end{equation} where $E:=\bigcup_F \d\rcal_F$ is the set of
transition points. Formula (\ref{ddbarb2}) says that  the current
$\ddbar b_P$ is  a $(1,1)$-form with piecewise smooth coefficients
obtained by differentiating $b_P$ on the regions $\rcal_F^\circ$.

We now let
\begin{equation}\label{psi} \psi_P=\frac{\sqrt{-1}}{2\pi}
\ddbar u_\infty=p\om_\FS - \frac{\sqrt{-1}}{2\pi} \ddbar
b_P\end{equation} on each of the regions $\rcal_F$. By
(\ref{ddbarb2}), the current $\psi_P$ is also a piecewise smooth
$(1,1)$-form; by (\ref{psi0})
\begin{equation}\label{weaklim0}N\inv\E_{|NP}(Z_f)
\to \psi_P\qquad \mbox {(weakly)}\;.\end{equation}

We shall show $\lcal^1_{\rm loc}$ convergence of (\ref{weaklim0})
when we prove Theorem \ref{simultaneous} below.  Continuing with
the proof of Theorem \ref{main}, we observe that (ii) is an
immediate consequence of (\ref{psi}), since $b_P=0$ on the
classically allowed region.

To  verify (iii), we again use the log coordinates $\zeta_j=\rho_j
+i\theta_j=\log z_j$, so that
\begin{equation}\label{uinfty}u_\infty=p\log\left(1+\sum
e^{2\rho_j}\right)-b_P(e^\rho)\;.\end{equation} Since $u_\infty$
is independent of the angle variables $\theta_j$, we have
\begin{equation}\label{d2rho}\psi_P=\frac{\sqrt{-1}}{2\pi}\ddbar u_\infty =
\frac{\sqrt{-1}}{8\pi}\sum_{jk} \frac{\d^2
u_\infty}{\d\rho_j\d\rho_k} d\zeta_j\wedge d\bar\zeta_k\ge
0\;.\end{equation}  Thus we must show that the Hessian of
$u_\infty(\rho)$ has rank $r$ at points $z^0\in\rcal^\circ_F$,
where $r=\dim F$. From (\ref{Db}) and (\ref{tauphase}), we have
\begin{equation} d_\rho b_P =-p d_\rho \log
\left.\left(1+\sum
e^{-\tau_j+2\rho_j}\right)\right|_{\tau=\tau_z}+ pd_\rho\log
\left(1+\sum e^{2\rho_j}\right) \label{Db1}\;.\end{equation} To
simplify (\ref{Db1}), we note that for $\sigma\in\R^m$ we have
\begin{equation}\label{dlog} d_\rho\log\left(1+\sum
e^{-\sigma_j+2\rho_j}\right) = \left(
\frac{2e^{-\sigma_1+2\rho_1}}{1 + \sum e^{-\sigma_j+2\rho_j}}
,\dots, \frac{2e^{-\sigma_1+2\rho_1}}{1 + \sum
e^{-\sigma_j+2\rho_j}}\right)=2\mu(e^{-\sigma/2}\cdot z)
\;.\end{equation}

Hence by (\ref{dlog}) with $\sigma =\tau_z$ and $0$, we have
\begin{equation}\label{Db2} d_\rho b_P = 2p(\mu -
\mu\circ\xi)\;,\end{equation} where
$\xi=(\xi_1,\dots,\xi_m):(\C^*)^m\to (\C^*)^m$ denotes the map given
by
\begin{equation}\label{xi}\xi(z)=e^{-\tau_z/2}\cdot z= e^{-\tau_z/2
+\rho+i\theta}\;.\end{equation} By (\ref{uinfty}), we then obtain
\begin{eqnarray} d_\rho u_\infty&=&2p\,\mu\circ\xi\;.
 \label{dbaru}\end{eqnarray}

Hence the Hessian $\hcal_\rho u_\infty$ equals the Jacobian of
$2p\mu\circ\xi$.  Since $\mu \circ \xi (z)=\frac{1}{p} q(z)$,
we easily see that
\begin{equation}\label{submersion}\mu\circ\xi:\rcal^\circ_F
\to \frac 1p F\;,\end{equation} so rank$\,D(\mu\circ\xi) \le r$.  In fact,
(\ref{submersion}) is a submersion, so that the rank equals $r$.
To see that it is a submersion, we recall from the proof of Lemma
\ref{claim2-3} that the `lipeomorphism' $\Phi\inv:\Si^\circ\to
\ncal^\circ$ restricts to a diffeomorphism
$$\textstyle\Si^\circ\supset \mu(\rcal^\circ_F) \buildrel{\approx}\over \to\frac 1p
F\times C_F\subset \ncal^\circ\;, \qquad \mu(z)\mapsto
\left(\frac{1}{p}q(z),\tau_z\right)
=(\mu\circ\xi(z),\tau_z)\;,$$ and therefore (\ref{submersion}) is
a submersion, completing the proof of (iii).

We now prove Theorem \ref{simultaneous}; the case $k=1$ of the
theorem will then yield part  (i) of Theorem~ \ref{main}. Let
$P_1,\dots, P_k$ be convex integral polytopes, as in Theorem
\ref{simultaneous}. We first show that
\begin{equation}\label{weaklim} N^{-k}\E_{|NP}(Z_{f_1,\dots,f_k})
\to \psi_k:=\psi_{P_1}\wedge \cdots\wedge \psi_{P_k}\qquad \mbox
{weakly}.\end{equation} For $1\le j\le k$, we let
$$u_N^j=\frac{1}{N}\log \Pi_{|NP_j}(z,z) +p_j \log(1+\|z\|^2)\;,$$
so that, recalling (\ref{u0}),
$$u_N^j(z)\to u_\infty^j(z):= p_j
\log(1+\|z\|^2) -b_{P_j}(z)\;.$$ By (\ref{psi}),
$$\frac{\sqrt{-1}}{2\pi} \ddbar u_\infty^j = \psi_{P_j}\;.$$

Recalling Proposition \ref{E-cond}, we introduce the $(k,k)$-forms
\begin{equation}\label{w1}\kappa_N:=\left(\frac{\sqrt{-1}}{2\pi}
\ddbar u_N^1\right) \wedge \cdots \wedge
\left(\frac{\sqrt{-1}}{2\pi} \ddbar u_N^k\right)
=N^{-k}\E_{|NP_1,\dots,NP_k}(Z_{f_1,\dots,f_k})\;.\end{equation}
Since $u_N^j\to u_\infty^j$ locally uniformly, it follows from the
Bedford-Taylor Theorem \cite{BT,Kl} that
\begin{equation}\label{w2}\kappa_N \to \left(\frac{\sqrt{-1}}{2\pi}
\ddbar u_\infty^1\right) \wedge \cdots \wedge
\left(\frac{\sqrt{-1}}{2\pi} \ddbar u_\infty^k\right) \qquad \mbox
{weakly}.\end{equation}  In fact, the limit current in
(\ref{w2}) is absolutely continuous, and hence is equal to the
locally bounded (piecewise smooth) $(k,k)$-form $\psi_k$.  To see
this, we note that by the Bedford-Taylor Theorem,
$$\bigwedge_{j=1}^k(\psi_{P_j} *\phi_\epsilon)=\bigwedge_{j=1}^k
\left[\frac{\sqrt{-1}}{2\pi}\ddbar (u_\infty^j *
\phi_\epsilon)\right] \to
\bigwedge_{j=1}^k\left(\frac{\sqrt{-1}}{2\pi}\ddbar
u_\infty^j\right)\;,$$ where $\phi_\epsilon$ denotes an
approximate identity. Since the currents $\psi_{P_j}$ have
coefficients in $\lcal^\infty_{\rm{loc}}$, the forms
$\bigwedge_{j=1}^k(\psi_P *\rho_\epsilon)$ have locally uniformly
bounded coefficients; absolute continuity of the limit current
follows. The weak limit (\ref{weaklim}) now follows from
(\ref{w1})--(\ref{w2}).

To complete the proof of Theorem \ref{simultaneous}, we must show
that $\kappa_N \to \psi_k$ in $\lcal^1(K)$ for all compact
$K\subset(\C^*)^m$.  Let $\ep>0$ be arbitrary, and choose a
nonnegative function $\eta\in \dcal((\C^*)^m)$ such that
$\eta\equiv 1$ on a neighborhood of $(E_1\cup\cdots\cup E_k)\cap
K$, where $E_j$ is the set of transition points for $P_j$, and
$$\int_{(\C^*)^m}
\psi_k\wedge \eta\om_\FS^{m-k}<\ep\;.$$ By
(\ref{weaklim})--(\ref{w1}),
$$\int_{(\C^*)^m}
\kappa_N\wedge \eta\om_\FS^{m-k}<2\ep\qquad \mbox{for }\ N\gg
0\;.$$ (Note that the above integrands are nonnegative.)
Therefore,
$$\left\|\kappa_N - \psi_k\right\|_{\lcal^1(K)}  \le
\left\|(1-\eta)(\kappa_N - \psi_k)\right\|_{\lcal^1(K)} +
\left\|\eta\kappa_N\right\|_{\lcal^1(K)}+
\left\|\eta\psi_k\right\|_{\lcal^1(K)}\;.
$$ Since $\kappa_N$ is a positive $(k,k)$-form,
$$\left\|\eta\kappa_N\right\|_{\lcal^1(K)}=
\int_K\kappa_N\wedge \eta\om_\FS^{m-k}<2\ep \;.$$ Similarly,
$\left\|\eta\psi_k\right\|_{\lcal^1(K)}=\int_K \psi_k\wedge
\eta\om_\FS^{m-k}<\ep$. Since
$$(1-\eta)(\kappa_N - \psi_k)\to 0\qquad \mbox{uniformly on }\ K\;,$$
it follows that
$$\left\|\kappa_N - \psi_k\right\|_{\lcal^1(K)} \to 0\;,$$
completing the proof of Theorems \ref{main} and
\ref{simultaneous}.\qed

\bigskip
We now derive our integral formula (\ref{b-action}) from
(\ref{Db2}):
\begin{prop}  The decay function $b_P$ of Theorem \ref{MASS} is given by
$$ b_P(z) =
\int_0^{\tau_z}\left[-q(e^{-\sigma/2}\cdot z) +
p\mu(e^{-\sigma/2}\cdot z)\right] \cdot d\sigma\;.$$
\label{b-prop}\end{prop}
\begin{proof}
Fix a point $z$ in the classically forbidden region. Then by
(\ref{Db2})  with  the change of variables $\rho_j=- \sigma_j/2 +
\log |z_j|$, we have
\begin{equation}b_P(e^{-\tau_z/2}\cdot z)) - b_P(z) = \int_0^{\tau_z} d_\sigma
b_P(e^{-\sigma/2}\cdot z) \cdot d\sigma = p \int_0^{\tau_z}
\big[\mu\circ\xi(e^{-\sigma/2}\cdot z) -\mu(e^{-\sigma/2}\cdot z)
 \big] \cdot
d\sigma\;,\label{FTC}\end{equation} where the map $\xi$ is given
by (\ref{xi}). By definition, $q(z)=p\mu\circ\xi(z)$, and hence
$p\mu\circ\xi(e^{-\sigma/2}\cdot z) = q(e^{-\sigma/2}\cdot z)$.
Furthermore $b_P(e^{-\tau_z/2}\cdot z)= 0$ since
$e^{-\tau_z/2}\cdot z=\xi(z)\in \d \acal_P$, and formula
(\ref{b-action}) then follows from (\ref{FTC}).\end{proof}

We note that the proof of Theorem \ref{simultaneous} also gives an asymptotic expansion
away from transition points:

\begin{theo}\label{expansion} Let $P_1,\dots,P_k$ be convex integral polytopes.  Let
$U$ be a relatively compact domain in $(\C^*)^m$ such that
$\overline U$ does not contain transition points for any of the
$P_j$. Then we have a complete asymptotic expansion of the form
$$\frac{1}{N^{k}}\E_{|NP_1,\dots,NP_k} (Z_{f_1, \dots, f_k}) \sim
\psi_{P_1}\wedge\cdots\wedge\psi_{P_k} +\frac{\phi_1}{N} + \cdots
+ \frac{\phi_n}{ N^n} +\cdots \qquad \mbox {on } \ U\;,$$ with
uniform $\ccal^\infty$ remainder estimates, where the $\phi_j$ are
smooth $(k,k)$-forms on $U$.
\end{theo}

\begin{proof} By (\ref{uN}) and (\ref{w1}) we have
$$\frac{1}{N^{k}}\E_{|NP_1,\dots,NP_k} (Z_{f_1, \dots, f_k})=
\bigwedge_{j=1}^k\frac{\sqrt{-1}}{2\pi} \ddbar
\left[\frac{1}{N}\log \Pi_{|NP_j}(z,z) +p_j
\log(1+\|z\|^2)\right]\;.$$ The conclusion now follows from the
asymptotic expansion of Theorem~\ref{SZEGO}.\end{proof}

\subsection{Proof of Theorem \ref{probK}}  Theorem \ref{probK} is
an easy consequence of our prior results.  Indeed, by
Proposition~\ref{E-cond}, we can write
$N^{-m}
\E_{|NP}(Z_{f_1,\dots,f_m}) = G_N\om^m$, where the $G_N$ are positive $\ccal^\infty$
functions on $(\C^*)^m$. We must show the convergence of the sequence of
measures
$$d\la_N:=N^{-m}
\E_{|NP}(Z_{f_1,\dots,f_m}) = G_N\om^m\;.$$
Now let $B$ be a Borel subset of $(\C^*)^m$. Then
$$\lambda_N(B) =N^{-m}\E_{|NP}\big(\#\{z\in B:
f_1(z)=\cdots=f_m(z)=0\}\big)\;.$$ We first consider
the case where
$B$ is  contained in a compact set  $K\subset(\C^*)^m$.  Then by
Theorem~\ref{simultaneous} with
$k=m$ and
$P_1=\cdots=P_k=P$, it follows that
$$G_N\om^m \to\psi_P^m \qquad \mbox{in }\ \lcal^1(K)\;.$$
By Theorem \ref{main}, $\psi_P^m=p^m\om^m_\FS$ on $\acal_P$, and $\psi_P^m=0$ on the
complement of $\acal_P$ since its rank is less than $m$ there.  Therefore,
\begin{equation}\la_N(B)= \int _K \chi_B G_N\om^m \to  \int_K
\chi_{B\cap\acal_P}\,p^m\om^m_\FS = m!p^m\,\vol_{\CP^m}(B\cap
\acal_P)\;.\label{lim}\end{equation}

For the general case, let $\ep>0$ be arbitrary and choose
$B'\subset B$ such that
$B'$ is relatively compact in $(\C^*)^m$ and $\vol_{\CP^m}(B\sm B') <\ep$. Then by the
above
$$\liminf_{N\to\infty} \la_N(B) \ge \liminf_{N\to\infty}
\la_N(B')
=m!p^m\,\vol_{\CP^m}(B'\cap\acal_P)
\ge m!p^m\,\vol_{\CP^m}(B\cap
\acal_P) -m!p^m\ep\;.$$
Since $\ep>0$ is arbitrary, we conclude that
\begin{equation}\label{liminf} \liminf_{N\to\infty}\la_N(B)\ge m!p^m\,\vol_{\CP^m}(B\cap
\acal_P)\;.\end{equation}
To obtain the reverse inequality, we recall from the Bernstein-Kouchnirenko Theorem
\cite{Be,Ku1,Ku2} that the number of common zeros of $\{f_1,\dots,f_m\}$ equals
$m!N^m\vol(P)$ for almost all $(f_1,\dots,f_m)\in\poly(P)^m$, and hence
\begin{equation}\label{BK}
 \la_N((\C^*)^m) = m!\vol(P)\qquad \mbox{for all }\ N\ge 1\;.\end{equation}
It is well-known and easy to verify that $\vol_{\CP^m}(\mu\inv(\Omega))=\vol(\Omega)$
for any $\Omega\subset\Si$. (This is a special case of a general fact about the moment
map in symplectic geometry; see (\ref{vol-deg}) in the Appendix.) Recalling that
$\acal_P=\mu\inv(\frac 1p P^\circ)$, we then have
\begin{equation}\label{totvol} \vol(\acal_P)=\frac {1}{p^m}\vol(P)\;.\end{equation}
It follows from (\ref{liminf}) with $B$ replaced by its complement $B^c=(\C^*)^m\sm B$
and (\ref{BK})--(\ref{totvol}) that
\begin{equation}\label{limsup} \limsup_{N\to\infty} \la_N(B) = m!\vol(P) -
\liminf_{N\to\infty} \la_N(B^c)\le m!p^m\,\vol_{\CP^m}(B\cap
\acal_P)\;.\end{equation}
Therefore, (\ref{lim}) holds for a general Borel set $B\subset (\C^*)^m$.\qed

\medskip
\subsection{Vanishing of $\psi_P$ along the normal flow}
The proof of Theorem \ref{main}(iii) gives us some more
information about the expected zero current in the classically
forbidden region.   Roughly speaking, let $f\in \poly(NP)$ be a
random polynomial with Newton polytope $NP$, with $N$ large, and
let $z^0\in |Z_f|\sm\acal_P$; then $|Z_f|$ is highly likely to be
close to being tangent to the orbit of the normal flow at $z^0$.
In particular, if $z^0$ is in the flow-out of an edge (dimension 1
face) of $P$, then $T_{z^0}(|Z_f|)$ is likely to be a good
approximation  to the tangent space of the normal flow through
$z^0$.

To make this statement precise, we define the {\it complexified
normal cone\/} of a face $F$,
$$\wt C_F := \{\tau+i\theta: \tau \in C_F,\ \theta\in T_F^\perp\}\;.$$  (Recall
that $C_F\subset T_F^\perp$.)
 We note that $\wt C_F$ is a semi-group, which
 acts on
$\rcal_F$ by the rule $\eta(z )= e^\eta \cdot z$; we call this
action the `(joint) normal flow.'  The (maximal) orbits of the
normal flow are of the form $\wt C_F\cdot z^0=\{e^{\eta}\cdot
z^0:\eta\in\wt C_F\}$, where  $z^0\in \mu\inv(\frac 1p F)$. We note that
the orbit $\wt C_F\cdot z^0$ is a complex $(m-r)$-dimensional
submanifold (with boundary) of $(\C^*)^m$.  (Indeed, $(\wt C_F\cdot
z^0)\cap \rcal^\circ_F$ is a submanifold without boundary in
$\rcal^\circ_F$.)

\begin{theo}\label{more} Let $P$ be a convex integral polytope and let $\psi_P$ be
the limit expected zero current of Theorem \ref{main}. Then
$\psi_P$ vanishes along the orbits of the normal flow.
\end{theo}

\begin{proof} Let $$O:=\wt C_F\cdot z^0=\{e^{\tau+i\theta}\cdot z^0: \tau
\in C_F,\ \theta\in T_F^\perp\}$$ be a maximal orbit of the normal
flow, where $p\mu(z^0)\in F$. For $z= e^{\tau+i\theta}\cdot z^0\in
O$, we have
$$\mu\circ\xi (z)= \frac{1}{p}q(z)=
\frac{1}{p}q(e^{\tau}\cdot z^0) = \frac{1}{p}q(z^0)
=\mu(z^0)$$ and hence $\mu\circ\xi$ is constant on $O$. It then
follows from (\ref{d2rho}) and (\ref{dbaru}) that $\psi_P|_O=0$.
\end{proof}

To relate Theorem \ref{more} to the explanation above, suppose for
example that $z^0$ is in the flow-out of a 1-dimensional face
(edge). Choose coordinates $w_1,\dots,w_m$  so that the normal
flow near $z^0$ is given by $w_1=\ $constant. Since $\psi_P$
vanishes along the orbits of the flow and has rank $1$ near $z^0$,
we have \begin{equation}\label{relate}\frac 1N \E_{|NP}
(Z_f)\to\psi_P =
 ic(w) dw_1\wedge d\bar
w_1\;,\qquad c(w)>0,\end{equation} near $z^0$. For a regular point
$z$ of $|Z_f|$, we let $\eta_f(z)\in T^{*1,0}((\C^*)^m)$ be a unit
vector (unique up to the $S^1$ action) annihilating
$T^{1,0}(|Z_f|)$.  We can then write $Z_f= \delta_{Z_f} (\frac i 2
\eta\wedge \bar \eta)\;,$ where $\delta_{Z_f}$ is the measure
given by $(\delta_{Z_f},\phi) = \int_{|Z_f|}\phi\, d\vol_{2m-2}$.
Writing $\eta= \sum_{j=1}^m a_j(w) dw_j$, we see by (\ref{relate})
that
$$\frac 1N \E_{|NP} \int_{|Z_f|}|a_j|^2\, d\vol_{2m-2} \to 0, \quad
\mbox{for }\ j\ge 2.$$ Thus the expected value of the average
distance between the tangent spaces of $|Z_f|$ and of the normal
flow approaches 0 as $N\to\infty$.

\subsection{Amoebas in the plane}\label{AM}  The term `amoeba' was introduced by Gelfand,
Kapranov and  Zelevinsky \cite{GKZ} to refer to the image under
the moment map of a zero set $Z_{f_1,\dots,f_k}$ of polynomials,
and have been studied in various contexts (see
\cite{FPT,GKZ,M1,PR} and the references in the survey article by
Mikhalkin \cite{M}).  The image of a zero set under the moment map
$\mu$ is called a {\it compact amoeba\/}, while the image under
the map
\begin{equation}\label{LOG}\Log:(\C^*)^m\to\R^m, \qquad (z_1,\dots,z_m)\mapsto
(\log|z_1|,\dots,
\log|z_m|)\end{equation} is
 a {\it noncompact amoeba\/}, or simply an {\it amoeba\/}.  Note that $\Log$ is the
moment map for the $\T$ action with respect to the Euclidean
symplectic form $\sum dx_j\wedge dy_j$, and $\Log= \half
\lcal\circ \mu$, where $\lcal:\Si^\circ\approx \R^m$ is the
diffeomorphism given by (\ref{lcal}).

To illustrate what our statistical results can say about amoebas,
we consider zero sets in $(\C^*)^2$. An amoeba in $\R^2$ is the
image of a plane algebraic curve under the Euclidean moment map
Log. An example of an amoeba of the form $\Log(Z_f)$, where $f$ is
a quartic polynomial in two variables with (full) Newton polytope
$4\Si$, is given in the illustration from \cite{Th} reproduced in
Figure~\ref{fig-amoeba} below.

\begin{figure}[htb]
\centerline{\includegraphics*[bb=194 274 417 469]{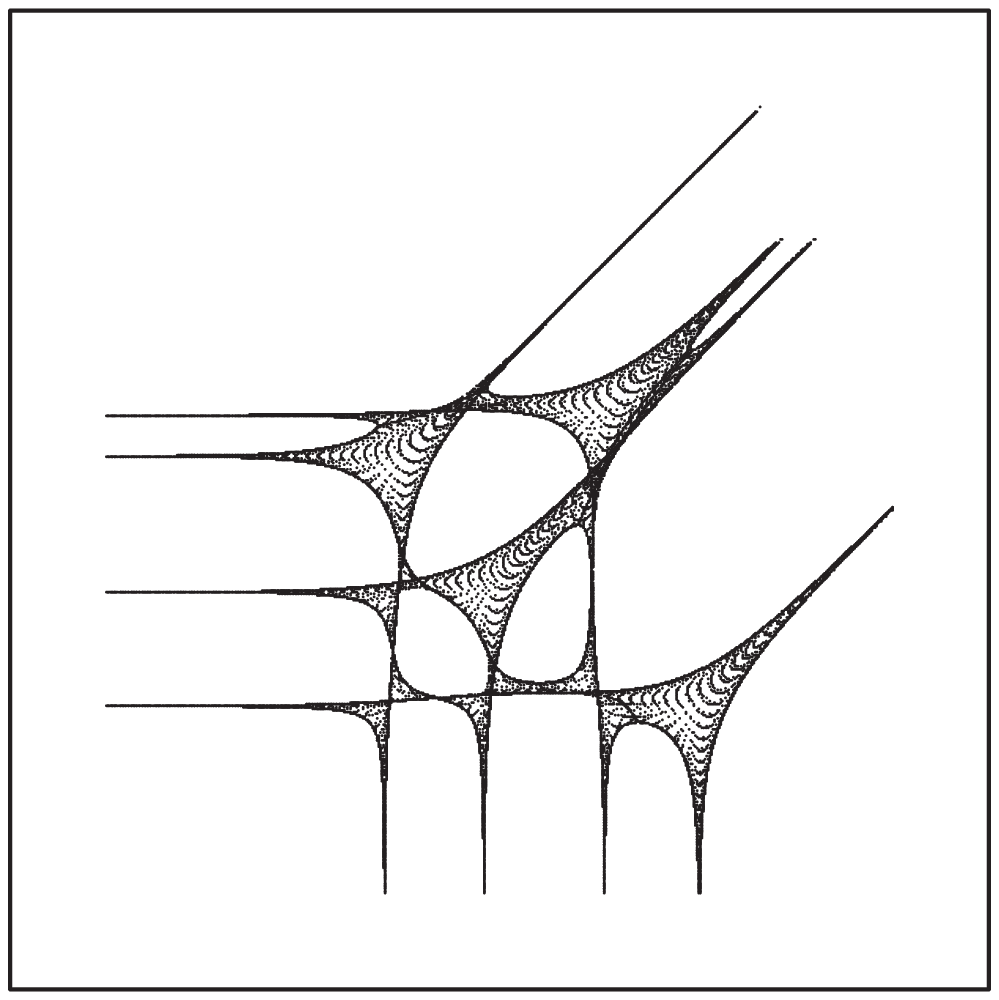}}
\caption{An amoeba with polytope $4\Si$}
\label{fig-amoeba}\end{figure}

One notices that this amoeba contains 12 `tentacles'. By
definition, a tentacle on a compact amoeba $A$ is a connected
component of a small neighborhood in $A$ of $A\cap\d\Si$; the
tentacles of a noncompact amoeba correspond to those of the
compact amoeba under the diffeomorphism  $\Si^\circ\approx \R^m$.

It is known that there is a natural injective map from the set of connected
components (which are convex sets) of the complement of a
noncompact amoeba $A$ to the set $P\cap\Z^2$ of lattice points of
the polytope, and that there are amoebas for which each lattice
point is assigned to a component of the complement. (This fact is
also valid in higher dimensions; see \cite{FPT, M1}.)  For a
generic 2-dimensional noncompact amoeba with polytope $P$, each
lattice point in $\d P$ corresponds to a distinct  unbounded
component of $\R^2\sm A$, and adjacent lattice points correspond
to adjacent unbounded components. (The correspondence is given in
\cite[\S 3.1]{M1} or \cite{FPT}.) Each tentacle thus corresponds to a segment
connecting 2 adjacent lattice points on $\d P$. Hence the number of tentacles of
$A$ equals the number of points of $\d P\cap\Z$; this number is
called the {\it length\/} of $\d P$. (For example, the 12
tentacles in Figure \ref{fig-amoeba} correspond to the 12 segments connecting the 12
lattice points in
$\d(4\Si)$.)

We can decompose $\d P$ into two pieces: $\d^\circ P=\d P\cap
p\Si^\circ$ and $\d^e P=P\cap \d (p\Si)$. Tentacles corresponding
to segments of $\d^\circ P$ end (in the compact picture $\Si$) at
a vertex of $\Si$, and tentacles corresponding to segments of
$\d^e P$ are free to end anywhere on the face of $\Si$ containing
the segment.  We call the latter {\it free tentacles\/}, and we
say that a free tentacle is a {\it classically allowed tentacle\/}
if its end is in the classically allowed region $\acal_P$. For an
amoeba $A$, we let $\nu_{\rm AT}(A)$ denote the number of
classically allowed tentacles of $A$. It is clear from the above
that
$$\nu_{\rm AT}(A) \le \#\{\mbox{free tentacles}\}=\mbox{Length}(\d^e P)$$ and
that this bound can be attained for any polytope $P$.  Here,
`Length' means the length in the above sense; i.e., the diagonal
face of $p\Si$ is scaled to have length $p$. As a consequence of
Theorem \ref{main} (for $m=1$), we conclude that this maximum is
asymptotically the average:

\begin{cor}\label{amoeba} For a convex integral polytope $P$ in $\R^2$, we have
$$\frac{1}{N}\E_{|NP}\Big(\nu_{\rm AT}\big(\Log
(Z_f)\big)\,\Big) \to \mbox{\rm Length}(\d^e P)\;.$$
\end{cor}

\begin{proof} Let $F_1,F_2,F_3$ denote the facets of $p\Si$ and apply Theorem
\ref{main}(ii) to the 1-dimensional polytopes $P\cap \bar F_j$,
$j=1,2,3$.\end{proof}

\section{Appendix: approach using toric geometry} \label{Appendix}
Although our results and proofs do not depend on the theory of toric varieties, the
subject motivated some of our  ideas. We describe briefly in this appendix how toric
varieties can be used to give a geometric derivation of the \szego kernel asymptotics of
Theorem~\ref{SZEGO}.

\subsection{Relationship to toric geometry}  A {\it toric variety\/} is a complex
algebraic variety $M$ containing the complex torus $(\C^*)^m$ as a
Zariski-dense open set such that the group action of $(\C^*)^m$
extends to a $(\C^*)^m$ action on $M$. A toric variety can be
constructed from a fan by gluing together the affine varieties
arising from the cones in the fan (see \cite[Chapter~1]{F}). The
toric variety $M_P$ constructed in this way from the fan of a
convex integral polytope $P$ is a projective variety.  If $P$ is
simple, then $M_P$ has orbifold singularities; if $P$ is Delzant,
then $M_P$ is smooth \cite{De} (see also \cite[\S2.1--2.2]{F}).

For example, $\CP^m$ is the toric variety corresponding to the
simplex $\Si\subset\R^m$. Hence, the toric variety corresponding
to the square from Example~1 in \S \ref{examples} is the product
$\CP^1\times \CP^1$ of projective lines. When $P$ is the trapezoid
in Example~2, the associated toric variety $M_P$ is the blow up of
$(0,0,1)\in\CP^2$, i.e., $M_P$ is the Hirzebruch surface $F_1$.
The toric variety corresponding to the polytope in Example~3 is
the Hirzebruch surface $F_n$ (see \cite[\S 1.1]{F}).

We shall now assume for simplicity that $P$ is Delzant  and
has nonempty interior. In this case, $M_P$ can be
given as the closure of the image of  a monomial embedding
$\Phi_P:(\C^*)^m\to\CP^{\#P-1}$ (see
\cite[Chapter~5]{GKZ} or
\cite{STZ1}). We also give $M_P$ the structure of a symplectic or \kahler
manifold with symplectic/\kahler form $\om_P:=\Phi_P^* \om_\FS$.  (The symplectic form
$\om_P$ depends on the choice of constants defining the monomial embedding $\Phi_P$;
see \cite{STZ1}.) We define the line
bundle
$L_P := \Phi_P^* {\mathcal O}(1)$,
where ${\mathcal O}(1)$ denotes the hyperplane section
 bundle on $\CP^{\#P-1}$, and we give $L_P$ the Hermitian metric obtained by
pulling back the Fubini-Study metric on ${\mathcal O}(1)$ so that $L_P$ has
curvature form $\om_P$.

For example, for the case $P=\Si$, the toric variety $M_\Si=\CP^m$,  $L_\Si=\ocal(1)$,
and $$\poly(N\Si)\cong H^0(\CP^m,\ocal(N))= H^0(M_\Si,L_\Si^N)\;.$$
(Recall that $H^0(M,L)$ denotes the space of holomorphic sections of $L$.)
More generally,
\begin{equation}\label{toric} \poly(NP) \simeq H^0(M_P, L_P^N) =\Phi_P^{*} H^0(\CP^{\#P
-1}, {\mathcal O}(1))\;. \end{equation}

The underlying real torus ${\bf T}^m$-action on the symplectic
manifold $(M_P,\om_P)$ is Hamiltonian with moment map $\mu_P : M_P
\to P$ fibering $M_P$ as a singular torus bundle over the polytope
$P$ (see e.g., \cite[\S 4.2]{F}). The volume form on $M_P$ can be
written in terms of the moment map: $\frac 1{m!} \om_P^m= dI
\wedge d\theta$, where $dI=\mu^*(dx_1\wedge\cdots dx_m)$.
Integrating, we obtain
\begin{equation}\label{vol-deg}\vol(P)=\int_{M_P}dI\wedge d\theta=
\int_{M_P}\frac{1}{m!}\om_P^m= \frac{1}{m!}\deg
[c_1(L_P)^m]\;,\end{equation}  which yields Kouchnirenko's Theorem, since
$\deg [c_1(L_P)^m]$ equals the number of points in the
intersection of the divisors of $m$ generic sections in
$H^0(M_P,L_P)\simeq \poly(P)$.  (This proof  of Kouchnirenko's Theorem  is
from Atiyah \cite{A}.)

\begin{rem}
Using the isomorphism
(\ref{toric}),  one can view $\chi_P(e^{i\theta})$ as the character of the torus ${\bf
T}^m$ action on the space $H^0(M_P, L_P)$ of holomorphic sections
of the line bundle $L_P$ over the toric variety $M_P$
(or equivalently as the equivariant index of the $\bar{\partial}$
operator on sections of $L_P$).  The  characters $\{\chi_{N P} (e^{i\theta})\}$ are
the characters of the powers $L_P^N$ of $L_P$ and are given by the
equivariant Riemann-Roch formula \cite{BV3, Gu2} (see also \cite[\S8]{BP}).  When
$\theta=0$, the Riemann-Roch formula gives another form of the Ehrhart formula for the
number of lattice points of $NP$:
\begin{eqnarray} \chi_{NP}(1)\ =\ \# (NP) &=&\dim \poly(NP)\ =\
\dim H^0(M_P,L_P^N)\nonumber \\&=&\chi(M_P,L_P^N)\ =\ \sum_{k=0}^m
\deg\big[c_1(L_P)^k\cup \mbox{Todd}_{m-k}(M_P)\big]
\frac{N^k}{k!}\nonumber \\&=&
 \frac{\deg [c_1(L_P)^m]}{m!} N^n +\dots + \deg
\mbox{Todd}_m(M_P)\,.\label{realRR}\end{eqnarray} (Note that equations
(\ref{RR}) and (\ref{realRR}) provide an alternate derivation of
(\ref{vol-deg}).)\end{rem}

\subsection{Geometric proof of Theorem \ref{SZEGO}.}

We outline here an alternative approach to evaluating the integral
(\ref{lattice3}) by lifting it to $M_P$, obtaining a complex
oscillatory integral over $\T\times M_P$. By (\ref{charN}), we
have
\begin{equation}\label{lattice4}\ical_N^1=
\frac{N^m}{(2\pi )^m} \int_{\T }d\phi \cdot \xi_1(\phi)
e^{N\Psi_0(\phi;\tau,z)} \todd(\fcal, N^{-1}
\partial/\partial h)  \left.\left(\int_{P(h)} e^{ N \langle \tau +i\phi, x
\rangle}\, dx\right)\right|_{h = 0}.\end{equation} We can analyze
this integral by lifting  to the toric variety $(M_h, \omega_h)$
associated to the deformed polytope $P(h)$ given by (\ref{P(h)}).
 The underlying complex variety $M_h = M_P$ is fixed,
while  the symplectic form is given by
\begin{equation*} \omega_h = \omega_P + \sum_{k = 1}^d h_k c_1(L_k)\;, \end{equation*}
which is affine in $h$.  Here, $L_k$ is the line bundle associated
to the divisor $\mu_P\inv(\bar F_k)$ in $M_P$, and $c_1(L_k)$ is a
$\T$-invariant $(1,1)$-form in the Chern class of $L_k$.

The moment map $\mu_h=(I_1^h,\dots,I_m^h):M_h\to P(h)$ is also an
affine function in $h$. Indeed, let $X_j$ denote the vector field
on $M_P$ giving the infinitesimal action on $\T$ corresponding to
the $j$-th unit vector in the Lie algebra $ {\go t}^m\approx
\R^m$. Then we have:
\begin{equation*}  d I_j^h  = d\mu_h(X_j) = \omega_h(X_j, \cdot)=
\omega_P(X_j, \cdot) + \sum_{k = 1}^d h_k c_1(L_k)(X_j, \cdot) =
dI_j + \sum_{k = 1}^d h_k c_1(L_k)(X_j, \cdot)\end{equation*} and
hence $$I_j^h = I_j + \sum_{k = 1}^d h_k J^k_j\;.$$

Using the fact that
$$\int_{P(h)}f(x)\,dx = \int_{M_P}f\circ \mu_h \frac 1{m!} \om_h^m\;,$$
 we then lift
(\ref{lattice4}) to $M_h$ to obtain:
\begin{equation} \label{lattice5} \ical_N^1=
\frac{N^m}{(2\pi )^m m!} \;  \todd(\fcal, N^{-1} \partial/\partial
h)|_{h = 0} \int_{{\bf T}^m} \int_{M_P} e^{N
\Psi_h(\phi,w;\tau,z)}\xi_1(\phi) \,\omega^m_h(w)\,
d\phi\;.
\end{equation} The phase is now the function on $\T \times M_P$ given by
\begin{equation}\Psi_h(\phi,w;\tau,z)= \langle \tau+ i\phi,\mu_P(w) \rangle
 +\sum_{k=1}^d h_k\langle \tau+i\phi,J^k(w)\rangle + p\log\left(\frac{ 1+\sum_{j=1}^m
e^{-\tau_j-i\phi_j}|z_j|^2}{1+\|z\|^2}\right)\;.\end{equation} The
amplitude is a polynomial in $h$:
\begin{equation} \om^m_h= \om_P^m + \sum_{1\le |\al|\le m}h^\al \ga_\al=G(w,h)
\om_P^m\;.
\end{equation}

Interchanging the order of integration and Todd differentiation in
(\ref{lattice5}) and then using (\ref{switch}) as before, we
obtain:
\begin{equation} \label{lattice6} \ical_N^1=
\frac{N^m}{(2\pi )^m m!} \;  \int_{{\bf T}^m} \int_{M_P} e^{N
\Psi(\phi,\mu_P(w);\tau,z)}A(N,\phi,w;\tau)
\,\omega_P^m(w)\, d\phi\;,
\end{equation}
with amplitude
$$A(N,\phi,w;\tau)= \xi_1(\phi)G(w,N\inv \d/\d q)\todd(\fcal,q)|_{q_k=\langle \tau+ i
\phi,J^k(w)\rangle}.$$ The phase is our familiar phase function
$\Psi$ given by (\ref{tauphase}) lifted to $M_P$.

We need to find the critical points where the phase $\Psi$ has
maximal real part. As we determined before, the phase has maximal
real part where $\phi=0$ and $\tau$ is in the normal cone to $P$
at $\mu_P(w)$. As before, we find that the critical point equation
$d_\phi\Psi(\phi,\mu_P(w);\tau,z)|_{\phi=0}=0$ is equivalent to
\begin{equation*} \frac{1}{p}\,
\mu_P(w)=\mu(e^{-\tau/2}\cdot z)\;.\end{equation*} The second
critical point equation $d_w\Psi(0,\mu_P(w);\tau,z)=0$
reduces to
\begin{equation*}   d_w  \langle\mu_P(w),\tau\rangle = 0. \end{equation*}
Since $\frac{\d}{\d w_j}\mu_P|_{w^0}$ is always tangent to
the face of $P$ at $\mu_P(w^0)$, this is just the condition
that $\tau$ is orthogonal to the face of $P$ containing
$\mu_P(w)$, which is automatically satisfied when $\tau$ is in
the normal cone to $P$ at $\mu_P(w)$.  Hence, in order to have
critical points, $\tau$ must equal $\tau_z$ as before, and we
obtain a critical submanifold
$$\ccal_z:=\{(0,w): \mu_P(w) =q(z)\}\;.$$  Note that $\dim \ccal_z = \codim F$,
where $F$ is the face of $P$ containing $q(z)$.

By a calculation using the methods of \cite[\S 3.2.2]{SZ2} (where
a more complicated phase function was used), one can show that the
normal Hessian is nondegenerate whenever $z$ is not a transition
point. The asymptotic expansion of Theorem~\ref{SZEGO} then
follows by the method of stationary phase with nondegenerate
critical submanifolds.

\begin{rem} We may also express polytope characters through the \szego kernel
$\Pi^{M_P}$  of $H^0(M_P, L_P^N)$  by means of the obvious
identity
\begin{equation} \label{SZCH} \chi_{NP}(e^{i \phi}) = \int_{M_P} \Pi_{N
}^{M_P}(e^{i\phi} \cdot w, w)\,  d\vol_{M_P}(w).
\end{equation}
In an article with T. Tate \cite{STZ1}, we describe an explicit
construction of the \szego kernel of a toric variety, which we use to obtain a
 formula for the polytope character  $\chi_{NP}(e^{i \phi})$ as a complex oscillatory
integral over the toric variety $M_P$.
\end{rem}

\twocolumn
\section*{Index of notation}

\begin{tabbing}
$\todd(\fcal, \partial/\partial h)$ \quad  \=\kill \\
${p\choose\al}$\> \eqref{multinom}\\
$\langle\cdot, \cdot\rangle$\> \eqref{IP}\\
$\ga_p$\> \eqref{G}\\
$\ga_{p|P}$\> \eqref{CG}\\
$\kappa_N$\> \eqref{w1}\\
$\mu$\> \eqref{momap}\\
$\Pi_p$\> \eqref{szego-proj}, \eqref{Sz}\\
$\Pi_{|P}$\> \eqref{Sz2}, \eqref{Sz2e}\\
$\Si$\> \S \ref{s-results}\\
$\tau_z$\> \eqref{cond1}--\eqref{cond2}\\
$\chi_\al$\> \eqref{chi}\\
$\|\chi_\al\|$\> \eqref{IP2}\\
$\wh\chi_\al^p$\> \eqref{chihat}, \eqref{mchi}\\
$\chi_{NP}$\> \eqref{char}\\
$\psi_P$\> \eqref{psiP}\\
$\Psi$\> \eqref{tauphase}\\
$\Psi_\acal$\> \eqref{phase}\\
$\om_\FS$\> \eqref{FSkahler}\\
$\acal_P$\> \eqref{AP}\\
$b_P$\> \eqref{b}\\
$b_{\{x\}}$\> \eqref{bx1}, \eqref{bx2}\\
$C_F$\> \eqref{ncone}\\
$C_x=C_x^P$\> \eqref{Cx}\\
$E_{|NP}$\> \eqref{EZ0}\\
$E_{|P_1,\dots,P_k}$\> \S \ref{s-distz}\\
$\ecal_Q$\> \eqref{EQ}\\
Flow$(x)$\> \eqref{flow}\\
$|f(z)|_\FS$\> \eqref{FSnorm}\\
$\ical_N^1,\ \ical_N^2$\> \eqref{lattice3}\\
$\jcal(x)$\> \eqref{J}\\
$\ell_j$\> \eqref{Pdef}\\
$\lcal$\> \eqref{lcal}\\
Log\> \eqref{LOG}\\
$\mcal_x$\> \eqref{mcal}\\
$\ncal(Q)$\> \S \ref{s-normalbundle}\\
$P_f$\> \eqref{NEWTONP}\\
$P(h)$\> \eqref{P(h)}\\
$\poly(P)$\> \eqref{polyP}\\
$q(z)$\> \eqref{cond1}--\eqref{cond2}\\
$\rcal_F$\> \eqref{RF}\\
$S(\ep)$\> \eqref{strip}\\
$S_f$\> \eqref{SUPPORT}\\
$\T$\> \eqref{torus}\\
$\todd(\fcal, \partial/\partial h)$\> \eqref{toddf}\\
$u_N$\> \eqref{uN}\\
$u_\infty$\> \eqref{u0}\\
$Z_{f_1,\dots,f_k}$\> \eqref{zerocurr}\\
$|Z_{f_1,\dots,f_k}|$\> \eqref{zeroset}\\
$\|Z_{f_1,\dots,f_k}\|$\> \eqref{zerovol}\\

\end{tabbing}

\onecolumn

\end{document}